\documentclass[pra,nofootinbib]{revtex4}

\usepackage{amsmath,amssymb}
\usepackage{mathrsfs}
\usepackage{mathbbol}
\usepackage{pst-func}
\usepackage{epsfig}
\usepackage{boldline}
\usepackage{cancel}
\usepackage[math]{cellspace}
\usepackage{listings}
\usepackage{tabularray}
\usepackage{enumitem}
\usepackage{siunitx}
\usepackage{chemarrow}

\allowdisplaybreaks

\begin{document}

\title{Eight-dimensional Octonion-like but Associative Normed Division Algebra}

\author{Joy Christian}

\email{jjc@bu.edu}

\affiliation{Einstein Centre for Local-Realistic Physics, Oxford OX2 6LB, United Kingdom}
\maketitle
\begin{center}
ABSTRACT:
\end{center}
\begingroup
\addtolength\leftmargini{0.24in}
\begin{quote}
We present an eight-dimensional even sub-algebra of the ${2^4=16}$-dimensional associative Clifford algebra ${\mathrm{Cl}_{4,0}}$ and show that its eight-dimensional multivectors ${\bf X}$ and ${\bf Y}$ respect the composition law ${||{\bf X}{\bf Y}||=||{\bf X}||\,||{\bf Y}||}$, thus forming an octonion-like but associative normed division algebra, where the norms are calculated using the fundamental geometric product instead of the usual scalar product so that the underlying coefficient algebra resembles split complex numbers instead of reals. The corresponding 7-sphere obtained from projecting this multivector-valued composition law to the scalar-valued composition law has a topology that differs from that of the octonionic 7-sphere. Just as the octonionic 7-sphere is parallelizable using the non-associative algebra of octonions, we demonstrate that the 7-sphere presented herein is parallelizable using the said associative algebra.   
\end{quote}
\endgroup

\vspace{0.2in}
\parindent 0pt
[\underbar{\bf Note added to proof}$\,$: The results of this paper are published in Section~2.8 of Ref.~\cite{RSOS-2} listed in the bibliography. \\[0.1cm]
In several publications between 2018 and 2022, I proposed this eight-dimensional associative normed division algebra in \cite{RSOS,RSOS-2,Local}, in analogy with the well-known non-associative octonionic normed division algebra. Unfortunately, some misunderstandings regarding this straightforward demonstration persist in the literature. Therefore, in Appendix~\ref{D} below, starting on page 12, I have added a streamlined proof of the algebra, with emphasis on pedagogical explanation.]
\vspace{0.2in}

\parskip 5pt
\parindent 15pt
\baselineskip 12pt

Consider the following eight-dimensional vector space with graded Clifford-algebraic basis and orientation ${\lambda=\pm1}$:
\begin{equation}
Cl^{\lambda}_{3,0}={\rm span}\!\left\{\,1,\;\lambda{\bf e}_x,\,\lambda{\bf e}_y,\,\lambda{\bf e}_z,\;\lambda{\bf e}_x{\bf e}_y,\,
\lambda{\bf e}_z{\bf e}_x,\,\lambda{\bf e}_y{\bf e}_z,\;\lambda{\bf e}_x{\bf e}_y{\bf e}_z\,\right\}\!. \label{cl30}
\end{equation}
As we shall see, the choice of orientation, ${\lambda=+1}$ or ${\lambda=-1}$ does not affect our argument. In what follows, we will use the language of Geometric Algebra, as used, for example, in Refs.~\cite{Clifford} and \cite{Dorst}. Accordingly, in the above definiton ${\{\,{\bf e}_x,\,{\bf e}_y,\,{\bf e}_z\}}$ is a set of anti-commuting orthonormal vectors in ${{\rm I\!R}^3}$ such that ${{\bf e}_j{\bf e}_i=-\,{\bf e}_i{\bf e}_j}$ for any ${i\not=j=x,\,y}$, or ${z}$. In general the vectors ${{\bf e}_i}$ satisfy the following geometric product in this {\it associative} but non-commutative algebra~\cite{Clifford,Dorst}:
\begin{equation}
{\bf e}_i\,{\bf e}_j\,=\,{\bf e}_i\cdot{\bf e}_j+\,{\bf e}_i\wedge\,{\bf e}_j\,, \label{gp}
\end{equation}
with
\begin{equation}
{\bf e}_i\cdot{\bf e}_j:=\,\frac{1}{2}\left\{{\bf e}_i{\bf e}_j+{\bf e}_j{\bf e}_i\right\} \label{sim-gp}
\end{equation}
being the symmetric inner product and
\begin{equation}
{\bf e}_i\wedge{\bf e}_j:=\,\frac{1}{2}\left\{{\bf e}_i{\bf e}_j-{\bf e}_j{\bf e}_i\right\}\, \label{ant-gp}
\end{equation}
being the anti-symmetric outer product, giving ${({\bf e}_i\wedge{\bf e}_j)^2=-1}$. There are thus basis elements of four different grades in this algebra: An identity element ${{\bf e}^2_i=1}$ of grade-0, three orthonormal vectors ${{\bf e}_i}$ of grade-1, three orthonormal bivectors ${{\bf e}_j{\bf e}_k}$ of grade-2, and a trivector ${I_3={\bf e}_i{\bf e}_j{\bf e}_k}$ of grade-3 representing a volume element in ${{\rm I\!R}^3}$. Since in ${{\rm I\!R}^3}$ there are ${2^3=8}$ ways to combine the vectors ${{\bf e}_i}$ using the geometric product (\ref{gp}) such that no two products are linearly dependent, the resulting algebra, ${Cl^{\lambda}_{3,0}\,}$, is a linear vector space of eight dimensions, spanned by these graded bases.

In this paper, we are interested in a certain conformal completion\footnote{\label{Conformal}The conformal space we are considering is an {\it in}-homogeneous version of the space usually studied in Conformal Geometric Algebra \cite{Dorst}. It can be viewed as an ${8}$-dimensional subspace of the ${32}$-dimensional representation space postulated in Conformal Geometric Algebra. The larger representation space results from a homogeneous freedom of the origin within ${{\mathbb E}^3}$, which does not concern us in this paper.} of this algebra, originally presented in Ref.~\cite{RSOS}. This is accomplished by using an additional vector, ${{\bf e}_{\infty}}$, to close the lines and volumes of the Euclidean space, giving
\begin{equation}
{\cal K}^{\lambda}=\,{\rm span}\!\left\{\,1,\,\lambda{\bf e}_x{\bf e}_y,\,\lambda{\bf e}_z{\bf e}_x,\,\lambda{\bf e}_y{\bf e}_z,\,\lambda{\bf e}_x{\bf e}_{\infty},\,\lambda{\bf e}_y{\bf e}_{\infty},\,\lambda{\bf e}_z{\bf e}_{\infty},\,\lambda I_3{\bf e}_{\infty}\,\right\}\!. \label{clind}
\end{equation}
With unit vector ${{\bf e}_{\infty}}$, this is an eight-dimensional {\it even} sub-algebra of the ${2^4=16}$-dimensional Clifford algebra ${Cl_{4,0}}$.
\begin{table}[t]
\hrule
\vspace{18pt}
\begin{center}
\SetTblrInner{rowsep=10pt}
\begin{tblr}[b]{|c|[2pt]c|c|c|c|c|c|c|c|} \hline
${*}$ &${1}$ &${\lambda\,{\bf e}_x{\bf e}_y}$ &${\lambda\,{\bf e}_z{\bf e}_x}$ &${\lambda\,{\bf e}_y{\bf e}_z}$ &${\lambda\,{\bf e}_x{\bf e}_{\infty}}$ &${\lambda\,{\bf e}_y{\bf e}_{\infty}}$ &${\lambda\,{\bf e}_z{\bf e}_{\infty}}$ &${\lambda\,I_3{\bf e}_{\infty}}$ \\ \hline[2pt]
${1}$ &${1}$ &${\lambda\,{\bf e}_x{\bf e}_y}$ &${\lambda\,{\bf e}_z{\bf e}_x}$ &${\lambda\,{\bf e}_y{\bf e}_z}$ &${\lambda\,{\bf e}_x{\bf e}_{\infty}}$ &${\lambda\,{\bf e}_y{\bf e}_{\infty}}$ &${\lambda\,{\bf e}_z{\bf e}_{\infty}}$ &${\lambda\,I_3{\bf e}_{\infty}}$ \\ \hline
${\lambda\,{\bf e}_x{\bf e}_y}$ &${\lambda\,{\bf e}_x{\bf e}_y}$ &${-1}$ &${{\bf e}_y{\bf e}_z}$ &${-{\bf e}_z{\bf e}_x}$ &${-{\bf e}_y{\bf e}_{\infty}}$ &${{\bf e}_x{\bf e}_{\infty}}$ &${I_3{\bf e}_{\infty}}$ &${-{\bf e}_z{\bf e}_{\infty}}$ \\ \hline
${\lambda\,{\bf e}_z{\bf e}_x}$ &${\lambda\,{\bf e}_z{\bf e}_x}$ &${-{\bf e}_y{\bf e}_z}$ &${-1}$ &${{\bf e}_x{\bf e}_y}$ &${{\bf e}_z{\bf e}_{\infty}}$ &${I_3{\bf e}_{\infty}}$ &${-{\bf e}_x{\bf e}_{\infty}}$ &${-{\bf e}_y{\bf e}_{\infty}}$ \\ \hline
${\lambda\,{\bf e}_y{\bf e}_z}$ &${\lambda\,{\bf e}_y{\bf e}_z}$ &${{\bf e}_z{\bf e}_x}$ &${-{\bf e}_x{\bf e}_y}$ &${-1}$ &${I_3{\bf e}_{\infty}}$ &${-{\bf e}_z{\bf e}_{\infty}}$ &${{\bf e}_y{\bf e}_{\infty}}$ &${-{\bf e}_x{\bf e}_{\infty}}$ \\ \hline
${\lambda\,{\bf e}_x{\bf e}_{\infty}}$ &${\lambda\,{\bf e}_x{\bf e}_{\infty}}$ &${{\bf e}_y{\bf e}_{\infty}}$ &${-{\bf e}_z{\bf e}_{\infty}}$ &${I_3{\bf e}_{\infty}}$ &${-1}$ &${-{\bf e}_x{\bf e}_y}$ &${{\bf e}_z{\bf e}_x}$ &${-{\bf e}_y{\bf e}_z}$ \\ \hline
${\lambda\,{\bf e}_y{\bf e}_{\infty}}$ &${\lambda\,{\bf e}_y{\bf e}_{\infty}}$ &${-{\bf e}_x{\bf e}_{\infty}}$ &${I_3{\bf e}_{\infty}}$ &${{\bf e}_z{\bf e}_{\infty}}$ &${{\bf e}_x{\bf e}_y}$ &${-1}$ &${-{\bf e}_y{\bf e}_z}$ &${-{\bf e}_z{\bf e}_x}$ \\ \hline
${\lambda\,{\bf e}_z{\bf e}_{\infty}}$ &${\lambda\,{\bf e}_z{\bf e}_{\infty}}$ &${I_3{\bf e}_{\infty}}$ &${{\bf e}_x{\bf e}_{\infty}}$ &${-{\bf e}_y{\bf e}_{\infty}}$ &${-{\bf e}_z{\bf e}_x}$ &${{\bf e}_y{\bf e}_z}$ &${-1}$ &${-{\bf e}_x{\bf e}_y}$ \\ \hline
${\lambda\,I_3{\bf e}_{\infty}}$ &${\lambda\,I_3{\bf e}_{\infty}}$ &${-{\bf e}_z{\bf e}_{\infty}}$ &${-{\bf e}_y{\bf e}_{\infty}}$ &${-{\bf e}_x{\bf e}_{\infty}}$ &${-{\bf e}_y{\bf e}_z}$ &${-{\bf e}_z{\bf e}_x}$ &${-{\bf e}_x{\bf e}_y}$ &${1}$ \\
\hline
\end{tblr}
\end{center}
\vspace{5pt}
\hrule
\caption{Multiplication Table for a ``Conformal Geometric Algebra\textsuperscript{\ref{Conformal}}" of ${{\mathbb E}^3}$. Here ${I_3={\bf e}_x{\bf e}_y{\bf e}_z}$, ${{\bf e}_{\infty}^2=+1}$, and ${\lambda =\pm 1}$.\break}
\vspace{5pt}
\hrule
\label{T+1}
\end{table}
Unlike the seven imaginaries of octonions, there are only six basis elements of ${\cal K}^{\lambda}$ that are imaginary. The seventh, ${\lambda I_3{\bf e}_{\infty}}$, squares to +1. This is evident from the multiplication table \ref{T+1}. We therefore call it an ``octonian-like" algebra. As an eight-dimensional linear vector space, ${{\cal K}^{\lambda}}$ has some remarkable properties \cite{RSOS}. To begin with, it is {\it closed} under multiplication. Suppose ${\bf X}$ and ${\bf Y}$ are two vectors in ${{\cal K}^{\lambda}}$. Then ${\bf X}$ and ${\bf Y}$ can be expanded in the graded basis of ${{\cal K}^{\lambda}}$:
\begin{equation}
{\bf X}=\,X_0+X_1\,\lambda{\bf e}_x{\bf e}_y+X_2\,\lambda{\bf e}_z{\bf e}_x+X_3\,\lambda{\bf e}_y{\bf e}_z+X_4\,\lambda{\bf e}_x{\bf e}_{\infty}+X_5\,\lambda{\bf e}_y{\bf e}_{\infty}+X_6\,\lambda{\bf e}_z{\bf e}_{\infty}+X_7\,\lambda I_3{\bf e}_{\infty} \label{X}
\end{equation}
and
\begin{equation}
{\bf Y}=\,Y_0+Y_1\,\lambda{\bf e}_x{\bf e}_y+Y_2\,\lambda{\bf e}_z{\bf e}_x+Y_3\,\lambda{\bf e}_y{\bf e}_z+Y_4\,\lambda{\bf e}_x{\bf e}_{\infty}+Y_5\,\lambda{\bf e}_y{\bf e}_{\infty}+Y_6\,\lambda{\bf e}_z{\bf e}_{\infty}+Y_7\,\lambda I_3{\bf e}_{\infty}\,.\label{Y}
\end{equation}
And using the definition ${||{\bf X}||^2:={\bf X}\cdot{\bf X}^{\dagger}}$ for the quadratic form ${{\cal Q}({\bf X})}$ (where ${\scriptstyle\dagger}$ represents the {\it reverse} operation \cite{Clifford} and ${{\bf X}\cdot{\bf X}^{\dagger}}$ represents the scalar part of the geometric product ${{\bf X}{\bf X}^{\dagger}}$) the multivectors ${\bf X}$ and ${\bf Y}$ can be normalized as
\begin{equation}
||{\bf X}||^2=\sum_{\mu\,=\,0}^{7} \;X_{\mu}^2\,=\,1\;\;\;\text{and}\;\;\;||{\bf Y}||^2=\sum_{\nu\,=\,0}^{7} \;Y_{\nu}^2\,=\,1\,.\label{now8}
\end{equation}
Now it is evident from the multiplication table above (Table \ref{T+1}) that if ${{\bf X},{\bf Y}\in{\cal K}^{\lambda}}$, then so is their product ${{\bf Z}={\bf X}{\bf Y}}$:
\begin{equation}
{\bf Z}=\,Z_0+Z_1\,\lambda{\bf e}_x{\bf e}_y+Z_2\,\lambda{\bf e}_z{\bf e}_x+Z_3\,\lambda{\bf e}_y{\bf e}_z+Z_4\,\lambda{\bf e}_x{\bf e}_{\infty}+Z_5\,\lambda{\bf e}_y{\bf e}_{\infty}+Z_6\,\lambda{\bf e}_z{\bf e}_{\infty}+Z_7\,\lambda I_3{\bf e}_{\infty}={\bf X}{\bf Y}. \label{Z}
\end{equation}
Thus ${{\cal K}^{\lambda}}$ remains closed under arbitrary number of multiplications of its elements. This is a powerful property. More importantly, we shall soon see that for vectors ${\bf X}$ and ${\bf Y}$ in ${{\cal K}^{\lambda}}$ (not necessarily unit) the following norm relation holds:
\begin{equation}
||{\bf X}{\bf Y}|| \,=\, ||{\bf X}||\;||{\bf Y}||\,, \label{norms}
\end{equation}
provided the norms are calculated employing the fundamental geometric product instead of the usual scalar product. In particular, this means that for any two unit vectors ${\bf X}$ and ${\bf Y}$ in ${{\cal K}^{\lambda}}$ with the geometric product ${{\bf Z}={\bf X}{\bf Y}}$ we have
\begin{equation}
||\,{\bf Z}\,||^2=\sum_{\rho\,=\,0}^{7} \;Z_{\rho}^2\,=\,1\,.
\end{equation}
\begin{figure}
\hrule
\scalebox{0.9}{
\begin{pspicture}(0.5,-3.5)(4.2,4.7)

\psline[linewidth=0.3mm,arrowinset=0.3,arrowsize=3pt 3,arrowlength=2]{->}(-1.27,-2.3)(6.0,-2.3)

\psline[linewidth=0.3mm,arrowinset=0.3,arrowsize=3pt 3,arrowlength=2]{->}(-0.77,-2.8)(-0.77,3.3)

\psline[linewidth=0.2mm,arrowinset=0.3,arrowsize=3pt 3,arrowlength=2]{->}(-1.27,-2.3)(4.0,-2.3)

\psline[linewidth=0.1mm,dotsize=2pt 3]{*-}(3.99,-2.3)(4.01,-2.3)

\psline[linewidth=0.2mm,arrowinset=0.3,arrowsize=3pt 3,arrowlength=2]{->}(-0.77,-2.8)(-0.77,1.53)

\psline[linewidth=0.1mm,dotsize=2pt 3]{*-}(-0.77,1.51)(-0.77,1.52)

\psline[linewidth=0.3mm,arrowinset=0.3,arrowsize=3pt 3,arrowlength=2]{->}(-0.77,-2.3)(4.0,1.5)

\psarc[linewidth=0.3mm,arrowinset=0.3,arrowsize=3pt 3,arrowlength=2]{<->}(-0.77,-2.3){6.1}{13}{65}

\put(2.4,3.15){\large ${S^7}$}

\psarc[linewidth=0.3mm,arrowinset=0.3,arrowsize=3pt 3,arrowlength=2]{<->}(-0.77,-2.3){3.8}{50}{130}

\put(-3.27,1.15){\large ${S^3}$}

\psline[linewidth=0.1mm,dotsize=2pt 3]{*-}(4.0,1.5)(4.0,1.55)

\psline[linewidth=0.1mm,dotsize=2pt 3]{*-}(-0.77,-2.31)(-0.77,-2.27)

\psline[linewidth=0.2mm,linestyle=dashed,arrowinset=0.3,arrowsize=3pt 3,arrowlength=2]{-}(-0.77,1.5)(4.0,1.5)

\psline[linewidth=0.2mm,linestyle=dashed,arrowinset=0.3,arrowsize=3pt 3,arrowlength=2]{-}(4.0,-2.3)(4.0,1.5)

\put(6.05,3.55){\large ${{\cal K}^{\lambda}}$}

\put(4.2,1.7){\large ${{\mathbb Q}_z =\, {\bf q}_r + {\bf q}_d\,\varepsilon}$}

\put(3.75,-2.80){\large ${{\bf q}_d\,\varepsilon}$}

\put(-1.32,1.75){\large ${{\bf q}_r}$}

\put(6.35,-2.45){\large ${{\rm I\!R}^4}$}

\put(-0.99,3.7){\large ${{\rm I\!R}^4}$}

\put(-1.23,-2.85){\large ${0}$}

\end{pspicture}}
\hrule
\caption{An illustration of the 8D plane of ${{\cal K}^{\lambda}}$, which may be interpreted as an Argand diagram for a pair of quaternions.}
\vspace{5pt}
\hrule
\label{fig1}
\end{figure}

Now, in order to prove the norm relation (\ref{norms}), it is convenient to express the elements of ${{\cal K}^{\lambda}}$ as ``dual" quaternions. The idea of dual numbers, ${z}$, analogous to complex numbers, was introduced by Clifford in his seminal work as follows:
\begin{equation}
z = r + d\,\varepsilon, \;\;\text{where}\;\,\varepsilon\not=0 \;\,\text{but}\;\, \varepsilon^2=0\,.
\end{equation}
Here ${\varepsilon}$ is the dual operator, ${r}$ is the real part, and ${d}$ is the dual part \cite{Dorst,Kenwright}. Similar to how the ``imaginary" operator ${i}$ is introduced in the complex number theory to distinguish the ``real" and ``imaginary" parts of a complex number, Clifford introduced the dual operator ${\varepsilon}$ to distinguish the ``real" and ``dual" parts of a dual number. The dual number theory can be extended to numbers of higher grades, including to numbers of composite grades, such as quaternions.

In analogy with dual numbers, but with ${\varepsilon^2=+1}$, it is convenient for our purposes to write the elements of ${{\cal K}^{\lambda}}$ as
\begin{equation}
{\mathbb Q}_z =\, {\bf q}_r + {\bf q}_d\,\varepsilon\,, \label{appen}
\end{equation}
where ${{\bf q}_r}$ and ${{\bf q}_d}$ are quaternions and ${{\mathbb Q}_z}$ may now be viewed as a ``dual"-quaternion (or in Clifford's terminology, as a bi-quaternion). Next, recall that the set of unit quaternions is a 3-sphere, which can be normalized to a radius ${\varrho_r}$ and written as the set
\begin{equation}
S^3=\left\{\,{\bf q}_r:=\,q_0 + q_1\,\lambda\,{\bf e}_x{\bf e}_y+q_2\,\lambda\,{\bf e}_z{\bf e}_x+q_3\,\lambda\,{\bf e}_y{\bf e}_z \;\Big|\;||{\bf q}_r|| = \sqrt{{\bf q}_{r}\,{\bf q}^{\dagger}_{r}\,} = \varrho_r \, \right\}. \label{s-r}
\end{equation}
Consider now a second, dual copy of the set of quaternions within ${{\cal K}^{\lambda}}$, corresponding to the fixed orientation ${\lambda=+1}$:
\begin{equation}
S^3=\left\{\,{\bf q}_d:=-q_7 + q_6\,{\bf e}_x{\bf e}_y+q_5\,{\bf e}_z{\bf e}_x+q_4\,{\bf e}_y{\bf e}_z \;\Big|\;||{\bf q}_d|| = \sqrt{{\bf q}_{d}\,{\bf q}^{\dagger}_{d}\,} = \varrho_d \, \right\}. \label{d-d}
\end{equation}
If we now identify ${\lambda\,I_3{\bf e}_{\infty}}$ appearing in (\ref{clind}) as the duality operator ${-\,\varepsilon}$, then (in the reverse additive order) we obtain
\begin{align}
&\,\varepsilon \equiv\,-\lambda\,I_3{\bf e}_{\infty},\;\;\text{with}\;\;\varepsilon^{\dagger}=\varepsilon\;\;\text{and}\;\;\varepsilon^2=+1\;\;(\text{since ${{\bf e}_{\infty}}$ is a unit vector in ${{\cal K}^{\lambda}}$}), \label{16-A} \\
&\text{and}\;\;{\bf q}_d\,\varepsilon \equiv\, -{\bf q}_d\,\lambda\,I_3{\bf e}_{\infty}=\,q_4\,\lambda\,{\bf e}_x{\bf e}_{\infty}+\,q_5\,\lambda\,{\bf e}_y{\bf e}_{\infty}+\,q_6\,\lambda\,{\bf e}_z{\bf e}_{\infty}+\,q_7\,\lambda\,I_3{\bf e}_{\infty}\,, \label{s-d}
\end{align}
which is a multi-vector ``dual" to the quaternion ${{\bf q}_d}$. Note that we write ${\varepsilon}$ as if it were a scalar because it commutes with all elements of ${{\cal K}^{\lambda}}$ in (\ref{clind}). Comparing (\ref{s-r}) and (\ref{s-d}) with (\ref{clind}) we can now write ${{\cal K}^{\lambda}}$ as a set of paired quaternions,
\begin{equation}
{\cal K}^{\lambda}=\,\left\{\,{\mathbb Q}_z :=\, {\bf q}_r + {\bf q}_d\,\varepsilon\;\Big|\; ||{\mathbb Q}_z|| = \sqrt{{\mathbb Q}_{z}\,{\mathbb Q}^{\dagger}_{z}\,} = \sqrt{\varrho^2_r + \varrho^2_d\,},\,0 < \varrho_r < \infty,\,0 < \varrho_d < \infty \,\right\}, \label{new18}
\end{equation}
in analogy with (\ref{s-r}) or (\ref{d-d}), with ${{\mathbb Q}_{z}\,{\mathbb Q}^{\dagger}_{z}}$ being the geometric product between ${{\mathbb Q}_{z}}$ and ${{\mathbb Q}^{\dagger}_{z}}$ (instead of the inner product ${{\mathbb Q}_{z}\cdot{\mathbb Q}^{\dagger}_{z}}$ used in (\ref{now8}) to calculate the value of ${||{\mathbb Q}_{z}||}$). But this definition ${||{\mathbb Q}_z|| = \sqrt{{\mathbb Q}_{z}\,{\mathbb Q}^{\dagger}_{z}\,} = \sqrt{\varrho^2_r + \varrho^2_d\,}}$ for the norm is possible only if we require ${{\bf q}_r\,{\bf q}^{\dagger}_d+{\bf q}_d\,{\bf q}^{\dagger}_r = 0}$, rendering every ${{\bf q}_r}$ orthogonal to its dual ${{\bf q}_d}$ (cf. Fig.~\ref{fig1}). In other words,
\begin{equation}
||{\mathbb Q}_z|| = \sqrt{{\mathbb Q}_{z}\,{\mathbb Q}^{\dagger}_{z}\,} = \sqrt{\varrho^2_r + \varrho^2_d\,}\;\;\Longleftrightarrow\;\;{\bf q}_r\,{\bf q}^{\dagger}_d+{\bf q}_d\,{\bf q}^{\dagger}_r = 0\,, \label{ken-norm}
\end{equation}
or equivalently, ${({\bf q}_r\,{\bf q}^{\dagger}_d)_s = 0}$; {\it i.e.}, ${{\bf q}_r\,{\bf q}^{\dagger}_d}$ must be a pure quaternion (for a pedagogical discussion of (\ref{ken-norm}) see section 7.1 of Ref.~\cite{Kenwright}). We can see this by working out the geometric product of ${{\mathbb Q}_{z}}$ with ${{\mathbb Q}^{\dagger}_{z}}$ while using ${\varepsilon^2=+1}$, which gives
\begin{equation}
{\mathbb Q}_{z}\,{\mathbb Q}^{\dagger}_{z}\,=\;\left({\bf q}_{r}\,{\bf q}^{\dagger}_{r}\,+\,{\bf q}_{d}\,{\bf q}^{\dagger}_{d}\right)\,+\,\left({\bf q}_{r}\,{\bf q}^{\dagger}_{d}\,+\,{\bf q}_{d}\,{\bf q}^{\dagger}_{r}\right)\,\varepsilon\,.
\end{equation}
Now, using definitions (\ref{s-r}) and (\ref{d-d}), it is easy to see that 
${{\bf q}_{r}\,{\bf q}^{\dagger}_{r}=\varrho^2_r}$ and  ${{\bf q}_{d}\,{\bf q}^{\dagger}_{d}=\varrho^2_d}$, reducing the above product to
\begin{equation}
{\mathbb Q}_{z}\,{\mathbb Q}^{\dagger}_{z}\,=\;
\varrho^2_r + \varrho^2_d\,+\,\left({\bf q}_{r}\,{\bf q}^{\dagger}_{d}\,+\,{\bf q}_{d}\,{\bf q}^{\dagger}_{r}\right)\,\varepsilon\,. \label{21}
\end{equation}
In terms of the coefficients of ${{\mathbb Q}_{z}}$ the quantity ${{\bf q}_{r}\,{\bf q}^{\dagger}_{d}+{\bf q}_{d}\,{\bf q}^{\dagger}_{r}}$ can be worked out and it turns out to be a scalar as well:
\begin{equation}
{\bf q}_{r}\,{\bf q}^{\dagger}_{d}+{\bf q}_{d}\,{\bf q}^{\dagger}_{r}\,=\,-2\,q_0q_7 + 2\lambda\,q_1q_6 + 2\lambda\,q_2q_5 + 2\lambda\,q_3q_4\,. \label{22}
\end{equation}
Consequently, since $\varepsilon$ appearing in (\ref{21}) is a pseudoscalar, the product ${\mathbb Q}_{z}\,{\mathbb Q}^{\dagger}_{z}$ between ${\mathbb Q}_{z}$ and ${\mathbb Q}^{\dagger}_{z}$ is always of the form
\begin{align}
{\mathbb Q}_{z}\,{\mathbb Q}^{\dagger}_{z}\,&=\,\text{(a scalar)} \,+\, \text{(a scalar)}\times\varepsilon\, \notag \\
&=\,\text{(a scalar)} \,+\, \text{(a pseudoscalar)}. \label{23} 
\end{align}
It is thus clear that for ${{\mathbb Q}_{z}\,{\mathbb Q}^{\dagger}_{z}}$ to be a scalar, ${{\bf q}_{r}\,{\bf q}^{\dagger}_{d}+{\bf q}_{d}\,{\bf q}^{\dagger}_{r}}$ must vanish, or equivalently ${{\bf q}_r}$ must be orthogonal to ${{\bf q}_d}$.

But there is more to the normalization condition ${{\bf q}_{r}\,{\bf q}^{\dagger}_{d}+{\bf q}_{d}\,{\bf q}^{\dagger}_{r}=0}$ than meets the eye. It also leads to the crucial norm relation (\ref{norms}), which is at the heart of the only known four normed division algebras ${\mathbb R}$, ${\mathbb C}$, ${\mathbb H}$ and ${\mathbb O}$ associated with the four parallelizable spheres ${S^0}$, ${S^1}$, ${S^3}$ and ${S^7\,}$, with octonions forming a non-associative algebra in addition to forming a non-commutative algebra \cite{Hurwitz-1898,Hurwitz-1923,Baez}. However, before we prove the norm relation (\ref{norms}), let us take a closer look at definition (\ref{ken-norm}) within Geometric Algebra. In (\ref{now8}) we used the following definition of norm for a general multivector:
\begin{equation}
||{\bf X}||=\sqrt{{\bf X}\cdot{\bf X}^{\dagger}}. \label{b}
\end{equation}
The left-hand side of this equation --- by definition --- is a scalar number; namely, $||{\bf X}||$. But what is important to recognize for our purpose of proving (\ref{norms}) is that there are two equivalent ways of working out this scalar number: 

(a) $\;\;\;\;||{\bf X}||$ = square-root of the scalar part ${{\bf X}\cdot{\bf X}^{\dagger}}$ of the geometric product ${{\bf X}{\bf X}^{\dagger}}$ between ${\bf X}$ and ${{\bf X}^{\dagger}}$,

or

(b) $\;\;\;\;||{\bf X}||$ = square-root of the geometric product ${{\bf X}{\bf X}^{\dagger}}$ with the non-scalar part of ${{\bf X}{\bf X}^{\dagger}}$ set to zero.

\pagebreak

\parindent 0pt

The above two definitions of the norm $||{\bf X}||$ are {\it entirely equivalent}. They give one and the same scalar value for the norm. Moreover, in general, given a product denoted by $*$, the quantity ${\bf X}^{\dagger}$ is said to be the conjugate of ${\bf X}$ if ${{\bf X}*{\bf X}^{\dagger}}$ happens to be equal to unity, ${{\bf X}*{\bf X}^{\dagger}}$ = 1, as in the case of quaternionic products in (\ref{s-r}) and (\ref{d-d}). On the other hand, in Geometric Algebra the fundamental product between any two multivectors ${\bf X}$ and ${\bf Y}$ is the geometric product, ${{\bf X}{\bf Y}}$, {\it not} the scalar product ${{\bf X}\cdot{\bf Y}}$ (or the wedge product ${{\bf X}\wedge{\bf Y}}$ for that matter). Therefore the product that must be used in computing the norm $||{\bf X}||$ that preserves the Clifford algebraic structure of ${\cal K}^{\lambda}$ is the geometric product ${{\bf X}{\bf X}^{\dagger}}$, not the scalar product ${{\bf X}\cdot{\bf X}^{\dagger}}$. To be sure, in practice, if one is interested only in working out the value of the norm $||{\bf X}||$, then it is often convenient to use the definition (a) above. However, our primary purpose here in working out the norms of ${\bf X}$ and ${\bf Y}$ is to preserve the algebraic structure of the space ${\cal K}^{\lambda}$ in the fundamental relation (\ref{norms}), and therefore the definition of the norm we must use is necessarily the second definition stated above; {\it i.e.}, definition (b).

\parindent 15pt

With the above comments in mind, we are now ready to prove the norm relation (\ref{norms}). To this end, suppose the multivectors ${\bf X}$ and ${\bf Y}$ belonging to ${\cal K}^{\lambda}$ as spelled out in (\ref{X}) and (\ref{Y}) are normalized using the definition (b) as follows:
\begin{align}
||{\bf X}||=\sqrt{{\bf X}{\bf X}^{\dagger}\;} &= \sqrt{(\varrho^2_{\!_{\bf X}} + 0\times\varepsilon)} = \varrho_{\!_{\bf X}}:=\sqrt{\varrho^2_{r1} + \varrho^2_{d1}}\\
\text{and}\;\;\;\;\;
||{\bf Y}||=\sqrt{{\bf Y}{\bf Y}^{\dagger}\;} &= \sqrt{(\varrho^2_{\!_{\bf Y}} + 0\times\varepsilon)} = \varrho_{\!_{\bf Y}}:=\sqrt{\varrho^2_{r2} + \varrho^2_{d2}}\,, 
\end{align}
where $\varrho_{\!_{\bf X}}$ and $\varrho_{\!_{\bf Y}}$ are fixed scalars. Then, according to (\ref{Z}), their product ${\bf X}{\bf Y}$ in ${\cal K}^{\lambda}$ is another multivector, giving
\begin{align}
||{\bf X}{\bf Y}||&=\sqrt{({\bf X}{\bf Y})({\bf X}{\bf Y})^{\dagger}\;} \\
&=\sqrt{{\bf X}{\bf Y}{\bf Y}^{\dagger}{\bf X}^{\dagger}\;} \\
&=\sqrt{{\bf X}\,\left(\varrho^2_{\!_{\bf Y}}+0\times\varepsilon\right)\,{\bf X}^{\dagger}\;} \\
&=\left(\sqrt{{\bf X}{\bf X}^{\dagger}\;}\,\right)\varrho_{\!_{\bf Y}} \\
&=\left(\sqrt{\left(\varrho^2_{\!_{\bf X}}+0\times\varepsilon\right)}\,\right)\,\varrho_{\!_{\bf Y}} \\
&=\varrho_{\!_{\bf X}}\,\varrho_{\!_{\bf Y}} \\
&=||{\bf X}||\,||{\bf Y}||.
\end{align}
Thus, at first sight, the norm relation (\ref{norms}) appears to be trivially true. However, the above simple proof is not quite satisfactory because we have assumed that ${\bf X}$ and ${\bf Y}$ are normalized using the definition (b), which requires us to set the non-scalar parts of the geometric products ${\bf X}{\bf X}^{\dagger}$ and ${\bf Y}{\bf Y}^{\dagger}$ equal to zero. That is not difficult to do for both ${\bf X}$ and ${\bf Y}$, but what is involved in the above proof is the geometric product ${{\bf X}{\bf Y}}$ and its conjugate ${({\bf X}{\bf Y})^{\dagger}}$, which makes the proof less convincing. It is therefore important to spell out the proof in full detail by assuming only that the non-scalar parts of the geometric products ${\bf X}{\bf X}^{\dagger}$ and ${\bf Y}{\bf Y}^{\dagger}$ are zero but {\it without} assuming {\it a priori} that the non-scalar part of the geometric product $({\bf X}{\bf Y})({\bf X}{\bf Y})^{\dagger}$ is zero; {\it i.e.}, we would like to derive the latter by assuming only the former. 

To that end, we first work out the right-hand side of the norm relation (\ref{norms}) in the notations of the condition (\ref{ken-norm}):
\begin{equation}
||{\mathbb Q}_{z1}||\;||{\mathbb Q}_{z2}||\,=\,
\left(\sqrt{\varrho^2_{r1}+\varrho^2_{d1}}\;\right)\left(\sqrt{\varrho^2_{r2}+\varrho^2_{d2}}\;\right)\,=\,\sqrt{\varrho^2_{r1}\,\varrho^2_{r2}+\varrho^2_{r1}\,\varrho^2_{d2}+\varrho^2_{d1}\,\varrho^2_{r2}+\varrho^2_{d1}\,\varrho^2_{d2}\,}\,. \label{s36}
\end{equation}
Now, to verify the left-hand side of the norm relation (\ref{norms}), consider a product of two distinct members of the set ${{\cal K}^{\lambda}}$,
\begin{equation}
{\mathbb Q}_{z1}\,{\mathbb Q}_{z2}\,=\;\left({\bf q}_{r1}\,{\bf q}_{r2}\,+\,{\bf q}_{d1}\,{\bf q}_{d2}\right)\,+\,\left({\bf q}_{r1}\,{\bf q}_{d2}\,+\,{\bf q}_{d1}\,{\bf q}_{r2}\right)\,\varepsilon\,,  \label{new24}
\end{equation}
together with their individual definitions
\begin{equation}
{\mathbb Q}_{z1} =\, {\bf q}_{r1} + {\bf q}_{d1}\,\varepsilon\;\;\;\text{and}\;\;\;{\mathbb Q}_{z2} =\, {\bf q}_{r2} + {\bf q}_{d2}\,\varepsilon\,.
\end{equation}
If we now use the fact that ${\varepsilon}$, along with ${\varepsilon^{\dagger}=\varepsilon}$ and ${\varepsilon^2=1}$, commutes with every element of ${{\cal K}^{\lambda}}$ defined in (\ref{clind}) and consequently with all ${{\bf q}_r}$, ${{\bf q}^{\dagger}_r}$, ${{\bf q}_d}$ and ${{\bf q}^{\dagger}_d}$, and work out ${{\mathbb Q}^{\dagger}_{z1}}$, ${{\mathbb Q}^{\dagger}_{z2}}$ and the products ${{\mathbb Q}_{z1}{\mathbb Q}^{\dagger}_{z1}}$, ${\,{\mathbb Q}_{z2}{\mathbb Q}^{\dagger}_{z2}}$ and ${({\mathbb Q}_{z1}\,{\mathbb Q}_{z2})^{\dagger}}$ as
\begin{align}
{\mathbb Q}^{\dagger}_{z1} &=\, {\bf q}^{\dagger}_{r1} + {\bf q}^{\dagger}_{d1}\,\varepsilon\,, \\
{\mathbb Q}^{\dagger}_{z2} &=\, {\bf q}^{\dagger}_{r2} + {\bf q}^{\dagger}_{d2}\,\varepsilon\,, \\
{\mathbb Q}_{z1}\,{\mathbb Q}^{\dagger}_{z1} &=\,\left({\bf q}_{r1}\,{\bf q}^{\dagger}_{r1}\,+\,{\bf q}_{d1}\,{\bf q}^{\dagger}_{d1}\right)\,+\,\left({\bf q}_{r1}\,{\bf q}^{\dagger}_{d1}\,+\,{\bf q}_{d1}\,{\bf q}^{\dagger}_{r1}\right)\,\varepsilon\,, \label{t-39} \\
{\mathbb Q}_{z2}\,{\mathbb Q}^{\dagger}_{z2} &=\,\left({\bf q}_{r2}\,{\bf q}^{\dagger}_{r2}\,+\,{\bf q}_{d2}\,{\bf q}^{\dagger}_{d2}\right)\,+\,\left({\bf q}_{r2}\,{\bf q}^{\dagger}_{d2}\,+\,{\bf q}_{d2}\,{\bf q}^{\dagger}_{r2}\right)\,\varepsilon\,, \label{t-40} \\
\text{and}\;\;\;({\mathbb Q}_{z1}\,{\mathbb Q}_{z2})^{\dagger}\,&=\;\left({\bf q}^{\dagger}_{r2}\,{\bf q}^{\dagger}_{r1}\,+\,{\bf q}^{\dagger}_{d2}\,{\bf q}^{\dagger}_{d1}\right)\,+\,\left({\bf q}^{\dagger}_{d2}\,{\bf q}^{\dagger}_{r1}\,+\,{\bf q}^{\dagger}_{r2}\,{\bf q}^{\dagger}_{d1}\right)\,\varepsilon\,, \label{new30}
\end{align}
then, using the {\it same} normalization condition ${{\bf q}_{r}\,{\bf q}^{\dagger}_{d}+{\bf q}_{d}\,{\bf q}^{\dagger}_{r}=0}$ of (\ref{ken-norm}), the norm relation (\ref{norms}) is not difficult to verify. To that end, we first work out the geometric product ${({\mathbb Q}_{z1}\,{\mathbb Q}_{z2})({\mathbb Q}_{z1}\,{\mathbb Q}_{z2})^{\dagger}}$ using expressions (\ref{new24}) and (\ref{new30}), which gives
\begin{align}
({\mathbb Q}_{z1}\,{\mathbb Q}_{z2})({\mathbb Q}_{z1}\,{\mathbb Q}_{z2})^{\dagger} =
&\left\{\left({\bf q}_{r1}\,{\bf q}_{r2}\,+\,{\bf q}_{d1}\,{\bf q}_{d2}\right)
\left({\bf q}^{\dagger}_{r2}\,{\bf q}^{\dagger}_{r1}\,+\,{\bf q}^{\dagger}_{d2}\,{\bf q}^{\dagger}_{d1}\right)
+\left({\bf q}_{r1}\,{\bf q}_{d2}\,+\,{\bf q}_{d1}\,{\bf q}_{r2}\right)
\left({\bf q}^{\dagger}_{d2}\,{\bf q}^{\dagger}_{r1}\,+\,{\bf q}^{\dagger}_{r2}\,{\bf q}^{\dagger}_{d1}\right)\right\} \notag \\
&+\left\{\left({\bf q}_{r1}\,{\bf q}_{d2}\,+\,{\bf q}_{d1}\,{\bf q}_{r2}\right)
\left({\bf q}^{\dagger}_{r2}\,{\bf q}^{\dagger}_{r1}\,+\,{\bf q}^{\dagger}_{d2}\,{\bf q}^{\dagger}_{d1}\right)
+\left({\bf q}_{r1}\,{\bf q}_{r2}\,+\,{\bf q}_{d1}\,{\bf q}_{d2}\right)
\left({\bf q}^{\dagger}_{d2}\,{\bf q}^{\dagger}_{r1}\,+\,{\bf q}^{\dagger}_{r2}\,{\bf q}^{\dagger}_{d1}\right)\right\}\varepsilon\,. \label{prod}
\end{align}
Now the ``real" part of the above product simplifies to (\ref{s32}) as follows:
\begin{align}
\left\{({\mathbb Q}_{z1}\,{\mathbb Q}_{z2})({\mathbb Q}_{z1}\,{\mathbb Q}_{z2})^{\dagger}\right\}_{real}
&={\bf q}_{r1}\,{\bf q}_{r2}\,{\bf q}^{\dagger}_{r2}\,{\bf q}^{\dagger}_{r1}
\,+\,
{\bf q}_{d1}\,{\bf q}_{d2}\,{\bf q}^{\dagger}_{r2}\,{\bf q}^{\dagger}_{r1}
\,+\,
{\bf q}_{r1}\,{\bf q}_{r2}\,{\bf q}^{\dagger}_{d2}\,{\bf q}^{\dagger}_{d1}
\,+\,
{\bf q}_{d1}\,{\bf q}_{d2}\,{\bf q}^{\dagger}_{d2}\,{\bf q}^{\dagger}_{d1} \notag \\
&\,\;\;\;\;\;\;\;+\,
{\bf q}_{r1}\,{\bf q}_{d2}\,{\bf q}^{\dagger}_{d2}\,{\bf q}^{\dagger}_{r1}
\,+\,
{\bf q}_{d1}\,{\bf q}_{r2}\,{\bf q}^{\dagger}_{d2}\,{\bf q}^{\dagger}_{r1}
\,+\,
{\bf q}_{r1}\,{\bf q}_{d2}\,{\bf q}^{\dagger}_{r2}\,{\bf q}^{\dagger}_{d1} \label{s30}
\,+\,
{\bf q}_{d1}\,{\bf q}_{r2}\,{\bf q}^{\dagger}_{r2}\,{\bf q}^{\dagger}_{d1} \\
&={\bf q}_{r1}\,{\bf q}_{r2}\,{\bf q}^{\dagger}_{r2}\,{\bf q}^{\dagger}_{r1}
\,+\,
{\bf q}_{d1}\,{\bf q}_{d2}\,{\bf q}^{\dagger}_{d2}\,{\bf q}^{\dagger}_{d1}
\,+\,
{\bf q}_{r1}\,{\bf q}_{d2}\,{\bf q}^{\dagger}_{d2}\,{\bf q}^{\dagger}_{r1} \label{s31}
\,+\,
{\bf q}_{d1}\,{\bf q}_{r2}\,{\bf q}^{\dagger}_{r2}\,{\bf q}^{\dagger}_{d1} \\
&=\varrho^2_{r1}\,\varrho^2_{r2}+\varrho^2_{d1}\,\varrho^2_{d2}+\varrho^2_{r1}\,\varrho^2_{d2}+\varrho^2_{d1}\,\varrho^2_{r2}\,. \label{s32}
\end{align}
Here (\ref{s31}) follows from (\ref{s30}) upon inserting the normalization condition (\ref{ken-norm}) in the form ${{\bf q}_r\,{\bf q}^{\dagger}_d =-\,{\bf q}_d\,{\bf q}^{\dagger}_r}$ into the second and third terms of (\ref{s30}), which then cancel out with the sixth and seventh terms of (\ref{s30}), respectively; and (\ref{s32}) follows from (\ref{s31}) upon inserting the normalization conditions ${||{\bf q}||^2={\bf q}{\bf q}^{\dagger}=\varrho^2}$ for the real and dual quaternions specified in (\ref{s-r}) and (\ref{d-d}), for each of the four terms of (\ref{s31}). Similarly, the ``dual" part of the product (\ref{prod}) simplifies${\;}$to
\begin{align}
\left\{({\mathbb Q}_{z1}\,{\mathbb Q}_{z2})({\mathbb Q}_{z1}\,{\mathbb Q}_{z2})^{\dagger}\right\}_{dual}
&=\Big\{{\bf q}_{r1}\,{\bf q}_{d2}\,{\bf q}^{\dagger}_{r2}\,{\bf q}^{\dagger}_{r1}
\,+\,
{\bf q}_{d1}\,{\bf q}_{r2}\,{\bf q}^{\dagger}_{r2}\,{\bf q}^{\dagger}_{r1}
\,+\,
{\bf q}_{r1}\,{\bf q}_{d2}\,{\bf q}^{\dagger}_{d2}\,{\bf q}^{\dagger}_{d1}
\,+\,
{\bf q}_{d1}\,{\bf q}_{r2}\,{\bf q}^{\dagger}_{d2}\,{\bf q}^{\dagger}_{d1} \notag \\
&\,\;\;\;\;\;\;\;+\,
{\bf q}_{r1}\,{\bf q}_{r2}\,{\bf q}^{\dagger}_{d2}\,{\bf q}^{\dagger}_{r1}
\,+\,
{\bf q}_{d1}\,{\bf q}_{d2}\,{\bf q}^{\dagger}_{d2}\,{\bf q}^{\dagger}_{r1}
\,+\,
{\bf q}_{r1}\,{\bf q}_{r2}\,{\bf q}^{\dagger}_{r2}\,{\bf q}^{\dagger}_{d1}
\,+\,
{\bf q}_{d1}\,{\bf q}_{d2}\,{\bf q}^{\dagger}_{r2}\,{\bf q}^{\dagger}_{d1}\Big\}\,\varepsilon \label{s33}\\
&= 0\,. \label{s34}
\end{align}
We can see this again by inserting into (\ref{s33}) the normalization condition (\ref{ken-norm}) in the form ${{\bf q}_r\,{\bf q}^{\dagger}_d =-\,{\bf q}_d\,{\bf q}^{\dagger}_r}$ and the normalization conditions ${||{\bf q}||^2={\bf q}{\bf q}^{\dagger}=\varrho^2}$ for the quaternions in (\ref{s-r}) and (\ref{d-d}), which cancels out the first four terms of (\ref{s33}) with the last four. Consequently, combining the results of (\ref{s32}) and (\ref{s34}), for the left-hand side of (\ref{norms}) we have
\begin{equation}
||{\mathbb Q}_{z1}\,{\mathbb Q}_{z2}||\,=\,\sqrt{\varrho^2_{r1}\,\varrho^2_{r2}+\varrho^2_{r1}\,\varrho^2_{d2}+\varrho^2_{d1}\,\varrho^2_{r2}+\varrho^2_{d1}\,\varrho^2_{d2}}\,. \label{s35}
\end{equation}
Thus, comparing the results in (\ref{s35}) and (\ref{s36}), we finally arrive at the relation 
\begin{equation}
||{\mathbb Q}_{z1}\,{\mathbb Q}_{z2}||\,=\;||{\mathbb Q}_{z1}||\;||{\mathbb Q}_{z2}||\,, \label{norm-final}
\end{equation}
which is evidently the same as the norm relation (\ref{norms}) in every respect apart from the appropriate change in notation. This result is facilitated by the definition (b) of the norm [or of the quadratic form ${{\cal Q}({\bf X})}$] explained below Eq.~(\ref{b}). We have thus proved that the finite-dimensional algebra ${{\cal K}^{\lambda}}$ over the reals can be equipped with a positive definite quadratic form ${\cal Q}$ (the square of the norm) such that ${{\cal Q}({\bf X}{\bf Y}) = {\cal Q}({\bf X})\,{\cal Q}({\bf Y})}$ for all ${\bf X}$ and ${\bf Y}$ in ${{\cal K}^{\lambda}}$. Consequently, a product ${{\bf X}{\bf Y}}$ would vanish if and only if ${\bf X}$ or ${\bf Y}$ vanishes. In other words, ${{\cal K}^{\lambda}}$, equipped with ${\cal Q}$, is a division algebra. In Appendix \ref{B} we prove the composition law $||{\mathbb Q}_{z1}\,{\mathbb Q}_{z2}||^2=||{\mathbb Q}_{z1}||^2\,||{\mathbb Q}_{z2}||^2$ in full generality without assuming (\ref{ken-norm}), and in Appendix \ref{C} we prove that the orthogonality of the quaternions ${{\bf q}_{r}}$ and ${{\bf q}_{d}}$ is preserved under multiplication.

Without loss of generality we can now restrict ${{\cal K}^{\lambda}}$ in (\ref{new18}) to a unit 7-sphere by setting the radii ${\varrho_r}$ and ${\varrho_d}$ to ${\frac{1}{\sqrt{2}}}$:
\begin{equation}
{\cal K}^{\lambda}\supset\,S^7:=\,\left\{\,{\mathbb Q}_z :=\, {\bf q}_r + {\bf q}_d\,\varepsilon\;\Big|\; ||{\mathbb Q}_z||=1\;\,\text{and}\,\;{\bf q}_{r}\,{\bf q}^{\dagger}_{d}+{\bf q}_{d}\,{\bf q}^{\dagger}_{r}=0\,\right\}, \label{spring}
\end{equation}
where ${\varepsilon =-\lambda\,I_3{\bf e}_{\infty}\,}$, ${\varepsilon^{\dagger}=\varepsilon}$, ${\,\varepsilon^2=\,{\bf e}^2_{\infty}=\,+1\,}$,
\begin{equation}
{\bf q}_r =\,q_0 + q_1\,\lambda\,{\bf e}_x{\bf e}_y+q_2\,\lambda\,{\bf e}_z{\bf e}_x+q_3\,\lambda\,{\bf e}_y{\bf e}_z\,,\;\;\,
\text{and}\;\;\,{\bf q}_d =\,-q_7 + q_6\,{\bf e}_x{\bf e}_y+q_5\,{\bf e}_z{\bf e}_x+q_4\,{\bf e}_y{\bf e}_z\,,
\end{equation}
so that
\begin{equation}
{\mathbb Q}_z=\,q_0+q_1\,\lambda{\bf e}_x{\bf e}_y+q_2\,\lambda{\bf e}_z{\bf e}_x+q_3\,\lambda{\bf e}_y{\bf e}_z+q_4\,\lambda{\bf e}_x{\bf e}_{\infty}+q_5\,\lambda{\bf e}_y{\bf e}_{\infty}+q_6\,\lambda{\bf e}_z{\bf e}_{\infty}+q_7\,\lambda I_3{\bf e}_{\infty}\,. \label{snons}
\end{equation}
Needless to say, since all Clifford algebras are associative algebras by definition, unlike the non-associative octonionic algebra the 7-sphere we have constructed here corresponds to an {\it associative} (but non-commutative) division algebra.

Note that in terms of the components of ${\bf q}_r$ and ${\bf q}_d$ the condition ${\bf q}_r {\bf q}^{\dagger}_d + {\bf q}_d {\bf q}^{\dagger}_r = 0$ is equivalent to the constraint
\begin{equation}
f_{\!_{\cal K}} = -\,q_0q_7 + \lambda\,q_1q_6 + \lambda\,q_2q_5 + \lambda\,q_3q_4 = 0. \label{50}
\end{equation}
This constraint reduces the space ${\cal K}^{\lambda}$ to the sphere $S^7$, thereby reducing the 8 dimensions of ${\cal K}^{\lambda}$ to the 7 dimensions of $S^7$ defined in (\ref{spring}). But the 7-sphere thus constructed has a topology \cite{Milnor} that is different from that of the octonionic 7-sphere, and the difference between the two is captured by the difference in the corresponding normalizing constraints
\begin{equation}
f(q_0, q_1, q_2, q_3, q_4, q_5, q_6, q_7) = 0.
\end{equation}
More precisely, the two normalizing constraints giving rise to the two topologically distinct 7-spheres of radius ${\rho}$ are:
\begin{equation}
f_{\!_{\mathbb O}} = q^2_0 + q^2_1 + q^2_2 + q^2_3 + q^2_4 + q^2_5 + q^2_6 + q^2_7 - \rho^2 = 0, \label{c-1}
\end{equation}
which reduces the set ${\mathbb O}$ of unit octonions to the sphere $S^7$ made up of eight-dimensional vectors of fixed length $\rho$, and
\begin{equation}
f_{\!_{\cal K}} = -\,q_0q_7 + \lambda\,q_1q_6 + \lambda\,q_2q_5 + \lambda\,q_3q_4 = 0, \label{c-2}\tag{\ref{50}}
\end{equation}
which reduces the set ${\cal K}^{\lambda}$ to the sphere $S^7$ made up of a different collection of eight-dimensional vectors of fixed length ${\rho=\sqrt{\varrho^2_r + \varrho^2_d}}$. Both constraints, (\ref{c-2}) and (\ref{c-1}), involve the same eight variables of the embedding space ${\rm I\!R}^8$, namely, $q_0$, $q_1$, $q_2$, $q_3$, $q_4$, $q_5$, $q_6$, and $q_7$, giving the same dimensions for the sphere $S^7$ of radius $\rho$, albeit respecting different topologies \cite{Milnor}. This difference arises because we have used the geometric product ${{\bf X}{\bf X}^{\dagger}}$ rather than the scalar product ${{\bf X}\cdot{\bf X}^{\dagger}}$ to derive the constraint $f_{\!_{\cal K}} = 0$. But both definitions of the norm $||{\bf X}||$ give identical results, as explained$\;$above. 

Given the quadratic form ${{\cal Q}({\bf X})}$ and the norm relation (\ref{norm-final}), we may now view the four {\it associative} normed division algebras in the only possible dimensions 1, 2, 4 and 8, respectively \cite{Hurwitz-1898,Hurwitz-1923}, as even sub-algebras of the Clifford algebras
\begin{align}
Cl^{\lambda}_{1,0}&={\rm span}\!\left\{\,1,\;\lambda{\bf e}_x\,\right\}\!, \\
Cl^{\lambda}_{2,0}&={\rm span}\!\left\{\,1,\;\lambda{\bf e}_x,\,\lambda{\bf e}_y,\,\lambda{\bf e}_x{\bf e}_y\,\right\}\!, \\
Cl^{\lambda}_{3,0}&={\rm span}\!\left\{\,1,\;\lambda{\bf e}_x,\,\lambda{\bf e}_y,\,\lambda{\bf e}_z,\;\lambda{\bf e}_x{\bf e}_y,\,
\lambda{\bf e}_z{\bf e}_x,\,\lambda{\bf e}_y{\bf e}_z,\;\lambda{\bf e}_x{\bf e}_y{\bf e}_z\,\right\}\!, \\
\text{and}\;\;Cl^{\lambda}_{4,0}&={\rm span}\!\big\{\,1,\;\lambda{\bf e}_x,\,\lambda{\bf e}_y,\,\lambda{\bf e}_z,\,\lambda{\bf e}_{\infty},\;\lambda{\bf e}_x{\bf e}_y,\,
\lambda{\bf e}_z{\bf e}_x,\,\lambda{\bf e}_y{\bf e}_z,\;\lambda{\bf e}_x{\bf e}_{\infty},\,\lambda{\bf e}_y{\bf e}_{\infty},\,\lambda{\bf e}_z{\bf e}_{\infty}, \notag \\
&\;\;\;\;\;\;\;\;\;\;\;\;\;\;\;\;\;\;\;\;\;\;\;\;\;\;\;\;\;\;\;\;\;\;\;\;\;\;\;\;\;\;\;\;\;\;\;\;\;\;\;\;\;\;\;\;\;\;\;\;\;\;\;\;\;\;\;\;\;\;\;\lambda{\bf e}_x{\bf e}_y{\bf e}_z,\,\lambda{\bf e}_x{\bf e}_y{\bf e}_{\infty},\,\lambda{\bf e}_z{\bf e}_x{\bf e}_{\infty},\,\lambda{\bf e}_y{\bf e}_z{\bf e}_{\infty},\;\lambda{\bf e}_x{\bf e}_y{\bf e}_z{\bf e}_{\infty}\,\big\}\!.
\end{align}
It is easy to verify that the even subalgebras of ${Cl^{\lambda}_{1,0}}$, ${Cl^{\lambda}_{2,0}}$ and ${Cl^{\lambda}_{3,0}}$ are indeed isomorphic to ${\mathbb R}$, ${\mathbb C}$ and ${\mathbb H}$, respectively.

In practice, the above eight-dimensional algebra sometimes appears in the guise of a `1d up' approach to Conformal Geometric Algebra in the engineering and computer vision applications \cite{Lasenby-1,Lasenby-2}. Such physical applications would benefit from explicitly using the quadratic form ${{\cal Q}({\bf X})}$ and the corresponding division algebra we have presented in this paper. For instance, it may help in removing the ``singularities" or non-zero zero divisors from occurring in such applications. An illustration of how that may work can be found in Ref.~\cite{RSOS} where we have applied the quadratic form ${{\cal Q}({\bf X})}$ and the corresponding division algebra to understand the geometrical origins of quantum correlations within the 7-sphere constructed in this paper. In the broader context of relativistic quantum theory, it is well known that between 1932 and 1952 Jordan attempted to use an alternative ring of octonions with non-associative multiplication rules to transfer the probabilistic interpretation of quantum theory to what is now known as exceptional Jordan algebra \cite{Jordan}. But as Dirac has noted \cite{Dirac}, Jordan's attempt to obtain a generalized quantum theory in this manner was not successful, because the non-associative multiplication rules are not compatible with any physically meaningful group of transformations such as the Poincar\'e group. However, the octonion-like algebra ${\cal K}^{\lambda}$ with six rather than seven imaginaries we have presented in this paper is {\it associative} by construction, and therefore it will be amenable to Jordan type application to quantum theory. Apart from these applications, in Section 5 of Ref.~\cite{Baez} Baez has discussed more mathematically oriented applications of the norm division algebras in four dimensions. These application can now be extended to eight dimensions, thanks to the associativity of ${\cal K}^{\lambda}$. The Clifford-algebraic investigations by Lounesto\break in normed division algebras and octonions may also benefit from the associativity of ${\cal K}^{\lambda}$ \cite{Lounesto}. To facilitate these applications, in Appendix \ref{A} we illustrate how non-zero zero divisors are precluded from the ${\cal K}^{\lambda}$ equipped with ${{\cal Q}({\bf X})}$.

\appendix

\section{Illustration of How the Definition (b) of the Norm Precludes Zero Divisors} \label{A}

According to Frobenius theorem \cite{Forbenius} --- which uses scalar products (instead of geometric products we have used) as an essential ingredient in its proof, a finite-dimensional associative division algebra over the reals is necessarily isomorphic to either ${\mathbb R}$, ${\mathbb C}$, or ${\mathbb H}$ in the 1, 2, and 4 dimensions, respectively.
Since Clifford algebras are finite-dimensional associative algebras, Frobenius theorem suggests that those Clifford algebras that are not isomorphic to ${\mathbb R}$, ${\mathbb C}$, or ${\mathbb H}$ may contain non-zero zero divisors or idempotent elements. It is therefore important to understand how the definition (b) of the norm leading to the quadratic form ${{\cal Q}({\bf X})}$ prevents non-zero zero divisors from occurring in the algebra ${\cal K}^{\lambda}$.

To that end, recall that the elements of ${{\cal K}^{\lambda}}$ are of the following general form in terms of quaternions ${{\bf q}_r}$ and ${{\bf q}_d}$:
\begin{equation}
{\mathbb Q}_z =\,{\bf q}_r + {\bf q}_d\,\varepsilon\,, \label{A-1}
\end{equation}
where $\varepsilon^{\dagger}=\varepsilon$ and $\varepsilon^2=+1$ as in (\ref{16-A}) and the normalization of ${\mathbb Q}_z$ requires that ${{\bf q}_r}$ and ${{\bf q}_d}$ must satisfy the condition
\begin{equation}
{\bf q}_r\,{\bf q}^{\dagger}_d+{\bf q}_d\,{\bf q}^{\dagger}_r = 0\;\;\;\Longleftrightarrow\;\;\;
{\bf q}_r\,{\bf q}^{\dagger}_d = -{\bf q}_d\,{\bf q}^{\dagger}_r. \label{A-2}
\end{equation}
This condition follows from the definition (b) of the norm $||{\mathbb Q}_z||$ discussed below Eq.~(\ref{b}). It respects the fundamental geometric product ${\mathbb Q}_z{\mathbb Q}_z^{\dagger}$ and gives the same scalar value for the norm $||{\mathbb Q}_z||$ as that calculated using definition (a).

Now, for the sake of argument, consider the following idempotent quantities as candidate non-zero zero divisors:
\begin{equation}
{\bf Z}_{\pm}=\frac{1}{2}(1 \pm \varepsilon). \label{A-3}
\end{equation}
We call the quantities ${\bf Z}_{\pm}$ idempotent quantities because they square to themselves, which can be easily verified:
\begin{equation}
{\bf Z}^2_{+}=
{\bf Z}_{+}{\bf Z}^{\dagger}_{+}={\bf Z}_{+}
\;\;\;\text{and}\;\;\;
{\bf Z}^2_{-}={\bf Z}_{-}{\bf Z}^{\dagger}_{-}={\bf Z}_{-}. \label{A-4}
\end{equation}
But ${\bf Z}_{+}$ and ${\bf Z}_{-}$ are also orthogonal to each other because their products vanish, which can also be easily verified:
\begin{equation}
{\bf Z}_{+}{\bf Z}_{-}={\bf Z}_{-}{\bf Z}_{+}=0. \label{A-5}
\end{equation}
Now consider two multivectors, ${\bf X}$ and ${\bf Y}$, confined to a two-dimensional subspace of ${\cal K}^{\lambda}$, by setting ${X_0 = Y_0 = \frac{1}{\sqrt{2}}}$, ${Y_7=-X_7=\frac{1}{\sqrt{2}}}$, and the remaining twelve coefficients equal to zero in the Eqs.~(\ref{X}) and (\ref{Y}), along with $\varepsilon \equiv\,-\lambda\,I_3{\bf e}_{\infty}$:
\begin{equation}
{\bf X} =\frac{1}{\sqrt{2}}(1 + \varepsilon)= \sqrt{2}\,{\bf Z}_{+}\;\;\;\text{and}\;\;\;{\bf Y}=\frac{1}{\sqrt{2}}(1 - \varepsilon)= \sqrt{2}\,{\bf Z}_{-}. \label{A-6}
\end{equation}
Then (\ref{A-5}) implies that $||{\bf X}{\bf Y}||=0$, which can be verified by substituting for the multivectors ${\bf X}$ and ${\bf Y}$ from (\ref{A-6}):
\begin{align}
||{\bf X}{\bf Y}||&=\left|\left|\frac{1}{\sqrt{2}}(1+\varepsilon)\,\frac{1}{\sqrt{2}}(1-\varepsilon)\right|\right| \label{maga} \\
&=\left|\left|\frac{1}{2}\left(1-\varepsilon^2\right)\right|\right| \\
&=||\,0\,|| \\
&=0, \label{12}
\end{align}
where $\varepsilon^2=1$ is used. Next, using $\varepsilon^{\dagger}=\varepsilon$ and $\varepsilon^2=1$, we evaluate the right-hand side of the norm relation (\ref{norms}), giving
\begin{align}
||{\bf X}||\,||{\bf Y}||&=\left|\left|\frac{1}{\sqrt{2}}(1+\varepsilon)\right|\right|\left|\left|\frac{1}{\sqrt{2}}(1-\varepsilon)\right|\right| \\
&=\left(\sqrt{\frac{1}{\sqrt{2}}(1+\varepsilon)\,\frac{1}{\sqrt{2}}(1+\varepsilon)^{\dagger}}\;\right)\left(\sqrt{\frac{1}{\sqrt{2}}(1-\varepsilon)\,\frac{1}{\sqrt{2}}(1-\varepsilon)^{\dagger}}\;\right) \\
&=\left(\sqrt{\frac{1}{2}(1+\varepsilon)(1+\varepsilon)}\;\right)\left(\sqrt{\frac{1}{2}(1-\varepsilon)(1-\varepsilon)}\;\right) \\
&=\left(\sqrt{(1+\varepsilon)}\,\right)\left(\sqrt{(1-\varepsilon)}\,\right) \label{30} \\
&=\sqrt{(1+\varepsilon)(1-\varepsilon)} \\
&=\sqrt{\left(1-\varepsilon^2\right)\,} \\
&=\sqrt{0\,} \\
&=0, \label{20}
\end{align}
where we have used geometric products to evaluate the norms. Comparing (\ref{12}) and (\ref{20}) we see that the norm relation $||{\bf X}{\bf Y}||=||{\bf X}||\,||{\bf Y}||$ is satisfied for the multivectors in (\ref{A-6}), despite (\ref{A-4}) and (\ref{A-5}). On the other hand, if we insist on using scalar products for evaluating the norms, then (\ref{12}) remains the same but instead of (\ref{20}) we obtain
\begin{align}
||{\bf X}||\,||{\bf Y}||&=\left|\left|\frac{1}{\sqrt{2}}(1+\varepsilon)\right|\right|\left|\left|\frac{1}{\sqrt{2}}(1-\varepsilon)\right|\right| \\
&=\left(\sqrt{\frac{1}{2}(1+\varepsilon)\cdot(1+\varepsilon)^{\dagger}}\;\right)\left(\sqrt{\frac{1}{2}(1-\varepsilon)\cdot(1-\varepsilon)^{\dagger}}\;\right) \label{not20} \\
&=\left(\sqrt{\frac{1}{2}(1+\varepsilon)\cdot(1+\varepsilon)}\;\right)\left(\sqrt{\frac{1}{2}(1-\varepsilon)\cdot(1-\varepsilon)}\;\right) \\
&=\left(\sqrt{\frac{1}{2}(1+1)}\,\right)\left(\sqrt{\frac{1}{2}(1+1)}\,\right) \label{30a} \\
&=\left(\sqrt{1}\right)\left(\sqrt{1}\right) \\
&=1, \label{20a}
\end{align}
which seems to imply that $||{\bf X}{\bf Y}||\not=||{\bf X}||\,||{\bf Y}||$. This is because the norms are evaluated inconsistently in arriving at the contradictory results (\ref{12}) and (\ref{20a}). While the product between ${\bf X}$ and ${\bf Y}$ for (\ref{12}) is evaluated using the geometric product giving ${{\bf X}{\bf Y}=0}$ so that $||{\bf X}{\bf Y}||=0$, the norms $||{\bf X}||$ and $||{\bf Y}||$ are evaluated for (\ref{20a}) using the scalar products giving $||{\bf X}||=1$ and $||{\bf Y}||=1$. But using scalar products to evaluate norms is inconsistent with the choices made in (\ref{A-6}) for the coefficients of ${\bf X}$ and ${\bf Y}$. To appreciate this, recall again that the fundamental product in Geometric Algebra is the geometric product, not the scalar product, and the geometric product such as ${\bf X}{\bf X}^{\dagger}$ is worked out in Eq.~(\ref{21}) above, from which it is clear that the norm $||{\bf X}||=\sqrt{{\bf X}{\bf X}^{\dagger}}$ can reduce to a scalar quantity {\it if and only if} $\,{\bf q}_{r}\,{\bf q}^{\dagger}_{d}\,+\,{\bf q}_{d}\,{\bf q}^{\dagger}_{r}=0$. And, as we saw in (\ref{22}), in terms of the coefficients of ${\bf X}$ this condition is equivalent to 
\begin{equation}
-X_0X_7 + \lambda\,X_1X_6 + \lambda\,X_2X_5 + \lambda\,X_3X_4 = 0. \label{50a}
\end{equation}
It is now easy to see that the choices of the coefficients in ${\bf X}$ and ${\bf Y}$ are incompatible with this condition. The choices of coefficients made in (\ref{A-6}) are $X_0=\frac{1}{\sqrt{2}}$ and $X_7=-\frac{1}{\sqrt{2}}$ for ${\bf X}$ and $Y_0=\frac{1}{\sqrt{2}}$ and $Y_7=\frac{1}{\sqrt{2}}$ for ${\bf Y}$, with all other coefficients set to zero. Substituting these values in (\ref{50a}) we immediately arrive at the contradictions $\pm\frac{1}{2}=0$, proving that the {\it ad hoc} coefficients chosen to define ${\bf X}$ and ${\bf Y}$ are not compatible with the use of scalar products to evaluate the norms. In other words, the {\it ad hoc} coefficients chosen in (\ref{A-6}) to define ${\bf X}$ and ${\bf Y}$ are not compatible with the condition (\ref{A-2}), and thus they are not compatible with the definition (b) for the norms. On the other hand, if the norms are evaluated consistently on both sides of Eq.~(\ref{norms}), as for (\ref{12}) and (\ref{20}), then the norm relation $||{\bf X}{\bf Y}|| = |||{\bf X}||\,||{\bf Y}||$ holds for any elements ${\bf X}$ and ${\bf Y}$ in ${{\cal K}^{\lambda}}$. This completes the illustration of how using the definition (b) for the norms precludes non-zero zero divisors from ${{\cal K}^{\lambda}}$, and leads us to the following general proof of the composition law for ${{\cal K}^{\lambda}}$.

\section{Proof of the Composition Law for ${\cal K}^{\lambda}$ without Assuming ${{\bf q}_{r}\,{\bf q}^{\dagger}_{d}+{\bf q}_{d}\,{\bf q}^{\dagger}_{r}=0}$} \label{B}

In the proof of the norm relations (\ref{norm-final}) we assumed the normalization condition ${{\bf q}_{r}\,{\bf q}^{\dagger}_{d}+{\bf q}_{d}\,{\bf q}^{\dagger}_{r}=0}$. It turns out, however, that the composition law holds for the algebra ${\cal K}^{\lambda}$ even without assuming this condition. In this appendix we first prove the composition law explicitly and then obtain the norm relation (\ref{norm-final}) as its special case. To this end, recall from Eq.~(\ref{23}) that the product ${\mathbb Q}_{z}\,{\mathbb Q}^{\dagger}_{z}$ between the general element ${\mathbb Q}_{z}$ in ${\cal K}^{\lambda}$ and its conjugate ${\mathbb Q}^{\dagger}_{z}$ is of the form
\begin{equation}
{\mathbb Q}_{z}\,{\mathbb Q}^{\dagger}_{z}\,=\,\text{(a scalar)} \,+\, \text{(a scalar)}\times\varepsilon\,, \label{B1}
\end{equation}
with $\varepsilon^2=+1$. Thus, the square of the norm resembles a split complex number \cite{Dray} rather than a real scalar. However, the composition law still holds. To prove this, we begin by evaluating the right-hand side of (\ref{norm-final}) using (\ref{t-39}) and (\ref{t-40}):
\begin{align}
||{\mathbb Q}_{z1}||^2\,||{\mathbb Q}_{z2}||^2 &= \left({\mathbb Q}_{z1}\,{\mathbb Q}^{\dagger}_{z1}\right)\left({\mathbb Q}_{z2}\,{\mathbb Q}^{\dagger}_{z2}\right) \\
&=\left\{\left({\bf q}_{r1}\,{\bf q}^{\dagger}_{r1}+{\bf q}_{d1}\,{\bf q}^{\dagger}_{d1}\right)+\left({\bf q}_{r1}\,{\bf q}^{\dagger}_{d1}+{\bf q}_{d1}\,{\bf q}^{\dagger}_{r1}\right)\varepsilon\right\}\left\{\left({\bf q}_{r2}\,{\bf q}^{\dagger}_{r2}+{\bf q}_{d2}\,{\bf q}^{\dagger}_{d2}\right)+\left({\bf q}_{r2}\,{\bf q}^{\dagger}_{d2}+{\bf q}_{d2}\,{\bf q}^{\dagger}_{r2}\right)\varepsilon\right\} \\
&= \left\{\left({\bf q}_{r1}\,{\bf q}^{\dagger}_{r1}+{\bf q}_{d1}\,{\bf q}^{\dagger}_{d1}\right) \left({\bf q}_{r2}\,{\bf q}^{\dagger}_{r2}+{\bf q}_{d2}\,{\bf q}^{\dagger}_{d2}\right)\,+\, \left({\bf q}_{r1}\,{\bf q}^{\dagger}_{d1}+{\bf q}_{d1}\,{\bf q}^{\dagger}_{r1}\right)\left({\bf q}_{r2}\,{\bf q}^{\dagger}_{d2}+{\bf q}_{d2}\,{\bf q}^{\dagger}_{r2}\right)
\right\} \notag \\
&\;\;\;\;\;\;\;\;\;\;+\left\{\left({\bf q}_{r1}\,{\bf q}^{\dagger}_{r1}+{\bf q}_{d1}\,{\bf q}^{\dagger}_{d1}\right) \left({\bf q}_{r2}\,{\bf q}^{\dagger}_{d2}+{\bf q}_{d2}\,{\bf q}^{\dagger}_{r2}\right)\,+\,\left({\bf q}_{r1}\,{\bf q}^{\dagger}_{d1}+{\bf q}_{d1}\,{\bf q}^{\dagger}_{r1}\right)\left({\bf q}_{r2}\,{\bf q}^{\dagger}_{r2}+{\bf q}_{d2}\,{\bf q}^{\dagger}_{d2}\right)
\right\}\varepsilon    \\
&=\left\{\left(\varrho^2_{r1}+\varrho^2_{d1}\right)\left(\varrho^2_{r2}+\varrho^2_{d2}\right) \,+\,\left({\bf q}_{r1}\,{\bf q}^{\dagger}_{d1}+{\bf q}_{d1}\,{\bf q}^{\dagger}_{r1}\right)\left({\bf q}_{r2}\,{\bf q}^{\dagger}_{d2}+{\bf q}_{d2}\,{\bf q}^{\dagger}_{r2}\right)\right\} \notag \\
&\;\;\;\;\;\;\;\;\;\;+\left\{\left( \varrho^2_{r1}\,+\,\varrho^2_{d1}\right) \left({\bf q}_{r2}\,{\bf q}^{\dagger}_{d2}+{\bf q}_{d2}\,{\bf q}^{\dagger}_{r2}\right)
+ \left({\bf q}_{r1}\,{\bf q}^{\dagger}_{d1}+{\bf q}_{d1}\,{\bf q}^{\dagger}_{r1}\right)\left( \varrho^2_{r2}\,+\,\varrho^2_{d2}\right)\right\}\,\varepsilon \\
&=\left\{\varrho^2_{r1}\,\varrho^2_{r2}+\varrho^2_{r1}\,\varrho^2_{d2}+\varrho^2_{d1}\,\varrho^2_{r2}+\varrho^2_{d1}\,\varrho^2_{d2} \,+\,\left({\bf q}_{r1}\,{\bf q}^{\dagger}_{d1}+{\bf q}_{d1}\,{\bf q}^{\dagger}_{r1}\right)\left({\bf q}_{r2}\,{\bf q}^{\dagger}_{d2}+{\bf q}_{d2}\,{\bf q}^{\dagger}_{r2}\right)\right\} \notag \\
&\;\;\;\;\;\;\;\;\;\;+\left\{\left( \varrho^2_{r1}\,+\,\varrho^2_{d1}\right) \left({\bf q}_{r2}\,{\bf q}^{\dagger}_{d2}+{\bf q}_{d2}\,{\bf q}^{\dagger}_{r2}\right)
+\left( \varrho^2_{r2}\,+\,\varrho^2_{d2}\right) \left({\bf q}_{r1}\,{\bf q}^{\dagger}_{d1}+{\bf q}_{d1}\,{\bf q}^{\dagger}_{r1}\right)\right\}\,\varepsilon\,, \label{B4}
\end{align}
where we have used ${\bf q}_{r1}\,{\bf q}^{\dagger}_{r1}=\varrho^2_{r1}$, ${\bf q}_{d1}\,{\bf q}^{\dagger}_{d1}=\varrho^2_{d1}$, {\it etc}. Recalling from (\ref{22}) that the quantities such as ${\bf q}_{r1}\,{\bf q}^{\dagger}_{d1}+{\bf q}_{d1}\,{\bf q}^{\dagger}_{r1}$ are scalar quantities, we see that the above product also resembles a split complex number similar to ${\mathbb Q}_{z}\,{\mathbb Q}^{\dagger}_{z}$ in (\ref{B1}).

Next, using (\ref{new24}), (\ref{new30}), and their product evaluated in (\ref{prod}), we can evaluate the left-hand side of (\ref{norm-final}) as follows:
\begin{align}
||{\mathbb Q}_{z1}\,{\mathbb Q}_{z2}||^2 &= \left({\mathbb Q}_{z1}\,{\mathbb Q}_{z2}\right)\left({\mathbb Q}_{z1}\,{\mathbb Q}_{z2}\right)^{\dagger} \\
&=\left\{\left({\bf q}_{r1}\,{\bf q}_{r2}+{\bf q}_{d1}\,{\bf q}_{d2}\right)+\left({\bf q}_{r1}\,{\bf q}_{d2}+{\bf q}_{d1}\,{\bf q}_{r2}\right)\,\varepsilon\right\}\left\{\left({\bf q}^{\dagger}_{r2}\,{\bf q}^{\dagger}_{r1}+{\bf q}^{\dagger}_{d2}\,{\bf q}^{\dagger}_{d1}\right)+\left({\bf q}^{\dagger}_{d2}\,{\bf q}^{\dagger}_{r1}+{\bf q}^{\dagger}_{r2}\,{\bf q}^{\dagger}_{d1}\right)\,\varepsilon\right\} \\ 
&=\left\{\left({\bf q}_{r1}\,{\bf q}_{r2}\,+\,{\bf q}_{d1}\,{\bf q}_{d2}\right)
\left({\bf q}^{\dagger}_{r2}\,{\bf q}^{\dagger}_{r1}\,+\,{\bf q}^{\dagger}_{d2}\,{\bf q}^{\dagger}_{d1}\right)
+\left({\bf q}_{r1}\,{\bf q}_{d2}\,+\,{\bf q}_{d1}\,{\bf q}_{r2}\right)
\left({\bf q}^{\dagger}_{d2}\,{\bf q}^{\dagger}_{r1}\,+\,{\bf q}^{\dagger}_{r2}\,{\bf q}^{\dagger}_{d1}\right)\right\} \notag \\
&\;\;\;\;\;\;+\left\{\left({\bf q}_{r1}\,{\bf q}_{r2}\,+\,{\bf q}_{d1}\,{\bf q}_{d2}\right)
\left({\bf q}^{\dagger}_{d2}\,{\bf q}^{\dagger}_{r1}\,+\,{\bf q}^{\dagger}_{r2}\,{\bf q}^{\dagger}_{d1}\right)+\left({\bf q}_{r1}\,{\bf q}_{d2}\,+\,{\bf q}_{d1}\,{\bf q}_{r2}\right)
\left({\bf q}^{\dagger}_{r2}\,{\bf q}^{\dagger}_{r1}\,+\,{\bf q}^{\dagger}_{d2}\,{\bf q}^{\dagger}_{d1}\right)
\right\}\varepsilon \\
&=\Big\{{\bf q}_{r1}\,{\bf q}_{r2}\,{\bf q}^{\dagger}_{r2}\,{\bf q}^{\dagger}_{r1}\,+\, {\bf q}_{r1}\,{\bf q}_{r2}\,{\bf q}^{\dagger}_{d2}\,{\bf q}^{\dagger}_{d1}
\,+\, {\bf q}_{d1}\,{\bf q}_{d2}\,{\bf q}^{\dagger}_{r2}\,{\bf q}^{\dagger}_{r1} \,+\, {\bf q}_{d1}\,{\bf q}_{d2}\,{\bf q}^{\dagger}_{d2}\,{\bf q}^{\dagger}_{d1} \notag \\
&\,\;\;\;\;\;\;\;+\, {\bf q}_{r1}\,{\bf q}_{d2}\,{\bf q}^{\dagger}_{d2}\,{\bf q}^{\dagger}_{r1} \,+\, 
{\bf q}_{r1}\,{\bf q}_{d2}\,{\bf q}^{\dagger}_{r2}\,{\bf q}^{\dagger}_{d1} \,+\,{\bf q}_{d1}\,{\bf q}_{r2}\,{\bf q}^{\dagger}_{d2}\,{\bf q}^{\dagger}_{r1} \,+\, {\bf q}_{d1}\,{\bf q}_{r2}\,{\bf q}^{\dagger}_{r2}\,{\bf q}^{\dagger}_{d1}\Big\} \notag \\
&\;\;\;\;\;\;\;\;\;\;\;+\Big\{{\bf q}_{r1}\,{\bf q}_{r2}\,{\bf q}^{\dagger}_{d2}\,{\bf q}^{\dagger}_{r1} \,+\, {\bf q}_{r1}\,{\bf q}_{r2}\,{\bf q}^{\dagger}_{r2}\,{\bf q}^{\dagger}_{d1} \,+\, {\bf q}_{d1}\,{\bf q}_{d2}\,{\bf q}^{\dagger}_{d2}\,{\bf q}^{\dagger}_{r1} \,+\, {\bf q}_{d1}\,{\bf q}_{d2}\,{\bf q}^{\dagger}_{r2}\,{\bf q}^{\dagger}_{d1} \notag \\
&\,\;\;\;\;\;\;\;\;\;\;\;\;\;\;\;+\, {\bf q}_{r1}\,{\bf q}_{d2}\,{\bf q}^{\dagger}_{r2}\,{\bf q}^{\dagger}_{r1} \,+\, {\bf q}_{r1}\,{\bf q}_{d2}\,{\bf q}^{\dagger}_{d2}\,{\bf q}^{\dagger}_{d1} \,+\,
{\bf q}_{d1}\,{\bf q}_{r2}\,{\bf q}^{\dagger}_{r2}\,{\bf q}^{\dagger}_{r1} \,+\, {\bf q}_{d1}\,{\bf q}_{r2}\,{\bf q}^{\dagger}_{d2}\,{\bf q}^{\dagger}_{d1}\Big\}\,\varepsilon \\
&=\Big\{\varrho^2_{r1}\,\varrho^2_{r2}\,+\, {\bf q}_{r1}\,{\bf q}_{r2}\,{\bf q}^{\dagger}_{d2}\,{\bf q}^{\dagger}_{d1}
\,+\, {\bf q}_{d1}\,{\bf q}_{d2}\,{\bf q}^{\dagger}_{r2}\,{\bf q}^{\dagger}_{r1} \,+\,\varrho^2_{d1}\,\varrho^2_{d2} \notag \\
&\,\;\;\;\;\;\;\;+\, \varrho^2_{r1}\,\varrho^2_{d2}\,+\, 
{\bf q}_{r1}\,{\bf q}_{d2}\,{\bf q}^{\dagger}_{r2}\,{\bf q}^{\dagger}_{d1} \,+\,{\bf q}_{d1}\,{\bf q}_{r2}\,{\bf q}^{\dagger}_{d2}\,{\bf q}^{\dagger}_{r1} \,+\,\varrho^2_{d1}\,\varrho^2_{r2}\Big\} \notag \\
&\;\;\;\;\;\;\;\;\;\;\;+\Big\{{\bf q}_{r1}\,{\bf q}_{r2}\,{\bf q}^{\dagger}_{d2}\,{\bf q}^{\dagger}_{r1} \,+\, \varrho^2_{r2}\,{\bf q}_{r1}\,{\bf q}^{\dagger}_{d1} \,+\, \varrho^2_{d2}\,{\bf q}_{d1}\,{\bf q}^{\dagger}_{r1} \,+\, {\bf q}_{d1}\,{\bf q}_{d2}\,{\bf q}^{\dagger}_{r2}\,{\bf q}^{\dagger}_{d1} \notag \\
&\,\;\;\;\;\;\;\;\;\;\;\;\;\;\;\;+\, {\bf q}_{r1}\,{\bf q}_{d2}\,{\bf q}^{\dagger}_{r2}\,{\bf q}^{\dagger}_{r1} \,+\, \varrho^2_{d2}\,{\bf q}_{r1}\,{\bf q}^{\dagger}_{d1} \,+\,\varrho^2_{r2}\,{\bf q}_{d1}\,{\bf q}^{\dagger}_{r1} \,+\, {\bf q}_{d1}\,{\bf q}_{r2}\,{\bf q}^{\dagger}_{d2}\,{\bf q}^{\dagger}_{d1}\Big\}\,\varepsilon \\
&=\left\{\varrho^2_{r1}\,\varrho^2_{r2}+\varrho^2_{r1}\,\varrho^2_{d2}+\varrho^2_{d1}\,\varrho^2_{r2}+\varrho^2_{d1}\,\varrho^2_{d2}+{\bf q}_{r1}\left({\bf q}_{r2}\,{\bf q}^{\dagger}_{d2}+{\bf q}_{d2}\,{\bf q}^{\dagger}_{r2}\right){\bf q}^{\dagger}_{d1}
\,+\, {\bf q}_{d1}\left({\bf q}_{r2}\,{\bf q}^{\dagger}_{d2}+
{\bf q}_{d2}\,{\bf q}^{\dagger}_{r2}\right){\bf q}^{\dagger}_{r1}\right\} \notag \\
&\;\;\;\;\;\;\;\;+\Big\{{\bf q}_{r1}\left({\bf q}_{r2}\,{\bf q}^{\dagger}_{d2}+{\bf q}_{d2}\,{\bf q}^{\dagger}_{r2}\right){\bf q}^{\dagger}_{r1} \,+\, \varrho^2_{r2}\left({\bf q}_{r1}\,{\bf q}^{\dagger}_{d1}+{\bf q}_{d1}\,{\bf q}^{\dagger}_{r1}\right) \notag \\
&\;\;\;\;\;\;\;\;\;\;\;\;\;\;\;\;\;\;+\varrho^2_{d2}\left({\bf q}_{r1}\,{\bf q}^{\dagger}_{d1}+{\bf q}_{d1}\,{\bf q}^{\dagger}_{r1}\right) + {\bf q}_{d1}\left({\bf q}_{r2}\,{\bf q}^{\dagger}_{d2}+{\bf q}_{d2}\,{\bf q}^{\dagger}_{r2}\right){\bf q}^{\dagger}_{d1}\Big\}\,\varepsilon \\
&=\left\{\varrho^2_{r1}\,\varrho^2_{r2}+\varrho^2_{r1}\,\varrho^2_{d2}+\varrho^2_{d1}\,\varrho^2_{r2}+\varrho^2_{d1}\,\varrho^2_{d2} \,+\,\left({\bf q}_{r1}\,{\bf q}^{\dagger}_{d1}+{\bf q}_{d1}\,{\bf q}^{\dagger}_{r1}\right)\left({\bf q}_{r2}\,{\bf q}^{\dagger}_{d2}+{\bf q}_{d2}\,{\bf q}^{\dagger}_{r2}\right)\right\} \notag \\
&\;\;\;\;\;\;\;\;\;\;+\left\{\left( \varrho^2_{r1}\,+\,\varrho^2_{d1}\right) \left({\bf q}_{r2}\,{\bf q}^{\dagger}_{d2}+{\bf q}_{d2}\,{\bf q}^{\dagger}_{r2}\right)
+\left( \varrho^2_{r2}\,+\,\varrho^2_{d2}\right) \left({\bf q}_{r1}\,{\bf q}^{\dagger}_{d1}+{\bf q}_{d1}\,{\bf q}^{\dagger}_{r1}\right)\right\}\,\varepsilon\,. \label{B8}
\end{align}
Again, since the quantities such as ${\bf q}_{r1}\,{\bf q}^{\dagger}_{d1}+{\bf q}_{d1}\,{\bf q}^{\dagger}_{r1}$ are scalar quantities, we see that the above product also resembles a split complex number similar to ${\mathbb Q}_{z}\,{\mathbb Q}^{\dagger}_{z}$ in (\ref{B1}). In other words, it is a sum of a scalar and a pseudoscalar, as in (\ref{23}).

More importantly, the right-hand sides of (\ref{B4}) and (\ref{B8}) are identical. We have thus proved that, although norms in ${{\cal K}^{\lambda}}$ resemble split complex numbers rather than scalars, the composition law continues to hold for the algebra ${{\cal K}^{\lambda}}$:
\begin{equation}
||{\mathbb Q}_{z1}\,{\mathbb Q}_{z2}||^2 = ||{\mathbb Q}_{z1}||^2\,||{\mathbb Q}_{z2}||^2. 
\end{equation}
Thus, the algebra ${\cal K}^{\lambda}$ is a normed division algebra even without assuming ${{\bf q}_{r}\,{\bf q}^{\dagger}_{d}+{\bf q}_{d}\,{\bf q}^{\dagger}_{r}=0}$. On the other hand, if we set ${{\bf q}_{r1}\,{\bf q}^{\dagger}_{d1}+{\bf q}_{d1}\,{\bf q}^{\dagger}_{r1}=0}$ and ${{\bf q}_{r2}\,{\bf q}^{\dagger}_{d2}+{\bf q}_{d2}\,{\bf q}^{\dagger}_{r2}=0}$ in (\ref{B4}) and (\ref{B8}) as a special case [as we did in proving the relation (\ref{norm-final})], then the above composition law reduces to the norm relation (\ref{norm-final}), confirming our main thesis above:
\begin{equation}
||{\mathbb Q}_{z1}\,{\mathbb Q}_{z2}|| \,=\, \sqrt{\varrho^2_{r1}\,\varrho^2_{r2}+\varrho^2_{r1}\,\varrho^2_{d2}+\varrho^2_{d1}\,\varrho^2_{r2}+\varrho^2_{d1}\,\varrho^2_{d2}\,} \,=\, ||{\mathbb Q}_{z1}||\;||{\mathbb Q}_{z2}||\,. 
\end{equation}
The advantage of the special case is that it reduces the norms from a split complex form to purely scalar quantities. 

\section{Orthogonality of the Quaternions ${{\bf q}_{r}}$ and ${{\bf q}_{d}}$ is Preserved under Multiplication} \label{C}

In this appendix, we prove that the 7-sphere we have constructed in this paper and defined in Eq.~(\ref{spring}), namely
\begin{equation}
{\cal K}^{\lambda}\supset\,S^7:=\,\left\{\,{\mathbb Q}_z :=\, {\bf q}_r + {\bf q}_d\,\varepsilon\;\Big|\; ||{\mathbb Q}_z||=1\;\,\text{and}\,\;{\bf q}_{r}\,{\bf q}^{\dagger}_{d}+{\bf q}_{d}\,{\bf q}^{\dagger}_{r}=0\,\right\}, \label{Cspring}
\end{equation}
remains closed under multiplication. This may not be obvious because of the orthogonality condition ${\bf q}_{r}\,{\bf q}^{\dagger}_{d}+{\bf q}_{d}\,{\bf q}^{\dagger}_{r}=0$ we have imposed for normalizing the elements ${\mathbb Q}_z$ of the algebra ${\cal K}^{\lambda}$. But the orthogonality of the quaternions ${{\bf q}_{r}}$ and ${{\bf q}_{d}}$ in ${\mathbb Q}_z$ is preserved under multiplication of the elements of ${\cal K}^{\lambda}$. To prove this, consider two distinct elements of ${{\cal K}^{\lambda}}$,
\begin{equation}
{\mathbb Q}_{z1} =\, {\bf q}_{r1} + {\bf q}_{d1}\,\varepsilon\;\;\;\text{and}\;\;\;{\mathbb Q}_{z2} =\, {\bf q}_{r2} + {\bf q}_{d2}\,\varepsilon\,,
\end{equation}
satisfying the orthogonality conditions 
${\bf q}_{r1}\,{\bf q}^{\dagger}_{d1}+{\bf q}_{d1}\,{\bf q}^{\dagger}_{r1}=0$ and ${\bf q}_{r2}\,{\bf q}^{\dagger}_{d2}+{\bf q}_{d2}\,{\bf q}^{\dagger}_{r2}=0$, together with their product
\begin{equation}
{\mathbb Q}_{z3}:=
{\mathbb Q}_{z1}\,{\mathbb Q}_{z2}\,=\,\left({\bf q}_{r1}\,{\bf q}_{r2}\,+\,{\bf q}_{d1}\,{\bf q}_{d2}\right)\,+\,\left({\bf q}_{r1}\,{\bf q}_{d2}\,+\,{\bf q}_{d1}\,{\bf q}_{r2}\right)\,\varepsilon\,, \label{c3}
\end{equation}
as considered in Eq.~(\ref{new24}). If we now define the quaternions
\begin{equation}
{\bf q}_{r3}\,:=\,\left({\bf q}_{r1}\,{\bf q}_{r2}\,+\,{\bf q}_{d1}\,{\bf q}_{d2}\right)\;\;\;\text{and}\;\;\;
{\bf q}_{d3}\,:=\,\left({\bf q}_{r1}\,{\bf q}_{d2}\,+\,{\bf q}_{d1}\,{\bf q}_{r2}\right),
\end{equation}
then the product (\ref{c3}) can be expressed as the following third element of ${\cal K}^{\lambda}$: 
\begin{equation}
{\mathbb Q}_{z3}=\,{\bf q}_{r3} + {\bf q}_{d3}\,\varepsilon\,.
\label{notso}
\end{equation}
It is now easy to prove that the quaternions ${{\bf q}_{r3}}$ and ${{\bf q}_{d3}}$ are also orthogonal, or equivalently, ${{\bf q}_{r3}\,{\bf q}^{\dagger}_{d3}+{\bf q}_{d3}\,{\bf q}^{\dagger}_{r3}=0}$:
\begin{align}
{\bf q}_{r3}\,{\bf q}^{\dagger}_{d3}+{\bf q}_{d3}\,{\bf q}^{\dagger}_{r3}\,&=\,
\left({\bf q}_{r1}\,{\bf q}_{r2}\,+\,{\bf q}_{d1}\,{\bf q}_{d2}\right)
\left({\bf q}_{r1}\,{\bf q}_{d2}\,+\,{\bf q}_{d1}\,{\bf q}_{r2}\right)^{\dagger}+\,
\left({\bf q}_{r1}\,{\bf q}_{d2}\,+\,{\bf q}_{d1}\,{\bf q}_{r2}\right)
\left({\bf q}_{r1}\,{\bf q}_{r2}\,+\,{\bf q}_{d1}\,{\bf q}_{d2}\right)^{\dagger} \\
&=\,
\left({\bf q}_{r1}\,{\bf q}_{r2}\,+\,{\bf q}_{d1}\,{\bf q}_{d2}\right)
\left({\bf q}_{d2}^{\dagger}\,{\bf q}_{r1}^{\dagger}\,+\,{\bf q}_{r2}^{\dagger}\,{\bf q}_{d1}^{\dagger}\right)+\,
\left({\bf q}_{r1}\,{\bf q}_{d2}\,+\,{\bf q}_{d1}\,{\bf q}_{r2}\right)
\left({\bf q}_{r2}^{\dagger}\,{\bf q}_{r1}^{\dagger}\,+\,{\bf q}_{d2}^{\dagger}\,{\bf q}_{d1}^{\dagger}\right) \\
&=\,
{\bf q}_{r1}\,{\bf q}_{r2}\,{\bf q}_{d2}^{\dagger}\,{\bf q}_{r1}^{\dagger}\,+\,{\bf q}_{r1}\,{\bf q}_{r2}\,{\bf q}_{r2}^{\dagger}\,{\bf q}_{d1}^{\dagger}\,+\,{\bf q}_{d1}\,{\bf q}_{d2}\,{\bf q}_{d2}^{\dagger}\,{\bf q}_{r1}^{\dagger}\,+\,{\bf q}_{d1}\,{\bf q}_{d2}\,{\bf q}_{r2}^{\dagger}\,{\bf q}_{d1}^{\dagger} \notag \\
&\;\;\;\;\;\;\;\;+
\,{\bf q}_{r1}\,{\bf q}_{d2}\,{\bf q}_{r2}^{\dagger}\,{\bf q}_{r1}^{\dagger}\,+\,{\bf q}_{r1}\,{\bf q}_{d2}\,{\bf q}_{d2}^{\dagger}\,{\bf q}_{d1}^{\dagger}\,+\,{\bf q}_{d1}\,{\bf q}_{r2}\,{\bf q}_{r2}^{\dagger}\,{\bf q}_{r1}^{\dagger}\,+\,{\bf q}_{d1}\,{\bf q}_{r2}\,{\bf q}_{d2}^{\dagger}\,{\bf q}_{d1}^{\dagger} \\
&=\,{\bf q}_{r1}\left\{{\bf q}_{r2}\,{\bf q}_{d2}^{\dagger}+{\bf q}_{d2}\,{\bf q}_{r2}^{\dagger}\right\}{\bf q}_{r1}^{\dagger}\,+\,\varrho_{r2}^2\left\{{\bf q}_{r1}\,{\bf q}_{d1}^{\dagger}+{\bf q}_{d1}\,{\bf q}_{r1}^{\dagger}\right\} \notag \\
&\;\;\;\;\;\;\;\;\;\;\;\;\;\;\;\;\;\;\;\;+\,{\bf q}_{d1}\left\{{\bf q}_{r2}\,{\bf q}_{d2}^{\dagger}+{\bf q}_{d2}\,{\bf q}_{r2}^{\dagger}\right\}{\bf q}_{d1}^{\dagger}\,+\,\varrho_{d2}^2\left\{{\bf q}_{r1}\,{\bf q}_{d1}^{\dagger}+{\bf q}_{d1}\,{\bf q}_{r1}^{\dagger}\right\} \label{c9}\\ 
&=\,\left(\varrho_{r1}^2+\varrho_{d1}^2\right)\left\{{\bf q}_{r2}\,{\bf q}_{d2}^{\dagger}+{\bf q}_{d2}\,{\bf q}_{r2}^{\dagger}\right\}\,+\,\left(\varrho_{r2}^2+\varrho_{d2}^2\right)\left\{{\bf q}_{r1}\,{\bf q}_{d1}^{\dagger}+{\bf q}_{d1}\,{\bf q}_{r1}^{\dagger}\right\} \label{c10}\\
&=\,0.
\end{align}
Here Eq.~(\ref{c10}) follows from Eq.~(\ref{c9}) because, as shown in Eq.~(\ref{22}), the quantities such as ${\bf q}_{r1}\,{\bf q}_{d1}^{\dagger}+{\bf q}_{d1}\,{\bf q}_{r1}^{\dagger}$ are scalar quantities, and therefore we can use 
${\bf q}_{r1}\,{\bf q}^{\dagger}_{r1}=\varrho^2_{r1}$, ${\bf q}_{d1}\,{\bf q}^{\dagger}_{d1}=\varrho^2_{d1}$, {\it etc}., as before, to reduce Eq.~(\ref{c9}) to
Eq.~(\ref{c10}). And the orthogonality conditions ${\bf q}_{r1}\,{\bf q}^{\dagger}_{d1}+{\bf q}_{d1}\,{\bf q}^{\dagger}_{r1}=0$ and ${\bf q}_{r2}\,{\bf q}^{\dagger}_{d2}+{\bf q}_{d2}\,{\bf q}^{\dagger}_{r2}=0$ then reduces the RHS of Eq.~(\ref{c10}) to zero.\break Consequently, the 7-sphere defined in (\ref{spring}) remains closed under multiplication, analogously to the octonionic 7-sphere.

\section{The Even Subalgebra $\,{\cal K}^{\lambda}$ of the Clifford algebra ${\mathrm{Cl}_{4,0}}$ is a Normed Division Algebra} \label{D}

In this appendix, I demonstrate why ${\cal K}^{\lambda}$ is an eight-dimensional normed division algebra despite the constructibility in it of multivectors such as ${\bf X}=\alpha\,(1+\,\varepsilon)$ and ${\bf Y}=\gamma\,(1-\,\varepsilon)$ with any real numbers $\alpha$ and $\gamma$ that lead to the vanishing of their geometric product, ${\bf X}{\bf Y}=0$, because the pseudoscalar $\varepsilon$ satisfies the property $\varepsilon^2=1$. Such mixed grade multivectors in ${\cal K}^{\lambda}$ may give the wrong impression of the appearance of ``nonzero zero divisors'' within the algebra, and lead one to suspect that perhaps ${\cal K}^{\lambda}$ is not a normed division algebra. However, it is not difficult to recognize that multivectors such as ${\bf X}=\alpha\,(1+\,\varepsilon)$ and ${\bf Y}=\gamma\,(1-\,\varepsilon)$ {\it cannot be normalized} to non-vanishing scalar magnitudes using any consistently implementable criterion for normalization in Geometric Algebra that does not violate the multiplicative closure of ${\cal K}^{\lambda}$ while generating new multivectors within the algebra, and therefore such multivectors are not a part of the {\it normed} division algebra ${\cal K}^{\lambda}$. By definition, a normed division algebra is a division algebra in which each of its members can be normalized, or assigned a scalar magnitude. While a detailed analysis of why the algebra ${\cal K}^{\lambda}$ is a normed division algebra is already published in \cite{RSOS,RSOS-2,Local} and discussed in the present paper, it is convenient to have the results from the first three physics-oriented publications presented coherently in one place, together with the results discussed in this paper. To this end, my goal in this appendix is to prove the following statements anew:
\begin{enumerate}[label=(\arabic*)]
    \item ${\cal K}^{\lambda}$ is eight-dimensional
    \item ${\cal K}^{\lambda}$ is a genuine algebra
    \item ${\cal K}^{\lambda}$ is an associative algebra
    \item ${\cal K}^{\lambda}$ is a normed algebra
    \item ${\cal K}^{\lambda}$ is a division algebra
\end{enumerate}
The first three properties in this list do not require special proofs. By definition, ${\cal K}^{\lambda}$ is an eight-dimensional even subalgebra of the $2^4=16$-dimensional associative Clifford algebra $\mathrm{Cl}_{4,0}$. As a subalgebra of the algebra $\mathrm{Cl}_{4,0}$, it is thus a {\it bona fide} associative algebra that remains closed under multiplication and addition, and possesses multiplicative and additive identities. The eight even dimensions of ${\cal K}^{\lambda}$ are also evident from the multiplication table~\ref{T+1}. Using that table, a general multivector ${\mathbb Q}_z$ in the graded basis of ${\cal K}^{\lambda}$ with appropriate notation for coefficients can be expanded~as
\begin{equation}
{\mathbb Q}_z=\,u^{\,}_0+u_z\,{\bf e}_x{\bf e}_y+u_y\,{\bf e}_z{\bf e}_x+u_x\,{\bf e}_y{\bf e}_z+v_x\,{\bf e}_x{\bf e}_{\infty}+v_y\,{\bf e}_y{\bf e}_{\infty}+v_z\,{\bf e}_z{\bf e}_{\infty}-v^{\,}_0\,I_3{\bf e}_{\infty}\,,\label{qinqq}
\end{equation}
where I have fixed the orientation of the vector space to $\lambda=+1$ for convenience,
so that it can also be expressed as
\begin{equation}
{\mathbb Q}_z =\,{\bf q}_r + {\bf q}_d\,\varepsilon\,, \label{AAA-1}
\end{equation}
with pseudoscalar defined as $\varepsilon:=-I_3{\bf e}_{\infty}$ with properties $\varepsilon^2=1$ and $\varepsilon^{\dagger}=\varepsilon$, in terms of the following two quaternions:
\begin{equation}
{\bf q}_{r}:=u_0+I_3{\bf u}\;\;\;\;\text{and}\;\;\;\;{\bf q}_{d}:=v_0+I_3{\bf v}\,.
\end{equation}
Here $u^{}_0$ and $v^{}_0$ are scalars, ${I_3={\bf e}_x{\bf e}_y{\bf e}_z}$ with the property $I_3^{\dagger}=-I_3$ is a unit pseudoscalar in three dimensions, and ${{\bf u}=u_{x}{\bf e}_x+u_{y}{\bf e}_y+u_{z}{\bf e}_z}$ and ${{\bf v}=v_{x}{\bf e}_x+v_{y}{\bf e}_y+v_{z}{\bf e}_z}$ are vectors in $\mathrm{I\!R}^3$. It is evident from (\ref{AAA-1}) that the algebra ${\cal K}^{\lambda}$ is a tensor product ${\mathbb H}\otimes\mathbb{C}'$ of quaternions with quantities such as $\alpha + \gamma\;\varepsilon\in\mathbb{C}'$ that resemble split complex numbers, with conjugation (or reverse) affecting only the quaternions because of the property $\varepsilon^{\dagger}=\varepsilon$, where $\alpha$ and $\gamma$ are scalars.

\renewcommand{\thesubsection}{D.\arabic{subsection}}

\subsection{What must be demonstrated to prove the properties (4) and (5) in the above list:}

The last two properties in the above list are contentious due to the constructibility of the so-called ``non-zero zero divisors'' in the {\bf un}-normed algebra ${\cal K}^{\lambda}$. We therefore devote the rest of this appendix to proving the last two listed properties. To that end, let me first explain what is meant by a ``normed algebra'' in the present context. Recall that mathematicians have long wondered: When is a product of two squares of real numbers itself a square: ${x^2\,y^2 = z^2}$? Any real number such as ${z}$ can be written as a product of two other real numbers, ${z = x\,y}$, if it is so factorizable; and then the factorization the squares of these numbers, ${z^2 = x^2\,y^2}$, sometimes also holds for the numbers ${x}$, ${y}$, and ${z}$. But what about the {\it sums} of the squares of real numbers? Historically, this question goes back to Euler, Hamilton, Degen, Graves, and Cayley, leading Hurwitz \cite{Hurwitz-1898} to prove in 1898 that such a factorization is possible only for $n=1$, 2, 4, and 8 real numbers, so that we can indeed decompose the sum of the squares of such numbers, say $z_1^2 + z_2^2 + z_3^2 + \dots + z_n^2$,~as
\begin{equation}
\left(z_1^2 + z_2^2 + z_3^2 + z_4^2 + z_5^2 + \dots + z_n^2\right) \,=\,
\left(x_1^2 + x_2^2 + x_3^2 + x_4^2 + x_5^2 + \dots + x_n^2\right)\left(y_1^2 + y_2^2 + y_3^2 + y_4^2 + y_5^2 + \dots + y_n^2\right) \label{hurhur}
\end{equation}
for $n=1$, 2, 4, or 8. We are interested in the case $n=8$ and in proving the following equality within our algebra ${\cal K}^{\lambda}$:
\begin{equation}
\left(z_1^2 + z_2^2 + z_3^2 + \dots + z_8^2\right) \,=\,
\left(x_1^2 + x_2^2 + x_3^2 + \dots + x_8^2\right)\left(y_1^2 + y_2^2 + y_3^2 + \dots + y_8^2\right), \label{sd5}
\end{equation}
If we now christen each of these three sums of squares of eight real numbers, respectively, by $\varrho^2_{\!_{\bf Z}}$, $\varrho^2_{\!_{\bf X}}$, and $\varrho^2_{\!_{\bf Y}}$ as follows,
\begin{align}
z_1^2 + z_2^2 + z_3^2 + z_4^2 + z_5^2 + z_6^2 + z_7^2 + z_8^2\;\, &= \,\varrho^2_{\!_{\bf Z}}\,, \\
x_1^2 + x_2^2 + x_3^2 + x_4^2 + x_5^2 + x_6^2 + x_7^2 + x_8^2 &= \,\varrho^2_{\!_{\bf X}}\,, \\
y_1^2 + y_2^2 + y_3^2 + y_4^2 + y_5^2 + y_6^2 + y_7^2 + y_8^2\; &= \,\varrho^2_{\!_{\bf Y}}\,,
\end{align}
then we can interpret these quadratic equations as defining equations, respectively, of three possible radii $\varrho^{}_{\!_{\bf Z}}$, $\varrho^{}_{\!_{\bf X}}$, and $\varrho^{}_{\!_{\bf Y}}$\break of a 7-sphere. Then, the equality (\ref{sd5}) among sums of squared numbers gives the following relation among these radii:
\begin{equation}
\varrho^2_{\!_{\bf Z}}=\varrho^2_{\!_{\bf X}}\,\varrho^2_{\!_{\bf Y}},
\end{equation}
which is equivalent to
\begin{equation}
\varrho_{\!_{\bf Z}}=\varrho_{\!_{\bf X}}\,\varrho_{\!_{\bf Y}}.
\end{equation}
The decomposition of the radius $\varrho_{\!_{\bf Z}}$ is thus multiplicative. As a result, we can rewrite the equality (\ref{sd5}) equivalently~as
\begin{equation}
\Bigg(\sqrt{z_1^2 + z_2^2 + z_3^2 + \dots + z_8^2\,}\;\Bigg)
=\Bigg(\sqrt{x_1^2 + x_2^2 + x_3^2 + \dots + x_8^2\,}\;\Bigg)\,\Bigg(\sqrt{y_1^2 + y_2^2 + y_3^2 + \dots + y_8^2\,}\;\Bigg). \label{eq-eq}
\end{equation}

Note that since each quantity in the parentheses in the above equality defines a radius of a 7-sphere, or $S^7$, it is the only possible spherical surface that can be embedded in the eight dimensions of $\mathrm{I\!R^8}$ whose radius can satisfy this equality. This also follows from Hurwitz's theorem \cite{Hurwitz-1898}, which states that only the radii of the spheres $S^0$, $S^1$, $S^3$, and $S^7$ can satisfy the equality (\ref{hurhur}), and therefore the above equality can be satisfied only by the radii of $S^7$. If we now think of the real numbers $z_1$, $z_2$, $z_3$, $\dots$, $z_8$ as coefficients of a vector ${\mathbf{Z}}$ in $\mathrm{I\!R}^8$, and analogously for the vectors ${\mathbf{X}}$ and ${\mathbf{Y}}$, and define their norms as
\begin{align}
\varrho_{\!_{\bf Z}}\;=||\mathbf{Z}|| &:= \;\sqrt{z_1^2 + z_2^2 + z_3^2 + z_4^2 + z_5^2 + z_6^2 + z_7^2 + z_8^2\,}\,, \\
\varrho_{\!_{\bf X}}=||\mathbf{X}|| &:= \sqrt{x_1^2 + x_2^2 + x_3^2 + x_4^2 + x_5^2 + x_6^2 + x_7^2 + x_8^2\,}\,, \\
\varrho_{\!_{\bf Y}}=||\mathbf{Y}|| &:= \;\sqrt{y_1^2 + y_2^2 + y_3^2 + y_4^2 + y_5^2 + y_6^2 + y_7^2 + y_8^2\,}\,,
\end{align}
then the relation (\ref{eq-eq}) we seek to prove for the even subalgebra ${{\cal K}^{\lambda}}$ of the algebra ${\mathrm{Cl}_{4,0}}$ can be expressed simply as  
\begin{equation}
||\mathbf{Z}||=||\mathbf{X}||\,||\mathbf{Y}||\,. \label{relD}
\end{equation}
Moreover, as noted above, ${{\cal K}^{\lambda}}$ remains closed under multiplication---a geometric product between two multivectors in ${{\cal K}^{\lambda}}$, say ${\bf X}$ and ${\bf Y}$, yields another multivector ${{\bf Z}={\bf X}{\bf Y}}$ in ${{\cal K}^{\lambda}}$. In terms of graded basis, ${\bf X}$ and ${\bf Y}$ can be expanded as 
\begin{equation}
{\bf X}=x_1+x_2\,{\bf e}_x{\bf e}_y+x_3\,{\bf e}_z{\bf e}_x+x_4\,{\bf e}_y{\bf e}_z+x_5\,{\bf e}_x{\bf e}_{\infty}+x_6\,{\bf e}_y{\bf e}_{\infty}+x_7\,{\bf e}_z{\bf e}_{\infty}+x_8\,\varepsilon \label{xD}
\end{equation}
and
\begin{equation}
{\bf Y}=\;y_1+y_2\,{\bf e}_x{\bf e}_y+y_3\,{\bf e}_z{\bf e}_x+y_4\,{\bf e}_y{\bf e}_z+y_5\,{\bf e}_x{\bf e}_{\infty}+y_6\,{\bf e}_y{\bf e}_{\infty}+y_7\,{\bf e}_z{\bf e}_{\infty}+y_8\,\varepsilon\,.\label{yD}
\end{equation}
It is then straightforward to verify using the multiplication table for ${\cal K}^{\lambda}$ (Table \ref{T+1}) that the geometric product ${{\bf Z}={\bf X}{\bf Y}}$ of the two multivectors ${\bf X}$ and ${\bf Y}$ is yet another multivector ${\bf Z}$ in ${\cal K}^{\lambda}$, which too can be expanded in the basis of ${\cal K}^{\lambda}$ as 
\begin{equation}
{\bf Z}=\,z_1+z_2\,{\bf e}_x{\bf e}_y+z_3\,{\bf e}_z{\bf e}_x+z_4\,{\bf e}_y{\bf e}_z+z_5\,{\bf e}_x{\bf e}_{\infty}+z_6\,{\bf e}_y{\bf e}_{\infty}+z_7\,{\bf e}_z{\bf e}_{\infty}+z_8\,\varepsilon\,. \label{zD}
\end{equation}
It is important to note that only a geometric product between multivectors in ${\cal K}^{\lambda}$ can give a new multivector in ${\cal K}^{\lambda}$, because geometric product, as defined in (\ref{gp}), is the only fundamental product in Geometric Algebra, unlike (\ref{sim-gp}) or (\ref{ant-gp}).\break This property of ${\cal K}^{\lambda}$ allows us to express the relation (\ref{relD}) we seek to prove here as the following composition law:
\begin{equation}
||{\bf X}{\bf Y}|| = ||{\bf X}||\;||{\bf Y}||\,. \label{comnorm}
\end{equation}
Thus, to prove that algebra ${\cal K}^{\lambda}$ is a {\it normed} algebra, we must prove that the scalar-valued norm $||{\bf X}{\bf Y}||$ is multiplicative.

\subsection{Proof of the Multivector-valued Composition Law $\,\widetilde{||{\bf X}{\bf Y}||} = \widetilde{||{\bf X}||}\;\widetilde{||{\bf Y}||}$}

We now proceed to prove the composition law (\ref{comnorm}) and the equivalent equality (\ref{eq-eq}) in two stages. The result we are about to prove in this subsection is the most general multivector-valued composition law possible for the algebra ${\cal K}^{\lambda}$. To that end, for any multivector ${\bf X}$ within ${\cal K}^{\lambda}$, we defined a multivector-valued norm using geometric product~as
\begin{equation}
\widetilde{||{\bf X}||}:=\sqrt{{\bf X}{\bf X}^{\dagger}}, \label{mulnorm}
\end{equation}
where ${\bf X}^{\dagger}$ is the reverse (or conjugate) of the multivector ${\bf X}$. This definition is consistent with how norms are defined for complex numbers $c$, quaternions ${\bf q}$, and octonions ${\bf O}$, as $||{c}||=\sqrt{{c}\,{c}^{\dagger}}$, $||{\bf q}||=\sqrt{{\bf q}{\bf q}^{\dagger}}$, and $||{\bf O}||=\sqrt{{\bf O}{\bf O}^{\dagger}}$, respectively. Moreover, a geometric product between ${\bf X}$ and ${\bf X}^{\dagger}$ in the above definition is necessary for maintaining consistency between the two sides of the composition law
\begin{equation}
\widetilde{||{\bf X}{\bf Y}||} = \widetilde{||{\bf X}||}\;\widetilde{||{\bf Y}||}\,, \label{mulmul}
\end{equation}
because, as noted above, only the geometric product, and not a scalar or a wedge product between the multivectors in ${\cal K}^{\lambda}$ can produce a new multivector ${\bf Z}={\bf X}{\bf Y}$ in ${\cal K}^{\lambda}$ appearing on the left-hand side of the composition law (\ref{mulmul}). On the other hand, because the algebra ${\cal K}^{\lambda}$ is a tensor product ${\mathbb H}\otimes\mathbb{C}'$ of quaternions with quantities such as $a+ b\,\varepsilon$, and thus resemble split complex numbers, its elements are of the form ${\mathbb Q}_z=\, {\bf q}_{r} + {\bf q}_{d}\,\varepsilon$ with $\varepsilon^2=1$ and $\varepsilon^{\dagger}=\varepsilon$, and consequently the quadratic form ${{\mathbb Q}_z{\mathbb Q}_z^{\dagger}}$ in general takes values in the split complex numbers $\mathbb{C}'$ instead of real numbers. 

Elsewhere \cite{RSOS,RSOS-2,Local}, I have proved that the multivector-valued composition law (\ref{mulmul}) holds, without exception, for arbitrary multivectors ${\bf X}$ and ${\bf Y}$ in ${\cal K}^{\lambda}$, and, as a special case, reduces to the composition law with scalar values for $||{\bf X}||$, $||{\bf Y}||$, and $||{\bf X}{\bf Y}||$; namely, to (\ref{comnorm}). The proof of (\ref{mulmul}) is straightforward. As noted above, the algebra ${\cal K}^{\lambda}$ is a tensor product $\mathbb{H}\otimes\mathbb{C}'$ of quaternions with split complex numbers \cite{Dray}, but with conjugation (or ``reverse") affecting only the quaternions, because $\varepsilon^{\dagger}=\varepsilon$. In other words, any multivectors ${\bf X}$ and ${\bf Y}$ within ${\cal K}^{\lambda}$ are of the form 
\begin{equation}
{\bf X}=\, {\bf q}_{r1} + {\bf q}_{d1}\,\varepsilon\;\;\;\;\text{and}\;\;\;\;
{\bf Y}=\, {\bf q}_{r2} + {\bf q}_{d2}\,\varepsilon\,, \label{def-XY}
\end{equation}
with $\varepsilon^2=1$ and $\varepsilon^{\dagger}=\varepsilon$, where ${\bf q}_{r1}$ and ${\bf q}_{d1}$ constituting ${\bf X}$, for example, are two independent quaternions, which can be expressed as ${\bf q}_{r1}=g_1+I_3{\bf u}_1$ and ${\bf q}_{d1}=h_1+I_3{\bf v}_1$, with $g_1$ and $h_1$ as scalars, ${I_3={\bf e}_x{\bf e}_y{\bf e}_z}$, with $I_3^{\dagger}=-I_3$, as the standard pseudoscalar or trivector in three dimensions, and ${{\bf u}_1=u_{1x}{\bf e}_x+u_{1y}{\bf e}_y+u_{1z}{\bf e}_z}$ and ${{\bf v}_1=v_{1x}{\bf e}_x+v_{1y}{\bf e}_y+v_{1z}{\bf e}_z}$ as vectors in $\mathrm{I\!R}^3$. As a result, the geometric product ${\bf X}{\bf X}^{\dagger}$ between the multivector ${\bf X}$ and its reverse ${\bf X}^{\dagger}$ works out~to~be
\begin{align}
{\bf X}{\bf X}^{\dagger}\,&=\left({\bf q}_{r1} + {\bf q}_{d1}\,\varepsilon\right)\left({\bf q}_{r1} + {\bf q}_{d1}\,\varepsilon\right)^{\dagger}
\label{sd10a}\\
&=\left({\bf q}_{r1}\,{\bf q}^{\dagger}_{r1}\,+\,{\bf q}_{d1}\,{\bf q}^{\dagger}_{d1}\right)+\left({\bf q}_{r1}\,{\bf q}^{\dagger}_{d1}\,+\,{\bf q}_{d1}\,{\bf q}^{\dagger}_{r1}\right)\varepsilon \label{sd10b} \\
&=(g_1^2 + {\bf u}_1\cdot{\bf u}_1 + h_1^2 + {\bf v}_1\cdot{\bf v}_1) + ( 2\,g_1h_1 + 2\,{\bf u}_1\cdot{\bf v}_1)\,\varepsilon \label{sd10bbb} \\
&=\,\text{(a scalar)}+\text{(a scalar)}\;\varepsilon\,.
\end{align}
And likewise, using ${\bf Y}=\, {\bf q}_{r2} + {\bf q}_{d2}\,\varepsilon$, with ${\bf q}_{r2}=g_2+I_3{\bf u}_2$ and ${\bf q}_{d2}=h_2+I_3{\bf v}_2$, the geometric product ${\bf Y}{\bf Y}^{\dagger}$ between ${\bf Y}$ and ${\bf Y}^{\dagger}$ works out to be
\begin{align}
{\bf Y}{\bf Y}^{\dagger}\,
&=\left({\bf q}_{r2} + {\bf q}_{d2}\,\varepsilon\right)\left({\bf q}_{r2} + {\bf q}_{d2}\,\varepsilon\right)^{\dagger}
\label{sd10ab}\\
&=\left({\bf q}_{r2}\,{\bf q}^{\dagger}_{r2}\,+\,{\bf q}_{d2}\,{\bf q}^{\dagger}_{d2}\right)+\left({\bf q}_{r2}\,{\bf q}^{\dagger}_{d2}\,+\,{\bf q}_{d2}\,{\bf q}^{\dagger}_{r2}\right)\varepsilon \label{sd10bc} \\
&=(g_2^2 + {\bf u}_2\cdot{\bf u}_2 + h_2^2 + {\bf v}_2\cdot{\bf v}_2) + ( 2\,g_2h_2 + 2\,{\bf u}_2\cdot{\bf v}_2)\,\varepsilon \label{sd10bcd} \\
&=\,\text{(a scalar)}+\text{(a scalar)}\;\varepsilon\,. \label{sd30}
\end{align}

\allowdisplaybreaks[0]

The proof of the multivector-valued composition law (\ref{mulmul}) is now straightforward. Evidently, the geometric products ${\bf X}{\bf X}^{\dagger}$ and ${\bf Y}{\bf Y}^{\dagger}$ resemble split complex numbers because for any two quaternions (such as ${\bf q}_{r1}$ and ${\bf q}_{d1}$) the quantities appearing in the parentheses in (\ref{sd10bbb}) and (\ref{sd10bcd}) are scalar quantities. Thus the products are of the form ${\bf X}{\bf X}^{\dagger}=a+b\,\varepsilon$ and ${\bf Y}{\bf Y}^{\dagger}=c+d\,\varepsilon$, where $a={\bf q}_{r1}\,{\bf q}^{\dagger}_{r1}+{\bf q}_{d1}\,{\bf q}^{\dagger}_{d1}$, $b={\bf q}_{r1}\,{\bf q}^{\dagger}_{d1}+{\bf q}_{d1}\,{\bf q}^{\dagger}_{r1}$, $c={\bf q}_{r2}\,{\bf q}^{\dagger}_{r2}+{\bf q}_{d2}\,{\bf q}^{\dagger}_{d2}$, and $d={\bf q}_{r2}\,{\bf q}^{\dagger}_{d2}+{\bf q}_{d2}\,{\bf q}^{\dagger}_{r2}$ are scalar quantities, with $\varepsilon^{\dagger}=\varepsilon$ instead of $-\,\varepsilon$. Hence, for the left-hand side of (\ref{mulmul}) we have
\begin{align}
\widetilde{||{\bf X}{\bf Y}||}&=\sqrt{({\bf X}{\bf Y})({\bf X}{\bf Y})^{\dagger}} \\
&=\sqrt{({\bf X}{\bf Y})\left({\bf Y}^{\dagger}{\bf X}^{\dagger}\right)} \\
&=\sqrt{{\bf X}\left({\bf Y}{\bf Y}^{\dagger}\right){\bf X}^{\dagger}} \\
&=\sqrt{{\bf X}(c+d\,\varepsilon){\bf X}^{\dagger}} \\
&=\sqrt{\left({\bf X}{\bf X}^{\dagger}\right)(c+d\,\varepsilon)} \\
&=\sqrt{(a+b\,\varepsilon)(c+d\,\varepsilon)} \\
&=\sqrt{(ac+bd)+(ad+bc)\,\varepsilon\,}, \label{sdLH}
\end{align}
because the pseudoscalar $\varepsilon$ satisfies $\varepsilon^2=1$ and commutes with every element of ${\cal K}^{\lambda}$, and consequently the identity $({\bf X}{\bf Y})^{\dagger}=\left({\bf Y}^{\dagger}{\bf X}^{\dagger}\right)$ for ${\bf X}$ and ${\bf Y}$ in ${\cal K}^{\lambda}$ is straightforward to verify, using the associativity of the algebra. On the other hand, the right-hand side of the multivector-valued composition law (\ref{mulmul}) also works out to give the same quantity:
\begin{equation}
\widetilde{||{\bf X}||}\;\widetilde{||{\bf Y}||}=\left(\sqrt{{\bf X}{\bf X}^{\dagger}}\;\right)\left(\sqrt{{\bf Y}{\bf Y}^{\dagger}}\;\right)=\sqrt{(a+b\,\varepsilon)(c+d\,\varepsilon)}=\sqrt{(ac+bd)+(ad+bc)\,\varepsilon\,}. \label{sdRH}
\end{equation}
Comparing the results (\ref{sdLH}) and (\ref{sdRH}), we see that both sides of the equation (\ref{mulmul}) have worked out to be identical to the square-root of the following multivector quantity, which also resembles a split complex (or hyperbolic) number, 
\begin{align}
\widetilde{||{\bf X}{\bf Y}||^2}
&=\left\{\left(\varrho^2_{r1}+\varrho^2_{d1}\right)\left(\varrho^2_{r2}+\varrho^2_{d2}\right) \,+\,\left({\bf q}_{r1}\,{\bf q}^{\dagger}_{d1}+{\bf q}_{d1}\,{\bf q}^{\dagger}_{r1}\right)\left({\bf q}_{r2}\,{\bf q}^{\dagger}_{d2}+{\bf q}_{d2}\,{\bf q}^{\dagger}_{r2}\right)\right\} \notag \\
&\;\;\;\;\;\;\;\;\;\;+\left\{\left( \varrho^2_{r1}\,+\,\varrho^2_{d1}\right) \left({\bf q}_{r2}\,{\bf q}^{\dagger}_{d2}+{\bf q}_{d2}\,{\bf q}^{\dagger}_{r2}\right)
+ \left({\bf q}_{r1}\,{\bf q}^{\dagger}_{d1}+{\bf q}_{d1}\,{\bf q}^{\dagger}_{r1}\right)\left( \varrho^2_{r2}\,+\,\varrho^2_{d2}\right)\right\}\varepsilon \\
&=\left\{\,\varrho^2_{\!_{\bf X}}\;\varrho^2_{\!_{\bf Y}}\,+\,\left({\bf q}_{r1}\,{\bf q}^{\dagger}_{d1}+{\bf q}_{d1}\,{\bf q}^{\dagger}_{r1}\right)\left({\bf q}_{r2}\,{\bf q}^{\dagger}_{d2}+{\bf q}_{d2}\,{\bf q}^{\dagger}_{r2}\right)\right\} \notag \\
&\;\;\;\;\;\;\;\;\;\;+\left\{\,\varrho^2_{\!_{\bf X}} \left({\bf q}_{r2}\,{\bf q}^{\dagger}_{d2}+{\bf q}_{d2}\,{\bf q}^{\dagger}_{r2}\right)
+ \left({\bf q}_{r1}\,{\bf q}^{\dagger}_{d1}+{\bf q}_{d1}\,{\bf q}^{\dagger}_{r1}\right)\varrho^2_{\!_{\bf Y}}\,\right\}\varepsilon \label{sd40}\\
&=\text{\{a scalar\}}+\text{\{a scalar\}}\,\varepsilon\,, \notag \\
&=\widetilde{||{\bf X}||^2}\;\widetilde{||{\bf Y}||^2} \label{sd41}
\end{align}
thus proving the multivector-valued composition law (\ref{mulmul}) we set out to prove. Here, I have identified the 3-sphere radii of arbitrary values as $\varrho^{}_{r1}=\sqrt{{\bf q}_{r1}\,{\bf q}^{\dagger}_{r1}}\,$, {\it etc.}, and Pythagorean distances from the origin of $\mathrm{I\!R}^8$ as $\varrho^{}_{\!_{\bf X}}=\sqrt{\varrho^2_{r1}+\,\varrho^2_{d1}\,}$ and $\varrho^{}_{\!_{\bf Y}}=\sqrt{\varrho^2_{r2}+\,\varrho^2_{d2}\,}$. It is also easy to recognize that the quantities appearing in the two curly brackets are scalar quantities because, as we saw, ${\bf q}_{r1}\,{\bf q}^{\dagger}_{d1}+{\bf q}_{d1}\,{\bf q}^{\dagger}_{r1}=2\,g_1h_1 + 2\,{\bf u}_1\cdot{\bf v}_1$, {\it etc.}, are scalar quantities. Since (\ref{sd40}) holds for any\break multivectors ${\bf X}$ and ${\bf Y}$ in ${\cal K}^{\lambda}$, the composition law (\ref{mulmul}) holds for {\it any} multivectors ${\bf X}$ and ${\bf Y}$ in ${\cal K}^{\lambda}$, {\it without exception}.

\subsection{Multivector-valued Composition Law holds Even for Non-zero Zero Divisors}

\renewcommand{\thesubsubsection}{D.3.\arabic{subsubsection}}

To verify the consistency of (\ref{mulmul}), consider the following two multivectors in ${\cal K}^{\lambda}$, where $\alpha$ and $\gamma$ are real numbers:
\begin{equation}
{\bf W} = \alpha + \alpha\;\varepsilon\;\;\;\text{and}\;\;\;{\bf Z}=\gamma - \gamma\;\varepsilon\,. \label{wz} 
\end{equation}
Using the property $\varepsilon^2=1$ of the pseudoscalar, it is easy to verify that the geometric product of ${\bf W}$ and ${\bf Z}$ vanishes:
\begin{equation}
{\bf W}{\bf Z}=(\alpha + \alpha\;\varepsilon)\;(\gamma - \gamma\;\varepsilon) =\alpha\,\gamma\,(1+\,\varepsilon)\,(1-\,\varepsilon) = \alpha\,\gamma\,(1-\,\varepsilon^2)=0\,. \label{norm-wz}
\end{equation}
Since the geometric product ${\bf W}{\bf Z}$ is the only way to construct another multivector in ${\cal K}^{\lambda}$, its vanishing may lead one to suspect that perhaps the multivector-valued composition law (\ref{mulmul}) we have proved may not always be respected. But that is not the case, as we can easily verify. To begin with, note that neither ${\bf W}$ nor ${\bf Z}$ can be normalized to a scalar magnitude unless the coefficients $\alpha$ and $\gamma$ of $\varepsilon$ are vanishing. In fact, their norms turn out to be multivector-valued:    
\begin{equation}
\widetilde{||{\bf W}||}=\alpha\,\widetilde{||(1+\,\varepsilon)||}=\alpha\,\sqrt{(1+\,\varepsilon)(1+\,\varepsilon)^{\dagger}}=\alpha\,\sqrt{(1+\,\varepsilon)(1+\,\varepsilon)}=\alpha\,(1+\varepsilon)= \alpha + \alpha\;\varepsilon, \label{30X}
\end{equation}
where the property $\varepsilon^{\dagger}=\varepsilon$ is used, and, similarly,\begin{equation}
\widetilde{||{\bf Z}||}=\gamma\,\widetilde{||(1-\,\varepsilon)||}=\gamma\,\sqrt{(1-\,\varepsilon)(1-\,\varepsilon)^{\dagger}}=\gamma\,\sqrt{(1-\,\varepsilon)(1-\,\varepsilon)}=\gamma\,(1-\varepsilon)=\gamma - \gamma\;\varepsilon. \label{30Y}
\end{equation}
Consequently, using the property $\varepsilon^2=1$, the right-hand side of (\ref{mulmul}) gives:
\begin{equation}
\widetilde{||{\bf W}||}\;\widetilde{||{\bf Z}||}=\alpha\,\gamma\,(1+\,\varepsilon)(1-\,\varepsilon)=\alpha\,\gamma\,(1-\varepsilon^2)=\alpha\,\gamma\,(1-1)=0, \label{20XY}
\end{equation}
which matches with the left-hand side of (\ref{mulmul}):
\begin{equation}
\widetilde{||{\bf W}{\bf Z}||}=\alpha\,\gamma\,\widetilde{||(1+\,\varepsilon)(1-\,\varepsilon)||}=\alpha\,\gamma\,\widetilde{\left|\left|1-\,\varepsilon^2\right|\right|}=\alpha\,\gamma\,\widetilde{\left|\left|1-1\right|\right|}=\alpha\,\gamma\,\widetilde{||\,0\,||}=0. \label{-ap}
\end{equation}
Thus, while non-zero zero divisors specified in (\ref{wz}) are not normalizable to scalar lengths using geometric products, they do not lead to inconsistency in the multivector-valued composition law (\ref{mulmul}). Moreover, we will see below that non-zero zero divisors, such as (\ref{wz}), reduce to ordinary zero divisors with vanishing scalar length upon normalization.

\subsection{Necessary and Sufficient Criterion for Computing Scalar-valued Norms}

The above consistency check of the multivector-valued composition law (\ref{mulmul}) would not have gone through had we followed the standard practice in Geometric Algebra for normalizing multivectors \cite{Clifford}. When normalizing a multivector such as ${\bf X}$, the standard practice in Geometric Algebra is to simply ignore the non-scalar part of the geometric product ${\bf X}{\bf X}^{\dagger}$ and declare that its scalar component is the squared magnitude of the multivector ${\bf X}$, rather than recognizing that, in fact, it is a {\it projection} of the multivector-valued magnitude of that multivector onto the scalar axis in the graded basis of the corresponding vector space (cf. Fig.~\ref{fig2} below). Although this practice is harmless for most computations \cite{Clifford}, it leads to inconsistencies when applied to non-trivial multivector-valued equations, such as our composition law (\ref{mulmul}), as I have explained in greater detail in Section~2.7 of the Royal Society paper listed below in the bibliography as Ref.~\cite{RSOS-2}. What is more, ignoring the non-scalar parts of the multivector ${\bf X}{\bf X}^{\dagger}$ amounts to inadequately abandoning the geometric product for the convenience of normalizing a multivector quantity, which is otherwise regarded as {\it the} fundamental product in Geometric Algebra. Indeed, the very purpose of a geometric product such as (\ref{gp}) is to retain both the scalar and non-scalar parts of the multivector-valued members of a graded algebra to maintain their algebraic integrity. It is therefore hardly appropriate to suddenly switch to the scalar product simply because the geometric product is inconvenient for computing a scalar-valued norm of a multivector-valued quantity. Instead of ignoring the non-scalar part of the geometric product ${\bf X}{\bf X}^{\dagger}$ for the purpose of finding the scalar magnitude of a multivector ${\bf X}$, a consistent strategy is to set the {\it coefficient} of the non-scalar part of ${\bf X}{\bf X}^{\dagger}$ to zero without compromising geometric product itself, and view the resulting equation as a constraint that reduces the dimensions of algebra by one; {\it i.e.}, reduces it to the dimensions of a 7-sphere embedded in the corresponding vector space $\mathrm{I\!R}^8$. Any other criterion would be inconsistent, as it would necessitate applying different product rules on the two sides of the composition law (\ref{comnorm}).

\subsection{Proof of the Scalar-valued Composition Law $\,||{\bf X}{\bf Y}|| = ||{\bf X}||\;||{\bf Y}||$ for ${\cal K}^{\lambda}$}

\begin{figure}
\hrule
\scalebox{1}{
\begin{pspicture}(5,-4.7)(6.2,5.9)

\psline[linewidth=0.3mm,arrowinset=0.3,arrowsize=3pt 3,arrowlength=2]{->}(-1.27,-2.3)(5.5,-2.3)

\psline[linewidth=0.3mm,arrowinset=0.3,arrowsize=3pt 3,arrowlength=2]{->}(-0.77,-2.8)(-0.77,3.0)

\psline[linewidth=0.3mm,arrowinset=0.3,arrowsize=3pt 3,arrowlength=2]{->}(-0.77,-2.3)(4.0,-2.3)

\psline[linewidth=0.3mm,arrowinset=0.3,arrowsize=3pt 3,arrowlength=2]{->}(-0.77,-2.3)(4.0,1.5)

\psarc[linewidth=0.3mm,arrowinset=0.3,arrowsize=3pt 3,arrowlength=2]{<->}(0.0,-2.5){3.98}{335}{30}

\psline[linewidth=0.1mm,dotsize=2pt 3]{*-}(4.0,1.5)(4.0,1.55)

\psline[linewidth=0.1mm,dotsize=2pt 3]{*-}(-0.77,-2.31)(-0.77,-2.27)

\psline[linewidth=0.2mm,linestyle=dashed,arrowsize=4pt 4,arrowlength=2]{<-}(4.0,-0.77)(4.0,1.5) 

\psline[linewidth=0.2mm,linestyle=dashed,arrowsize=3pt 3,arrowlength=2]{-}(4.0,-2.3)(4.0,-0.77) 

\rput*[0.0]{269}(4.5,-0.4){${\mathbf q}_{\mathbf r}\,{\mathbf q}^{\dagger}_{\mathbf d}+{\mathbf q}_{\mathbf d}\,{\mathbf q}^{\dagger}_{\mathbf r}\,=\,0$}

\put(1.0,-2.83){\large ${||{\mathbb Q}_z||=\varrho}$}

\put(4.2,1.7){\large ${\sqrt{{\mathbb Q}_z{\mathbb Q}_z^{\dagger}}}$}

\put(8.65,-2.04){${{\bf q}_{r}\,{\bf q}^{\dagger}_{d}+{\bf q}_{d}\,{\bf q}^{\dagger}_{r}\,=\,0}$}

\put(5.65,-2.4){\large $1$}

\put(-0.89,3.2){\Large $\varepsilon$}

\put(-1.23,-2.85){\large ${0}$}

\psline[linewidth=0.3mm,arrowinset=0.3,arrowsize=3pt 3,arrowlength=2]{->}(8.5,-1.0)(11.0,-1.0)

\psline[linewidth=0.4mm,arrowinset=0.3,arrowsize=3pt 3,arrowlength=2]{->}(8.5,-1.0)(8.5,1.3)

\put(-1.55,0.5){dual}

\put(-1.52,0.1){axis}

\put(6.1,-2.1){real}

\put(6.1,-2.5){axis}

\put(9.05,-0.45){${\SI{90}{\degree}}$}

\put(11.1,-1.05){\large ${{\bf q}_r}$}

\put(8.3,1.5){\large ${{\bf q}_d}$}

\put(9.9,0.3){${{\bf q}_d\,=v_0 + I_3\,{\bf v}}$}

\put(9.9,0.7){${{\bf q}_r=u_0 + I_3\,{\bf u}}$}

\put(8.2,-1.6){Normalization condition:}

\put(7.8,-2.5){$u_0v_0 + u_xv_x+u_yv_y+u_zv_z=0$}

\psline[linewidth=0.2mm,arrowinset=0.3,arrowsize=3pt 3,arrowlength=2]{-}(8.5,-0.5)(9.0,-0.5)

\psline[linewidth=0.2mm,arrowinset=0.3,arrowsize=3pt 3,arrowlength=2]{-}(9.0,-1.0)(9.0,-0.5)

\put(4.0,-3.7){\Large ${S^7\hookrightarrow{\mathrm{I\!R}^8}}$}

\put(0.3,3.9){$\boxed{\,\varepsilon\cdot\left({\mathbb Q}_{z}{\mathbb Q}^{\dagger}_{z}\right):=\,{\bf q}_{r}\,{\bf q}^{\dagger}_{d}+{\bf q}_{d}\,{\bf q}^{\dagger}_{r}\,=\,2\,(u^{}_0\,v^{}_0 + u^{}_x\,v^{}_x+u^{}_y\,v^{}_y+u^{}_z\,v^{}_z)\;\quad\quad\quad\;0\,}$}

\psline[linewidth=0.2mm,arrowinset=0.3,arrowsize=3pt 3,arrowlength=2]{->}(9.6,4)(10.55,4)

\put(0.1,4.8){\Large {\bf Projection of $\widetilde{||{\mathbb Q}_z||}$ onto the Real dimensions}}

\end{pspicture}}
\hrule
\vspace{-5pt}
\caption{Illustration of a projection of the non-scalar quantity $\sqrt{{\mathbb Q}_z{\mathbb Q}_z^{\dagger}}$ onto the real axis of the graded sub-basis $\{1,\,\varepsilon\}$ in ${\cal K}^{\lambda}$, yielding the scalar magnitude $\varrho=||{\mathbb Q}_z||$ as the Pythagorean distance form origin in the vector space $\mathrm{I\!R}^8$ corresponding to ${\cal K}^{\lambda}$.}
\vspace{7pt}
\hrule
\label{fig2}
\end{figure}

We now have all the necessary ingredients to prove our central result; namely, the scalar-valued composition law (\ref{comnorm}). We begin with the generic multivector ${\mathbb Q}_z$ we expanded in (\ref{qinqq}) in the graded basis of ${\cal K}^{\lambda}$. With appropriate notation for its coefficients, it can be expressed as in (\ref{AAA-1}), ${\mathbb Q}_z=\, {\bf q}_{r} + {\bf q}_{d}\,\varepsilon$, within ${\cal K}^{\lambda}$, with ${\bf q}_{r}:=u_0+I_3{\bf u}$ and ${\bf q}_{d}:=v_0+I_3{\bf v}$. The quadratic form in the definition (\ref{defnorm}) with its image resembling a split-complex number is then
\begin{align}
{\mathbb Q}_{z}{\mathbb Q}^{\dagger}_{z}
&=\left({\bf q}_{r}\,{\bf q}^{\dagger}_{r}+{\bf q}_{d}\,{\bf q}^{\dagger}_{d}\right)+\left({\bf q}_{r}\,{\bf q}^{\dagger}_{d}+{\bf q}_{d}\,{\bf q}^{\dagger}_{r}\right)\varepsilon \label{ddd48} \\
&=\left(\varrho^2_r + \varrho^2_d\,\right)+\left( 2\,u^{}_0v^{}_0 + 2\,{\bf u}\cdot{\bf v}\right)\,\varepsilon \label{sd47} \\
&=\left(u_0^2+u_x^2+u_y^2+u_z^2 + v_0^2+v_x^2+v_y^2+v_z^2\right) + 2\left(u^{}_0\,v^{}_0 + u^{}_x\,v^{}_x+u^{}_y\,v^{}_y+u^{}_z\,v^{}_z\right)\,\varepsilon \\
&=\,\text{(a scalar)}+\text{(a scalar)}\;\varepsilon\,. \label{ddd51}
\end{align}
where $\varrho_{r}=\sqrt{{\bf q}_{r}\,{\bf q}^{\dagger}_{r}}$ and  $\varrho_{d}=\sqrt{{\bf q}_{d}\,{\bf q}^{\dagger}_{d}}$ are radii of 3-spheres in $\mathrm{I\!R}^4\hookrightarrow\mathrm{I\!R}^8$. The coefficient of the pseudoscalar $\varepsilon$ is thus
\begin{equation}
\varepsilon\cdot\left({\mathbb Q}_{z}{\mathbb Q}^{\dagger}_{z}\right):=\,{\bf q}_{r}\,{\bf q}^{\dagger}_{d}+{\bf q}_{d}\,{\bf q}^{\dagger}_{r}\,=\,2\,(u^{}_0\,v^{}_0 + u^{}_x\,v^{}_x+u^{}_y\,v^{}_y+u^{}_z\,v^{}_z)\,.
\end{equation}
In accordance with the criterion discussed above, we now set this scalar coefficient of $\varepsilon$ to zero to compute the norm: 
\begin{equation}
\varepsilon\cdot\left({\mathbb Q}_{z}{\mathbb Q}^{\dagger}_{z}\right):=\,{\bf q}_{r}\,{\bf q}^{\dagger}_{d}+{\bf q}_{d}\,{\bf q}^{\dagger}_{r}\,=\,2\,(u^{}_0\,v^{}_0 + u^{}_x\,v^{}_x+u^{}_y\,v^{}_y+u^{}_z\,v^{}_z)\,=\,0\,. \label{consttt}
\end{equation}
This normalization criterion is illustrated in Fig.~\ref{fig2} for clarity. It amounts to projecting the multivector-valued quantity ${\mathbb Q}_{z}{\mathbb Q}^{\dagger}_{z}$ onto the real dimensions. We can now define the scalar-valued norm of any multivector ${\mathbb Q}_{z}$ within ${\cal K}^{\lambda}$ as follows:
\begin{equation}
||{\mathbb Q}_z||:=\sqrt{\left({\mathbb Q}_{z}{\mathbb Q}^{\dagger}_{z}\right)\!\Big|_{\,\varepsilon\cdot\left({\mathbb Q}_{z}{\mathbb Q}^{\dagger}_{z}\right)\,=\,0}\,}\;. \label{defnorm}
\end{equation}
Since---as we can see from (\ref{sd47})---the scalar part of the geometric product ${\mathbb Q}_{z}{\mathbb Q}^{\dagger}_{z}$ is always of the form $\varrho^2_r + \varrho^2_d\,$, this definition of norm always yields a {\it positive definite} value for the norm of any multivector ${\mathbb Q}_{z}$ in ${\cal K}^{\lambda}$. Moreover, it gives the same value for the norm as that would be obtained by employing the traditional {\it ad hoc} criterion usually employed in practice in Geometric Algebra, but without leading to any inconsistencies (cf. Section 2.7 in \cite{RSOS}). It is also easy to see that definition (\ref{defnorm}) gives the value for the norm as a Pythagorean distance from the origin in eight dimensions:
\begin{equation}
||{\mathbb Q}_z||=\sqrt{\varrho_r^2+\varrho_d^2\,}=\sqrt{u_0^2+u_x^2+u_y^2+u_z^2+v_0^2+v_x^2+v_y^2+v_z^2\,}=\varrho\,. 
\end{equation}
Since this equation is equivalent to the following quadratic equation in eight variables and thus defines a 7-sphere,
\begin{equation}
u_0^2+u_x^2+u_y^2+u_z^2+v_0^2+v_x^2+v_y^2+v_z^2=\varrho^2,
\end{equation}
it is natural to identify $\varrho$ as a radius of this 7-sphere. However, Lasenby has claimed in Ref.~[14] of \cite{RSOS-2} that viewing $\varrho$ as a radius of a 7-sphere and choosing appropriate units to set $\varrho=1$ amounts to imposing another constraint, in addition to (\ref{consttt}), and thus reduces the dimensions of ${\cal K}^{\lambda}$ to $8-2=6$ rather than $8-1=7$. But this is incorrect. The constraint (\ref{consttt}) is not only necessary but also sufficient to determine $\varrho$ as a scalar-valued radius of a 7-sphere embedded in ${\cal K}^{\lambda}$. However, to prevent a Lasenby-type misinterpretation, let us not identify $\varrho$ at this stage as a radius of a 7-sphere, but view it simply as an arbitrary Pythagorean distance from the origin in $\mathrm{I\!R}^8$ that is variable within~it.

Next, we apply the definition (\ref{defnorm}) of normalization to the multivectors ${\bf X}$ and ${\bf Y}$ specified in (\ref{def-XY}) to determine their scalar-valued lengths using the set of equations (\ref{sd10a}) to (\ref{sd30}) we have previously worked out. Together with
\begin{equation}
{\bf q}_{r1}\,{\bf q}^{\dagger}_{d1}+{\bf q}_{d1}\,{\bf q}^{\dagger}_{r1}=0\,, \label{sd55}
\end{equation}
(\ref{defnorm}) then gives the scalar-valued norm of ${\bf X}$ as
\begin{align}
||{\bf X}||=\sqrt{\left({\mathbf{X}}{\mathbf{X}}^{\dagger}\right)\!\Big|_{\,{\bf q}_{r1}\,{\bf q}^{\dagger}_{d1}+{\bf q}_{d1}\,{\bf q}^{\dagger}_{r1}\,=\,0}\,}\,
&=\,\sqrt{g_1^2 + {\bf u}_1\cdot{\bf u}_1 + h_1^2 + {\bf v}_1\cdot{\bf v}_1\,} \\
&=\,\sqrt{g_1^2+u_{1x}^2+u_{1y}^2+u_{1z}^2+h_1^2+v_{1x}^2+v_{1y}^2+v_{1z}^2\,} \label{sd57} \\
&=\,\varrho^{}_{\!_{\bf X}}\,.
\end{align}
Similarly, together with
\begin{equation}
{\bf q}_{r2}\,{\bf q}^{\dagger}_{d2}+{\bf q}_{d2}\,{\bf q}^{\dagger}_{r2}=0\,, \label{sd59}
\end{equation}
(\ref{defnorm}) gives the scalar-valued norm of ${\bf Y}$ as
\begin{align}
||{\bf Y}||=\sqrt{\left({\mathbf{Y}}{\mathbf{Y}}^{\dagger}\right)\!\Big|_{\,{\bf q}_{r2}\,{\bf q}^{\dagger}_{d2}+{\bf q}_{d2}\,{\bf q}^{\dagger}_{r2}\,=\,0}\,}\,
&=\,\sqrt{g_2^2 + {\bf u}_2\cdot{\bf u}_2 + h_2^2 + {\bf v}_2\cdot{\bf v}_2\,} \\
&=\,\sqrt{g_2^2+u_{2x}^2+u_{2y}^2+u_{2z}^2+h_2^2+v_{2x}^2+v_{2y}^2+v_{2z}^2\,} \label{sd61} \\
&=\,\varrho^{}_{\!_{\bf Y}}\,.
\end{align}
Thus, on the right-hand side of the scalar-valued composition law (\ref{comnorm}), we have the following product of scalars:
\begin{equation}
||{\bf X}||\;||{\bf Y}|| = \varrho^{}_{\!_{\bf X}}\,\varrho^{}_{\!_{\bf Y}}\,. \label{sd63}
\end{equation}
Now, to evaluate the left-hand side of the scalar-valued composition law (\ref{comnorm}), we compute $({\bf X}{\bf Y})({\bf X}{\bf Y})^{\dagger}$ as follows:
\begin{align}
({\bf X}{\bf Y})({\bf X}{\bf Y})^{\dagger}
&=({\bf X}{\bf Y})({\bf Y}^{\dagger}{\bf X}^{\dagger}) \label{sd64} \\
&=\,{\bf X}\,({\bf Y}{\bf Y}^{\dagger})\,{\bf X}^{\dagger} \label{sd65}\\
&=\,{\bf X}\left[\left({\bf q}_{r2}\,{\bf q}^{\dagger}_{r2}+{\bf q}_{d2}\,{\bf q}^{\dagger}_{d2}\right)+\left({\bf q}_{r2}\,{\bf q}^{\dagger}_{d2}+{\bf q}_{d2}\,{\bf q}^{\dagger}_{r2}\right)\varepsilon\,\right]{\bf X}^{\dagger} \label{sd66} \\
&=\,{\bf X}\left[\left(g_2^2 + {\bf u}_2\cdot{\bf u}_2 + h_2^2 + {\bf v}_2\cdot{\bf v}_2\right)+\,0\,\right]{\bf X}^{\dagger} \label{sd67} \\
&=\,\left({\bf X}{\bf X}^{\dagger}\right)\left(g_2^2+u_{2x}^2+u_{2y}^2+u_{2z}^2+h_2^2+v_{2x}^2+v_{2y}^2+v_{2z}^2\right) \\
&=\,\left[\left({\bf q}_{r1}\,{\bf q}^{\dagger}_{r1}+{\bf q}_{d1}\,{\bf q}^{\dagger}_{d1}\right)+\left({\bf q}_{r1}\,{\bf q}^{\dagger}_{d1}+{\bf q}_{d1}\,{\bf q}^{\dagger}_{r1}\right)\varepsilon\,\right]\left(g_2^2+u_{2x}^2+u_{2y}^2+u_{2z}^2+h_2^2+v_{2x}^2+v_{2y}^2+v_{2z}^2\right) \label{sd69} \\
&=\,\left[\left(g_1^2 + {\bf u}_1\cdot{\bf u}_1 + h_1^2 + {\bf v}_1\cdot{\bf v}_1\right) +\,0\,\right] \left(g_2^2+u_{2x}^2+u_{2y}^2+u_{2z}^2+h_2^2+v_{2x}^2+v_{2y}^2+v_{2z}^2\right) \label{sd70} \\
&=\,\left(g_1^2+u_{1x}^2+u_{1y}^2+u_{1z}^2+h_1^2+v_{1x}^2+v_{1y}^2+v_{1z}^2\right)\left(g_2^2+u_{2x}^2+u_{2y}^2+u_{2z}^2+h_2^2+v_{2x}^2+v_{2y}^2+v_{2z}^2\right) \\
&=\,\varrho^2_{\!_{\bf X}}\,\varrho^2_{\!_{\bf Y}}\,. \label{sd72}
\end{align}
Here, step (\ref{sd65}) follows from the step (\ref{sd64}) by the associativity of the geometric product in ${\cal K}^{\lambda}$, in step (\ref{sd66}) I have used (\ref{sd10bc}) to substitute for ${\bf Y}{\bf Y}^{\dagger}$, in step (\ref{sd67}) I have used the orthogonality condition (\ref{sd55}), in step (\ref{sd69}) I have used (\ref{sd10b}) to substitute for ${\bf X}{\bf X}^{\dagger}$, and in step (\ref{sd70}) I have used the orthogonality condition (\ref{sd59}). Consequently, the right-hand side of (\ref{sd72}) is automatically a real number. In other words, the coefficient of $\varepsilon$ is already zero, giving
\begin{equation}
||{\bf X}{\bf Y}||:=\sqrt{\left[\,({\bf X}{\bf Y})({\bf X}{\bf Y})^{\dagger}\,\right]\!\Big|_{\,\varepsilon\cdot\left[\,({\bf X}{\bf Y})({\bf X}{\bf Y})^{\dagger}\right]\;=\;0}\,}\;\;=\;\varrho^{}_{\!_{\bf X}}\,\varrho^{}_{\!_{\bf Y}}\,. \label{sd73}
\end{equation}
Comparing (\ref{sd73}) with (\ref{sd63}), we conclude that
\begin{equation}
||{\bf X}{\bf Y}|| = \,\varrho^{}_{\!_{\bf X}}\,\varrho^{}_{\!_{\bf Y}} \!=\, ||{\bf X}||\,||{\bf Y}||\,, \label{sd74} 
\end{equation}
which is the scalar-valued composition law (\ref{comnorm}) we set out to prove. Eq.~(\ref{sd74}) proves that ${\cal K}^{\lambda}$ is a {\it normed} algebra. Thus, the multivector-valued composition law (\ref{mulmul}) reduces to the scalar-valued composition law (\ref{comnorm}) for ${\bf q}_{r}\perp{\bf q}_{d}\,$:
\makeatletter 
\renewcommand*{\rightarrowfill@}{%
  \arrowfill@\relbar\relbar\chemarrow} 
\makeatother
\begin{equation}
\widetilde{||{\bf X}{\bf Y}||} =\, \widetilde{||{\bf X}||}\;\widetilde{||{\bf Y}||}\;\;\xrightarrow[\;\;\;{\bf q}_{r}\,{\bf q}^{\dagger}_{d}\,+\;{\bf q}_{d}\,{\bf q}^{\dagger}_{r}\;\longrightarrow\;\,0\;\;\;]{{\bf q}_{r}\,\perp\;{\bf q}_{d}}\;\;||{\bf X}{\bf Y}||= \,||{\bf X}||\;||{\bf Y}||\,. \label{ed77}
\end{equation}
We can also see this projection of the composition law (\ref{mulmul}) onto real dimensions immediately from (\ref{sd40}) (cf. Fig.~\ref{fig2}).

Next, recall that ${\cal K}^{\lambda}$ remains closed under multiplication, giving ${\bf X}{\bf Y}=\,{\bf Z} \in {\cal K}^{\lambda}$, which can be parameterized as 
\begin{equation}
{\bf Z}=\,{\bf q}_{r3} + {\bf q}_{d3}\,\varepsilon\;\;\;\;\text{with}\;\;\;\;{\bf q}_{r3}=g_3+I_3{\bf u}_3\;\;\;\;\text{and}\;\;\;\;{\bf q}_{d3}=h_3+I_3{\bf v}_3\,,
\end{equation}
where $g_3$ and $h_3$ as scalars and ${{\bf u}_3=u_{3x}{\bf e}_x+u_{3y}{\bf e}_y+u_{3z}{\bf e}_z}$ and ${{\bf v}_3=v_{3x}{\bf e}_x+v_{3y}{\bf e}_y+v_{3z}{\bf e}_z}$ are vectors in $\mathrm{I\!R}^3$. Then, just as for ${\bf X}{\bf X}^{\dagger}$ and ${\bf Y}{\bf Y}^{\dagger}$, the geometric product ${\bf Z}{\bf Z}^{\dagger}$ between the multivector ${\bf Z}$ and its reverse ${\bf Z}^{\dagger}$ works out~to~be
\begin{align}
{\bf Z}{\bf Z}^{\dagger}\,
&=\left({\bf q}_{r3} + {\bf q}_{d3}\,\varepsilon\right)\left({\bf q}_{r3} + {\bf q}_{d3}\,\varepsilon\right)^{\dagger} \label{10mnmnmn}\\
&=\left({\bf q}_{r3}\,{\bf q}^{\dagger}_{r3}+{\bf q}_{d3}\,{\bf q}^{\dagger}_{d3}\right)+\left({\bf q}_{r3}\,{\bf q}^{\dagger}_{d3}+{\bf q}_{d3}\,{\bf q}^{\dagger}_{r3}\right)\varepsilon \label{10zxzxzxz} \\
&=\left(g_3^2 + {\bf u}_3\cdot{\bf u}_3 + h_3^2 + {\bf v}_3\cdot{\bf v}_3\right) + \left({\bf q}_{r3}\,{\bf q}^{\dagger}_{d3}+{\bf q}_{d3}\,{\bf q}^{\dagger}_{r3}\right)\varepsilon \,. \label{opopop}
\end{align}
Consequently, together with
\begin{equation}
{\bf q}_{r3}\,{\bf q}^{\dagger}_{d3}+{\bf q}_{d3}\,{\bf q}^{\dagger}_{r3}=0\,,
\end{equation}
(\ref{defnorm}) gives the scalar-valued norm of ${\bf Z}$ as
\begin{align}
||{\bf Z}||=\sqrt{\left({\mathbf{Z}}{\mathbf{Z}}^{\dagger}\right)\!\Big|_{\,{\bf q}_{r3}\,{\bf q}^{\dagger}_{d3}+{\bf q}_{d3}\,{\bf q}^{\dagger}_{r3}\,=\,0}\,}\,
&=\,\sqrt{g_3^2 + {\bf u}_3\cdot{\bf u}_3 + h_3^2 + {\bf v}_3\cdot{\bf v}_3\,} \\
&=\,\sqrt{g_3^2+u_{3x}^2+u_{3y}^2+u_{3z}^2+h_3^2+v_{3x}^2+v_{3y}^2+v_{3z}^2\,} \label{sd81} \\
&=\,\varrho^{}_{\!_{\bf Z}}\,. \label{sd82}
\end{align}
Using the results from (\ref{sd74}) and (\ref{sd82}), together with ${\bf Z}={\bf X}{\bf Y}$, or equivalently $||{\bf Z}||=||{\bf X}{\bf Y}||$, we can conclude that
\begin{equation}
\varrho^{}_{\!_{\bf Z}}=\,\varrho^{}_{\!_{\bf X}} \,\varrho^{}_{\!_{\bf Y}}\,, \label{sd83}
\end{equation}
which is equivalent to
\begin{equation}
\varrho^2_{\!_{\bf Z}}=\,\varrho^2_{\!_{\bf X}} \,\varrho^2_{\!_{\bf Y}}\,. \label{sd84}
\end{equation}
Then, using the results from (\ref{sd57}), (\ref{sd61}), (\ref{sd81}), (\ref{sd83}), and (\ref{sd84}), we can further conclude the following equality:
\begin{align}
\big(g_3^2+u_{3x}^2+&\,u_{3y}^2+u_{3z}^2+h_3^2+v_{3x}^2+v_{3y}^2+v_{3z}^2\big) \notag \\
&=\left(g_1^2+u_{1x}^2+u_{1y}^2+u_{1z}^2+h_1^2+v_{1x}^2+v_{1y}^2+v_{1z}^2\right)
\left(g_2^2+u_{2x}^2+u_{2y}^2+u_{2z}^2+h_2^2+v_{2x}^2+v_{2y}^2+v_{2z}^2\right). \label{sd85}
\end{align}

Apart from different notation, this equality is identical to the one shown in (\ref{sd5}). It exhibits a sum of eight squares, factorized into a product of two other sums of eight squares. This is the equality discovered, independently, by Degen, Graves, and Cayley, and later proved by Hurwitz \cite{Hurwitz-1898} to be one of the only four possible equalities, in dimensions 1, 2, 4, and 8, corresponding to the parallelizable spheres $S^0$, $S^1$, $S^3$, and $S^7$. Recall now that, so far in this appendix, I have refrained from identifying the Pythagorean distances $\varrho^{}_{\!_{\bf Z}}$, $\varrho^{}_{\!_{\bf X}}$, and $\varrho^{}_{\!_{\bf Y}}$ from the origin in $\mathrm{I\!R}^8$ as radii of a 7-sphere. But thanks to Hurwitz's theorem, this identification is now inescapable. Only the radii of a 7-sphere can satisfy the factorization displayed in (\ref{sd85}). This refutes Lasenby's claim in Ref.~[14] of \cite{RSOS-2}. No such factorization is possible in any dimension other than eight, such as, for example, in seven dimensions, embedding a six-dimensional sphere, $S^6$, as claimed by Lasenby. As a result, we can now identify the spherical surface we have arrived at, as a 7-sphere of radius ${\varrho=\sqrt{\varrho_r^2+\varrho_d^2\,}}$ embedded in the eight-dimensional vector space $\mathrm{I\!R}^8$ corresponding to the normed algebra ${\cal K}^{\lambda}$:
\begin{equation}
{\cal K}^{\lambda}\hookleftarrow S^7=\,\left\{\,{\mathbb Q}_z=\,{\bf q}_r + {\bf q}_d\,\varepsilon\;\,\bigg|\;\,||{\mathbb Q}_z||=\sqrt{\left({\mathbb Q}_{z}{\mathbb Q}^{\dagger}_{z}\right)\!\Big|_{\,{\bf q}_{r}\,{\bf q}^{\dagger}_{d}+{\bf q}_{d}\,{\bf q}^{\dagger}_{r}\,=\,0}\,}\,=\sqrt{\varrho_r^2+\varrho_d^2\,}=\varrho\,\right\}. \label{sevsp}
\end{equation}
In words, this 7-sphere is a set of all orthogonal pairs $\{{\bf q}_r,\,{\bf q}_d\}$  of quaternions ${\bf q}_r$ and ${\bf q}_d$ of radii ${\varrho^{}_r}$ and ${\varrho^{}_d}$, respectively, with their graded sum  ${\mathbb Q}_z=\,{\bf q}_r + {\bf q}_d\,\varepsilon$ representing the multivectors in ${\cal K}^{\lambda}$, constrained by the orthogonality condition
\begin{equation}
{\bf q}_{r}\,{\bf q}^{\dagger}_{d}+{\bf q}_{d}\,{\bf q}^{\dagger}_{r}\,=\,2\,(u^{}_0\,v^{}_0 + u^{}_x\,v^{}_x+u^{}_y\,v^{}_y+u^{}_z\,v^{}_z)\,=\,0\,, \label{sd87}
\end{equation}
which reduces the eight dimensions of ${\cal K}^{\lambda}$ to the seven dimensions of $S^7$, with topology different from the octonionic~$S^7$. 

\subsection{Normed Algebra ${\cal K}^{\lambda}$ is a Division Algebra}

It is now straightforward to prove that there are no non-zero zero divisors in the normed algebra ${\cal K}^{\lambda}$ because it respects the scalar-valued composition law (\ref{comnorm}) for all multivectors within it that can be normalized to yield scalar magnitudes, as illustrated in Fig.~\ref{fig2}. If two multivectors ${\bf X}$ and ${\bf Y}$ within ${\cal K}^{\lambda}$ are such that their geometric product vanishes to give ${\bf X}{\bf Y}=0$, then it is evident from the definition (\ref{defnorm}) of scalar-valued norms that it gives $||{\bf X}{\bf Y}||=0$. And then from the scalar-valued composition law (\ref{comnorm}) we proved in (\ref{sd74}), it follows that the product of the scalar-valued norms of ${\bf X}$ and ${\bf Y}$ given by (\ref{defnorm}) would also vanish, $||{\bf X}||\;||{\bf Y}||=0$, which implies that $||{\bf X}||=0$ or $||{\bf Y}||=0$, because, as we noted earlier, the norms defined by (\ref{defnorm}) are always positive definite. Thus, the conclusion stated in the paper is inescapable --- once correctly normalized, no ``non-zero zero divisors'' remain in the normed algebra ${\cal K}^{\lambda}$. 

We can also arrive at the same conclusion somewhat differently, and more systematically, as follows: Let ${\bf X}$ and ${\bf Y}$ be any two multivectors in ${\cal K}^{\lambda}$, subject to the normalization condition (\ref{sd87}), giving scalar-valued norms obtained by applying definition (\ref{defnorm}). Then their norms $||{\bf X}||$ and $||{\bf Y}||$ are scalar numbers. Therefore, the product of their norms $||{\bf X}||\;||{\bf Y}||$ can vanish if and only if either $||{\bf X}||=0$ or $||{\bf Y}||=0$. The positive definiteness of the norms defined by (\ref{defnorm}) then implies that this is possible if and only if either ${\bf X}=0$ or ${\bf Y}=0$. In other words, if both ${\bf X}\not=0$ and ${\bf Y}\not=0$, then $||{\bf X}||\;||{\bf Y}||\not=0$. But that implies ${||{\bf X}{\bf Y}||\not=0}$, because, according to (\ref{sd74}), $||{\bf X}||\;||{\bf Y}||=||{\bf X}{\bf Y}||$. The positive definiteness of the norms defined by (\ref{defnorm}) then implies ${{\bf X}{\bf Y}\not=0}$. Thus, for arbitrary  ${\bf X}\not=0$ and ${\bf Y}\not=0$ we have proved\break that ${{\bf X}{\bf Y}\not=0}$. Therefore, the normed algebra ${\cal K}^{\lambda}$ subject to the normalization condition (\ref{sd87}) is a division algebra.

Finally, it is worth noting that definition (\ref{defnorm}) of norm is valid also for the real and complex numbers, quaternions, and octonions, because for them the coefficient of $\varepsilon$ is naturally zero, where the complex numbers and quaternions are isomorphic, respectively, to the even sub-algebras of the Clifford algebras ${\mathrm{Cl}_{2,0}}$ and ${\mathrm{Cl}_{3,0}}$, whereas the algebra ${\cal K}^{\lambda}$ is isomorphic to the even subalgebra of the sixteen-dimensional algebra ${\mathrm{Cl}_{4,0}}$. Thus, the even sub-algebras of the Clifford algebras ${\mathrm{Cl}_{1,0}}$, ${\mathrm{Cl}_{2,0}}$, ${\mathrm{Cl}_{3,0}}$, and ${\mathrm{Cl}_{4,0}}$ form {\it associative} norm division algebras with respect to the norm defined in (\ref{defnorm}), isomorphic, respectively, to ${\mathbb R}$, ${\mathbb C}$, ${\mathbb H}$, and ${\cal K}^{\lambda}$, in the only possible dimensions 1, 2, 4, and 8 as proved in \cite{Hurwitz-1898}. 

\subsection{Conclusion}

In this appendix, I have proved anew that the eight-dimensional even subalgebra ${\cal K}^{\lambda}$ of the Clifford algebra ${\mathrm{Cl}_{4,0}}$ is an associative normed division algebra \cite{RSOS}. Since the algebra ${\cal K}^{\lambda}$ remains closed under multiplication, the components of any eight-dimensional multivector ${\bf Z}={\bf X}{\bf Y}$ in it, constructed from a geometric product of two other multivectors ${\bf X}$ and ${\bf Y}$ in it, satisfies the following equality, once the correct criterion to obtain scalar-valued norms is implemented:
\begin{equation}
\Bigg(\sqrt{z_1^2 + z_2^2 + z_3^2 + \dots + z_8^2\,}\;\Bigg)
=\Bigg(\sqrt{x_1^2 + x_2^2 + x_3^2 + \dots + x_8^2\,}\;\Bigg)\,\Bigg(\sqrt{y_1^2 + y_2^2 + y_3^2 + \dots + y_8^2\,}\;\Bigg).
\end{equation}
In fact, any such triple of the eight-dimensional multivectors in ${\cal K}^{\lambda}$ respects the composition law $\widetilde{||{\bf X}{\bf Y}||} = \widetilde{||{\bf X}||}\;\widetilde{||{\bf Y}||}$ for multivector-valued norms, {\it without exception}, including for the apparent ``non-zero zero divisors'', and reduces to the familiar composition law $||{\bf X}{\bf Y}|| = ||{\bf X}||\;||{\bf Y}||$ for scalar-valued norms once the said criterion to obtain scalar-valued norms is correctly employed. This is because the apparent ``non-zero zero divisors'' reduce to the trivial ``zero divisors'' once the correct criterion for normalizing multivector-valued elements in ${\cal K}^{\lambda}$ is recognized and properly implemented.

In summary:

\begin{enumerate}[label=(\arabic*),leftmargin=1.2cm]
    \item The eight-dimensional even subalgebra ${\cal K}^{\lambda}$ of the Clifford algebra ${\mathrm{Cl}_{4,0}}$ is a \underbar{normed division algebra.}
    \item The multivector-valued composition law (\ref{mulmul}) holds within ${\cal K}^{\lambda}$ without exception, even for ${\bf W}{\bf Z}=0$.
    \item \underbar{Non-zero zero divisors} consistently normalizable to yield scalar-valued magnitudes \underbar{do not exist} in ${\cal K}^{\lambda}$.
    \item The scalar-valued composition law (\ref{comnorm}), $||{\bf X}{\bf Y}||=||{\bf X}||\,||{\bf Y}||$, also holds within ${\cal K}^{\lambda}$ without exception.
    \item The associative algebra ${\cal K}^{\lambda}$ ensconces a 7-sphere of topology different from its octonionic counterpart.
    \item Nevertheless, $S^7$ defined in (\ref{sevsp}) is \underbar{parallelizable} using multivectors in ${\cal K}^{\lambda}$, as shown in Appendix~\ref{E}.
\end{enumerate}

\section{The 7-sphere Defined in (\ref{sevsp}) is Parallelizable using the Associative Algebra ${\cal K}^{\lambda}$} \label{E}

The goal of this appendix is to demonstrate explicitly that, just as the octonionic 7-sphere is parallelizable using the non-associative algebra of octonions, the 7-sphere defined in (\ref{sevsp}) is parallelizable using the associative algebra of multivectors in ${\cal K}^{\lambda}$. To appreciate this, let ${T_q\,S^7}$ denote the tangent space to ${S^7}$ at the tip of a multivector ${{\mathbb Q}_z \in {\cal K}^{\lambda}}$,
\begin{equation}
{\mathbb Q}_z={\bf q}_r + {\bf q}_d\;\varepsilon\,=\,q_0+q_1\,{\bf e}_x{\bf e}_y+q_2\,{\bf e}_z{\bf e}_x+q_3\,{\bf e}_y{\bf e}_z+q_4\,{\bf e}_x{\bf e}_{\infty}+q_5\,{\bf e}_y{\bf e}_{\infty}+q_6\,{\bf e}_z{\bf e}_{\infty}+q_7\,I_3{\bf e}_{\infty}\,. \label{genvec-2}
\end{equation}
As before, I have set $\lambda=+1$, ${\varepsilon =-\,I_3{\bf e}_{\infty}\,}$,
\begin{equation}
{\bf q}_r =\,q_0 + q_1\,{\bf e}_x{\bf e}_y+q_2\,{\bf e}_z{\bf e}_x+q_3\,{\bf e}_y{\bf e}_z\,,\;\;\;\;\,
\text{and}\;\;\;\;\,{\bf q}_d =\,-q_7 + q_6\,{\bf e}_x{\bf e}_y+q_5\,{\bf e}_z{\bf e}_x+q_4\,{\bf e}_y{\bf e}_z\,.
\end{equation}
In this notation, the normalization of ${\mathbb Q}_z\,$, defined in the previous appendix by 
\begin{equation}
||{\mathbb Q}_z||:=\sqrt{\left({\mathbb Q}_{z}{\mathbb Q}^{\dagger}_{z}\right)\!\Big|_{\,\varepsilon\cdot\left({\mathbb Q}_{z}{\mathbb Q}^{\dagger}_{z}\right)\,=\,0}\,}\;\;, \label{E3}
\end{equation}
gives the following scalar-valued radius of the 7-sphere defined in (\ref{sevsp}): 
\begin{equation}
\varrho = ||{\mathbb Q}_z|| = \sqrt{q_0^2+q_1^2+q_2^2+q_3^2+q_4^2+q_5^2+q_6^2+q_7^2\,}\,,
\end{equation}
which is deduced by setting to zero the coefficient of the pseudoscalar ${\varepsilon =-\,I_3{\bf e}_{\infty}\,}$ in the product ${\mathbb Q}_{z}{\mathbb Q}^{\dagger}_{z}\,$, as in (\ref{consttt}):
\begin{equation}
\varepsilon\cdot\left({\mathbb Q}_{z}{\mathbb Q}^{\dagger}_{z}\right)
=\,{\bf q}_{r}\,{\bf q}^{\dagger}_{d}+{\bf q}_{d}\,{\bf q}^{\dagger}_{r}
\,= 2\left(-\,q_0q_7 + q_1q_6 + q_2q_5 + q_3q_4\right) = 0\,. \label{E4}
\end{equation}

We are now in a position to define the space $T_q\,S^7$ tangent to the 7-sphere at the tip of the multivector ${{\mathbb Q}_z \in {\cal K}^{\lambda}}$~as
\begin{equation}
T_q\,S^7:=\left\{({\mathbb Q}_z,\,{\mathbb t}_{{\mathbb Q}_z})\;\Big|\;
{\mathbb Q}_z\,,\;{\mathbb t}_{{\mathbb Q}_z}\!\in {\cal K}^{\lambda},\,\;||{{\mathbb Q}_z}||=\,\varrho\,,\,\;\langle\,{\mathbb Q}_z\,,\,{\mathbb t}_{{\mathbb Q}_z}\,\rangle =\,0\right\},
\label{tanspace}
\end{equation}
where ${\langle\,{\mathbb Q}_z\,,\,{\mathbb t}_{{\mathbb Q}_z}\,\rangle}$ represents the usual inner product $\langle\,\cdot\;,\,\cdot\,\rangle$ in ${\mathrm{I\!R}^8}$ between ${\mathbb Q}_z$ and any multivector ${\mathbb t}_{{\mathbb Q}_z\!} \in
{\cal K}^{\lambda}$ tangent to ${{\mathbb Q}_z}$ at the tip of ${{\mathbb Q}_z}$. Then, by denoting the tip of ${{\mathbb Q}_z}$ as ${q={\mathbb Q}_z \cap S^7}$, the tangent bundle of ${S^7}$ can be expressed as
\begin{equation}
{\rm T}\,S^7\,=\!\bigcup_{\,q\,\in\, S^7}\{\,q\,\}\times T_q\,S^7.\label{tanbundle}
\end{equation}
As we will soon see, this tangent bundle turns out to be a trivial product bundle: ${\mathrm T}S^7\,\equiv\,S^7\!\times{\rm I\!R}^7.\label{cfeq3}$ The triviality of the bundle ${{\rm T}S^7}$ implies that the 7-sphere we have defined in (\ref{sevsp}) is parallelizable --- {\it i.e.}, it admits seven multivectors that are linearly independent everywhere and vanishing nowhere. We can use these multivectors to specify a basis for the tangent space at each point $q$, thereby defining a global anholonomic frame on $S^7$ that fixes each point uniquely.

To verify this, let a tangent space at the tip of the multivector 
${\,{\mathbb Q}_0=(1,\,0,\,0,\,0,\,0,\,0,\,0,\,0)}$ be spanned by the basis
\begin{align}
&\left\{\,\beta_1({\mathbb Q}_0),\;\beta_2({\mathbb Q}_0),\;\beta_3({\mathbb Q}_0),\;\beta_4({\mathbb Q}_0),\;\beta_5({\mathbb Q}_0),\;\beta_6({\mathbb Q}_0),\;\beta_7({\mathbb Q}_0)\,\right\}\equiv
\left\{\,{\bf e}_x{\bf e}_y\,,\,{\bf e}_z{\bf e}_x\,,\,{\bf e}_y{\bf e}_z\,,\,{\bf e}_x{\bf e}_{\infty}\,,\,{\bf e}_y{\bf e}_{\infty}\,,\,{\bf e}_z{\bf e}_{\infty}\,,\,I_3{\bf e}_{\infty}\,\right\}.
\end{align}
This set of graded bases in the subalgebra ${\cal K}^{\lambda}$ can also be expressed as the following seven orthonormal multivectors:
\begin{subequations}
\begin{align}
\beta_1({\mathbb Q}_0)&=(\,0,\,\;1,\,\;0,\,\;0,\,\;0,\,\;0,\,\;0,\,\;0\,), \\
\beta_2({\mathbb Q}_0)&=(\,0,\,\;0,\,\;1,\,\;0,\,\;0,\,\;0,\,\;0,\,\;0\,), \\
\beta_3({\mathbb Q}_0)&=(\,0,\,\;0,\,\;0,\,\;1,\,\;0,\,\;0,\,\;0,\,\;0\,), \\
\beta_4({\mathbb Q}_0)&=(\,0,\,\;0,\,\;0,\,\;0,\,\;1,\,\;0,\,\;0,\,\;0\,), \\
\beta_5({\mathbb Q}_0)&=(\,0,\,\;0,\,\;0,\,\;0,\,\;0,\,\;1,\,\;0,\,\;0\,), \\
\beta_6({\mathbb Q}_0)&=(\,0,\,\;0,\,\;0,\,\;0,\,\;0,\,\;0,\,\;1,\,\;0\,), \\
\beta_7({\mathbb Q}_0)&=(\,0,\,\;0,\,\;0,\,\;0,\,\;0,\,\;0,\,\;0,\,\;1\,).
\end{align}
\end{subequations}
These multivectors satisfy the usual inner product in space $\mathrm{I\!R}^8$  defined by the map $\langle\,\cdot\;,\,\cdot\,\rangle\!: T_q\,S^7\times T^*_qS^7 \rightarrow {\rm I\!R}$~so~that
\begin{align}
\langle\beta_{\mu}({\mathbb Q}_0),\,\beta_{\nu}({\mathbb Q}_0)\rangle\,=\,\delta_{\mu\nu}\,,\label{orige-flat}
\end{align}
where indices ${\mu}$, ${\nu}$ run from 1 to 7 and ${T^{*}_qS^7}$ is a cotangent space at ${q}$. Evidently, these multivectors forming the basis of the tangent space are orthogonal also to the multivector ${\,{\mathbb Q}_0=(1,\,0,\,0,\,0,\,0,\,0,\,0,\,0)}$. Then, the tangent space basis 
\begin{equation}
\left\{\,\beta_1({\mathbb Q}_z),\;\beta_2({\mathbb Q}_z),\;\beta_3({\mathbb Q}_z),\;\beta_4({\mathbb Q}_z),\;\beta_5({\mathbb Q}_z),\;\beta_6({\mathbb Q}_z),\;\beta_7({\mathbb Q}_z)\right\} \label{E17}
\end{equation}
at the tip of any arbitrary multivector with same normalization $||{\mathbb Q}_z||=||{\mathbb Q}_0||=||(1,\,0,\,0,\,0,\,0,\,0,\,0,\,0)||=1$, such as
\begin{align}
{\mathbb Q}_z &=\,q_0+q_1\,{\bf e}_x{\bf e}_y+q_2\,{\bf e}_z{\bf e}_x+q_3\,{\bf e}_y{\bf e}_z+q_4\,{\bf e}_x{\bf e}_{\infty}+q_5\,{\bf e}_y{\bf e}_{\infty}+q_6\,{\bf e}_z{\bf e}_{\infty}+q_7\,I_3{\bf e}_{\infty} \notag \\
&=(q_0,\;q_1,\;q_2,\;q_3,\;q_4,\;q_5,\;q_6,\;q_7) \label{E18}
\end{align}
that sweeps $S^7\!\hookrightarrow\!{\cal K}^{\lambda}$, can be found by taking its geometric product with the previous tangent basis elements~as~follows:
\begin{align}
\beta_1({\mathbb Q}_z)=\beta_1({\mathbb Q}_0)\,{\mathbb Q}_z &= ({{\bf e}_x}\,{{\bf e}_y})\,{\mathbb Q}_z \notag \\
&=\,-\,q_1+q_0\,{\bf e}_x{\bf e}_y-q_3\,{\bf e}_z{\bf e}_x+q_2\,{\bf e}_y{\bf e}_z+q_5\,{\bf e}_x{\bf e}_{\infty}-q_4\,{\bf e}_y{\bf e}_{\infty}-q_7\,{\bf e}_z{\bf e}_{\infty}+q_6\,I_3{\bf e}_{\infty} \notag \\
&=(-\,q_1,\;q_0,\,-q_3,\;q_2,\;q_5,\,-q_4,\,-q_7,\;q_6), \label{E19}\\
\beta_2({\mathbb Q}_z)=\beta_2({\mathbb Q}_0)\,{\mathbb Q}_z &= ({{\bf e}_z}\,{{\bf e}_x})\,{\mathbb Q}_z \notag \\
&=\,-\,q_2+q_3\,{\bf e}_x{\bf e}_y+q_0\,{\bf e}_z{\bf e}_x-q_1\,{\bf e}_y{\bf e}_z-q_6\,{\bf e}_x{\bf e}_{\infty}-q_7\,{\bf e}_y{\bf e}_{\infty}+q_4\,{\bf e}_z{\bf e}_{\infty}+q_5\,I_3{\bf e}_{\infty} \notag \\
&=(-\,q_2,\;q_3,\;q_0,\,-q_1,\,-q_6,\,-q_7,\;q_4,\;\,q_5), \\
\beta_3({\mathbb Q}_z)=\beta_3({\mathbb Q}_0)\,{\mathbb Q}_z &= ({{\bf e}_y}\,{{\bf e}_z})\,{\mathbb Q}_z \notag \\
&=\,-\,q_3-q_2\,{\bf e}_x{\bf e}_y+q_1\,{\bf e}_z{\bf e}_x+q_0\,{\bf e}_y{\bf e}_z-q_7\,{\bf e}_x{\bf e}_{\infty}+q_6\,{\bf e}_y{\bf e}_{\infty}-q_5\,{\bf e}_z{\bf e}_{\infty}+q_4\,I_3{\bf e}_{\infty} \notag \\
&=(-q_3,\,-q_2,\;q_1,\;q_0,\,-q_7,\;q_6,\,-q_5,\;q_4), \\
\beta_4({\mathbb Q}_z)=\beta_4({\mathbb Q}_0)\,{\mathbb Q}_z &= ({\bf e}_x\,{\bf e}_{\infty})\,{\mathbb Q}_z \notag \\
&=\,-\,q_4-q_5\,{\bf e}_x{\bf e}_y+q_6\,{\bf e}_z{\bf e}_x-q_7\,{\bf e}_y{\bf e}_z+q_0\,{\bf e}_x{\bf e}_{\infty}+q_1\,{\bf e}_y{\bf e}_{\infty}-q_2\,{\bf e}_z{\bf e}_{\infty}+q_3\,I_3{\bf e}_{\infty} \notag \\
&=(-\,q_4,\,-q_5,\;q_6,\,-q_7,\;q_0,\;q_1,\,-q_2,\;q_3), \\
\beta_5({\mathbb Q}_z)=\beta_5({\mathbb Q}_0)\,{\mathbb Q}_z &= ({{\bf e}_y}\,{{\bf e}_{\infty}})\,{\mathbb Q}_z \notag \\
&=\,-\,q_5+q_4\,{\bf e}_x{\bf e}_y-q_7\,{\bf e}_z{\bf e}_x-q_6\,{\bf e}_y{\bf e}_z-q_1\,{\bf e}_x{\bf e}_{\infty}+q_0\,{\bf e}_y{\bf e}_{\infty}+q_3\,{\bf e}_z{\bf e}_{\infty}+q_2\,I_3{\bf e}_{\infty} \notag \\
&=(-\,q_5,\;q_4,\,-q_7,\,-q_6,\,-q_1,\;q_0,\;q_3,\;q_2), \\
\beta_6({\mathbb Q}_z)=\beta_6({\mathbb Q}_0)\,{\mathbb Q}_z &= ({{\bf e}_z}\,{{\bf e}_{\infty}})\,{\mathbb Q}_z \notag \\
&=\,-\,q_6-q_7\,{\bf e}_x{\bf e}_y-q_4\,{\bf e}_z{\bf e}_x+q_5\,{\bf e}_y{\bf e}_z+q_2\,{\bf e}_x{\bf e}_{\infty}-q_3\,{\bf e}_y{\bf e}_{\infty}+q_0\,{\bf e}_z{\bf e}_{\infty}+q_1\,I_3{\bf e}_{\infty} \notag \\
&=(-\,q_6,\,-q_7,\,-q_4,\;\,q_5,\;q_2,\,-q_3,\;q_0,\;q_1), \\
\beta_7({\mathbb Q}_z)=\beta_7({\mathbb Q}_0)\,{\mathbb Q}_z &= (I_3\,{\bf e}_{\infty})\,{\mathbb Q}_z \notag \\
&=\,q_7-q_6\,{\bf e}_x{\bf e}_y-q_5\,{\bf e}_z{\bf e}_x-q_4\,{\bf e}_y{\bf e}_z-q_3\,{\bf e}_x{\bf e}_{\infty}-q_2\,{\bf e}_y{\bf e}_{\infty}-q_1\,{\bf e}_z{\bf e}_{\infty}+q_0\,I_3{\bf e}_{\infty} \notag \\
&=(\,q_7,\,-q_6,\,-q_5,\,-q_4,\,-q_3,\,-q_2,\,-q_1,\,q_0). \label{E25}
\end{align}
It is now straightforward to verify that the basis elements ${\left\{\,\beta_1({\mathbb Q}_z),\;\beta_2({\mathbb Q}_z),\;\beta_3({\mathbb Q}_z),\;\beta_4({\mathbb Q}_z),\;\beta_5({\mathbb Q}_z),\;\beta_6({\mathbb Q}_z),\;\beta_7({\mathbb Q}_z)\right\}}$ of the tangent space at point $q\in S^7$ are not only mutually orthogonal but also orthogonal to the arbitrary multivector ${\mathbb Q}_z$ we have chosen for this demonstration, with respect to the standard inner product in ${{\rm I\!R}^8}$, and thus define a {\it bona fide}\break tangent space ${T_q\,S^7=\,{\rm I\!R}^7}$ at the tip of that ${\mathbb Q}_z\in {\cal K}^{\lambda}$. For this purpose, it is necessary to verify only 28 unordered inner products among the 56 pairs of the eight multivectors displayed in (\ref{E18}) through (\ref{E25}), including the multivector~${\mathbb Q}_z$:
\begin{equation}
\binom{8}{2} = \frac{8!}{2!\left(8-2\right)!} = \frac{8\left(8-1\right)}{2} = 28.
\end{equation}
Among these 28 pairs, the inner products $\langle\beta_{\mu}({\mathbb Q}_z),\,\beta_{\nu}({\mathbb Q}_z)\rangle$ of 24 pairs with respect to the usual metric in $\mathrm{I\!R}^8$ vanish, which can be verified using components from (\ref{E18}) to (\ref{E25}). The inner products of the remaining 4 pairs also vanish because their computation leads precisely to the vanishing coefficient of the pseudoscalar, $\varepsilon\cdot\left({\mathbb Q}_{z}{\mathbb Q}^{\dagger}_{z}\right)=0$,~shown~in~(\ref{E4}):
\begin{align}
-\,\langle\,{\mathbb Q}_z\,,\,\beta_7({\mathbb Q}_z)\rangle
=\langle\beta_1({\mathbb Q}_z),\,\beta_6({\mathbb Q}_z)\rangle
=\langle\beta_2({\mathbb Q}_z),\,\beta_5({\mathbb Q}_z)\rangle
=\langle\beta_3({\mathbb Q}_z),\,\beta_4({\mathbb Q}_z)\rangle
&=2\left(-\,q_0q_7 + q_1q_6 + q_2q_5 + q_3q_4\right) \notag \\
&= 0\,.
\end{align}
Thus, at least for these 4 pairs in inner products, orthogonality and linear independence of the tangent basis (\ref{E17}) are dictated by the normalization criterion (\ref{E3}). More importantly, all 28 unordered pairs of the multivectors displayed in (\ref{E18}) through (\ref{E25}) turn out to be mutually orthogonal. This is hardly surprising, since we have already established the scalar-valued composition law (\ref{ed77}) for the algebra ${\cal K}^{\lambda}$, which can be used to verify further the linear independence of the tangent basis (\ref{E17}). The basis multivectors (\ref{E19}) through (\ref{E25}) thus form a well defined tangent space to $S^7$ at the tip of the multivector ${\mathbb Q}_z\in {\cal K}^{\lambda}$. Since this multivector was chosen arbitrarily and variations of its components shown in (\ref{E18}) sweep all points of $S^7$, and since ${\cal K}^{\lambda}$ and ${S^7\hookrightarrow {\cal K}^{\lambda}}$ are closed under multiplication, this procedure of finding orthogonal tangent basis at a given point in ${S^7}$ can be repeated {\it ad infinitum} to construct a continuous
field of absolutely parallel tangent basis at every point in ${S^7}$. That is, given an orthogonal tangent basis such as (\ref{E17}) at the tip of a multivector ${\mathbb Q}_z\in {\cal K}^{\lambda}$, the tangent basis at the tip of another multivector ${\mathbb P}_z\in {\cal K}^{\lambda}$ can be found by computing 
\begin{align}
&\left\{\,\beta_1({\mathbb Q}_z)\,{\mathbb P}_z,\;\beta_2({\mathbb Q}_z)\,{\mathbb P}_z,\;\beta_3({\mathbb Q}_z)\,{\mathbb P}_z,\;\beta_4({\mathbb Q}_z)\,{\mathbb P}_z,\;\beta_5({\mathbb Q}_z)\,{\mathbb P}_z,\;\beta_6({\mathbb Q}_z)\,{\mathbb P}_z,\;\beta_7({\mathbb Q}_z)\,{\mathbb P}_z\right\} \notag \\
&\qquad\qquad\qquad\qquad\qquad\qquad\qquad\qquad\qquad\quad = \left\{\,\beta_1({\mathbb P}_z),\;\beta_2({\mathbb P}_z),\;\beta_3({\mathbb P}_z),\;\beta_4({\mathbb P}_z),\;\beta_5({\mathbb P}_z),\;\beta_6({\mathbb P}_z),\;\beta_7({\mathbb P}_z)\right\}, \label{E28}
\end{align}
and so on for all points in ${S^7}$. This generates a continuous, orthogonality-preserving translation of the tangent basis at the tip of ${\mathbb Q}_z$ to the tangent basis at the tip of ${\mathbb P}_z$, for all pairs $\{{\mathbb Q}_z,\,{\mathbb P}_z\}$ in ${\cal K}^{\lambda}$. In other words, for every multivector ${\mathbb P}_z$, the set (\ref{E28}) forms a basis of ${T_p\,S^7}$, and thus each point in ${S^7}$ is characterized by a tangent multivector in the basis (\ref{E28}), representing a smooth flowing motion of that point, devoid of discontinuities, singularities, or fixed points.

So far, we have focused on the orthogonality and linear independence of the basis for ${T_q\,S^7}$ shown in (\ref{E19}) to (\ref{E25}).\break It is also straightforward to verify that the multivectors constituting this basis are normalized. Consider, for example, 
\begin{align}
\{\beta_1({\mathbb Q}_z)\}\{\beta_1({\mathbb Q}_z)\}^{\dagger}
&= \{({{\bf e}_x}\,{{\bf e}_y})\,{\mathbb Q}_z\}\{({{\bf e}_x}\,{{\bf e}_y})\,{\mathbb Q}_z\}^{\dagger} \\
&=\{({{\bf e}_x}\,{{\bf e}_y})\,{\mathbb Q}_z\}\left\{{\mathbb Q}_z^{\dagger}\,({{\bf e}_x}\,{{\bf e}_y})^{\dagger}\right\} \\
&=\,({{\bf e}_x}\,{{\bf e}_y})\left\{{\mathbb Q}_z\,{\mathbb Q}_z^{\dagger}\right\}({{\bf e}_y}\,{{\bf e}_x})\,. \label{eee25}
\end{align}
Now, as in (\ref{ddd48}) to (\ref{ddd51}), in the current notation the geometric product ${\mathbb Q}_{z}{\mathbb Q}^{\dagger}_{z}$ appearing in (\ref{eee25}) works out to be
\begin{align}
{\mathbb Q}_{z}{\mathbb Q}^{\dagger}_{z}
&=\left({\bf q}_{r}\,{\bf q}^{\dagger}_{r}+{\bf q}_{d}\,{\bf q}^{\dagger}_{d}\right)+\left({\bf q}_{r}\,{\bf q}^{\dagger}_{d}+{\bf q}_{d}\,{\bf q}^{\dagger}_{r}\right)\varepsilon \\
&=\left(\varrho^2_r + \varrho^2_d\,\right)+ 2\left(-\,q_0q_7 + q_1q_6 + q_2q_5 + q_3q_4\right)\,\varepsilon \label{eee27} \\
&=\,\text{(a scalar)}+\text{(a scalar)}\;\varepsilon\,.
\end{align}
But since $\varepsilon$ commutes with all element of the algebra ${\cal K}^{\lambda}$, substituting the expansion (\ref{eee27}) into (\ref{eee25}) reduces~(\ref{eee25})~to
\begin{align}
\{\beta_1({\mathbb Q}_z)\}\,\{\beta_1({\mathbb Q}_z)\}^{\dagger}
&= \left\{{\mathbb Q}_z\,{\mathbb Q}_z^{\dagger}\right\}({{\bf e}_x}\,{{\bf e}_y})\,({{\bf e}_y}\,{{\bf e}_x}) \\
&= \left\{{\mathbb Q}_z\,{\mathbb Q}_z^{\dagger}\right\}\,{{\bf e}_x}\,({{\bf e}_y}{{\bf e}_y})\,{{\bf e}_x} \\
&=\,{\mathbb Q}_z\,{\mathbb Q}_z^{\dagger}\,.
\end{align}
Consequently, the scalar-valued norm of the multivector $\beta_1({\mathbb Q}_z)$ in the tangent basis can be easily computed to give
\begin{align}
||\beta_1({\mathbb Q}_z)|| &=\sqrt{\Big(\{\beta_1({\mathbb Q}_z)\}\,\{\beta_1({\mathbb Q}_z)\}^{\dagger}\Big)\Big|_{\,\varepsilon\cdot\big(\{\beta_1({\mathbb Q}_z)\}\,\{\beta_1({\mathbb Q}_z)\}^{\dagger}\big)\,=\;0}\,} \\
&=\sqrt{\left({\mathbb Q}_{z}{\mathbb Q}^{\dagger}_{z}\right)\!\Big|_{\,\varepsilon\cdot\left({\mathbb Q}_{z}{\mathbb Q}^{\dagger}_{z}\right)\,=\;0}\,} \\
&=\sqrt{\left(\varrho^2_r + \varrho^2_d\,\right)\,+\,0\;} \label{E33} \\
&=\,\varrho \\
&=\,1\,.
\end{align}
Here, $\varrho=1$ because in (\ref{E19}) we started out with $||{\mathbb Q}_z||=||{\mathbb Q}_0||=||(1,\,0,\,0,\,0,\,0,\,0,\,0,\,0)||=1$, and the normalization criterion (\ref{E4}) is used for (\ref{E33}). Similar verification goes through for all basis elements of the tangent space, including $\beta_7({\mathbb Q}_z)=(I_3\,{\bf e}_{\infty})\,{\mathbb Q}_z$. Thus, the tangent space basis (\ref{E17}) expressed in (\ref{E19}) through (\ref{E25}) form orthonormal basis. 

One consequence of the above construction of the orthonormal basis for the
tangent space at every point ${{p}\in S^7}$ is that it renders the Riemann curvature tensor of ${S^7}$ to vanish with respect to asymmetric Weitzenb\"ock connection~$\Omega_{\gamma\,\beta}^{\alpha}$:
\begin{equation}
{R}^{\,\alpha}_{\;\;\,\beta\,\gamma\,\delta}\,=\,\partial_\gamma\,\Omega_{\delta\,\beta}^{\alpha}\,-\,
\partial_\delta\,\Omega_{\gamma\,\beta}^{\alpha}\,+\,\Omega_{\delta\,\beta}^{\rho}\,\Omega_{\gamma\,\rho}^{\alpha}\,-\,
\Omega_{\gamma\,\beta}^{\rho}\,\Omega_{\delta\,\rho}^{\alpha}\,=\,0\,.
\end{equation}
The resulting teleparallel geometry of $S^7$ is then entirely characterized by a totally anti-symmetric torsion tensor. This vanishing of the curvature tensor renders the parallelism on ${S^7}$ {\it absolute} --- {\it i.e.}, it guarantees the path-independence of parallel transport within ${S^7}$. A parallel transport of any arbitrary multivector in ${T_p\,S^7}$ is then simply a translation of that multivector in ${S^7}$, and this translation is path-independent because of the vanishing of the curvature tensor.

\section{Absence of Non-zero Zero Divisors in Normed Algebra ${\cal K}^{\lambda}$ with Scalar-valued Norms}

In this appendix, we prove that the apparent non-zero zero divisors in the algebra ${\cal K}^{\lambda}$ reduce to zero divisors when ${\cal K}^{\lambda}$ is normed. To that end, recall from Appendix~\ref{D} that such divisors in the algebra ${\cal K}^{\lambda}$ are exclusively of the form 
\begin{equation}
{\bf W} = \alpha\,(1+\,\varepsilon)\;\;\;\;\text{and}\;\;\;\;{\bf Z}=\gamma\,(1-\,\varepsilon)\,, \label{ff1}
\end{equation}
where $\alpha$ and $\gamma$ are arbitrary real numbers, and $\varepsilon\equiv-I_3\,{\bf e}_{\infty}$ is a pseudoscalar with the properties $\varepsilon^{\dagger}=\varepsilon$ and $\varepsilon^2=1$. Now, by definition, a ``normed algebra'' is an algebra equipped with a definition of norm that yields a scalar-valued norm for any given element of the algebra, such as that defined in (\ref{E3}) for the algebra ${\cal K}^{\lambda}$ under consideration. On the other hand, in the notation of Appendix~\ref{E}, a general multivector ${\mathbb Q}_z$ in the graded basis of algebra ${\cal K}^{\lambda}$ is expressed~as
\begin{equation}
{\mathbb Q}_z=\,q_0+q_1\,{\bf e}_x{\bf e}_y+q_2\,{\bf e}_z{\bf e}_x+q_3\,{\bf e}_y{\bf e}_z+q_4\,{\bf e}_x{\bf e}_{\infty}+q_5\,{\bf e}_y{\bf e}_{\infty}+q_6\,{\bf e}_z{\bf e}_{\infty}+q_7\,I_3{\bf e}_{\infty}\,, \label{ff2}
\end{equation}
so that, using $\varepsilon\equiv-I_3\,{\bf e}_{\infty}$, the geometric product ${\mathbb Q}_{z}{\mathbb Q}^{\dagger}_{z}$ takes the following simple form, as in (\ref{21}), (\ref{sd47}), and (\ref{eee27}):
\begin{equation}
{\mathbb Q}_{z}{\mathbb Q}^{\dagger}_{z}
=(\text{a scalar})\,+\, 2\left(-\,q_0q_7 + q_1q_6 + q_2q_5 + q_3q_4\right)\,\varepsilon\,. \label{ff3}
\end{equation}
But, comparing (\ref{ff1}) with (\ref{ff2}), we recognize that, in the notation of (\ref{ff2}), the non-zero zero divisors are of the form:
\begin{equation}
{\bf U} = q^{}_0\,-\,q^{}_7\,\varepsilon\,,\;\;\;\;\text{with}\;\;\;\;q^{}_7=\mp\,q^{}_0\;\;\;\;\text{and}\;\;\;\;q^{}_1=q^{}_2=q^{}_3=q^{}_4=q^{}_5=q^{}_6\equiv0\,,\label{ff4}
\end{equation}
so that their geometric product
\begin{equation}
{\bf U}\,{\bf U}^{\dagger}=\,(q^{}_0\,-\,q^{}_7\,\varepsilon)\,(q^{}_0\,-\,q^{}_7\,\varepsilon)^{\dagger} = \,\left(q_0^2+q_7^2\right) - 2\,\left(q^{}_0\,q^{}_7\right)\,\varepsilon \,. \label{ff5}
\end{equation}
On the other hand, as discussed in Appendix~\ref{D}, the necessary and sufficient condition for scalar-valued norm $||{\mathbb Q}_z||$ is
\begin{equation}
\varepsilon\cdot\left({\mathbb Q}_{z}{\mathbb Q}^{\dagger}_{z}\right)
= \,2\left(-\,q_0q_7 + q_1q_6 + q_2q_5 + q_3q_4\right) = 0\,.
\end{equation}
In light of (\ref{ff4}) and (\ref{ff5}), this condition gives $q^{}_7=\mp{q^{}_0=0}$ as necessary and sufficient conditions for the scalar-valued\break norm of ${\bf U}$, which, using (\ref{ff5}), leads to $||{\bf U}||=0$. The positive definiteness of the norm defined in (\ref{E3}) then necessitates that ${\bf U}=\,0$, reducing it to a harmless zero divisor. Thus, there are no non-zero zero divisors in the normed algebra~${\cal K}^{\lambda}$.

\vspace{-0.6cm}

\section*{Acknowledgments}

\vspace{-0.4cm}

I thank Tevian Dray for his comments on the previous version of this preprint, which led to the proof in Appendix~\ref{B}.


\vspace{-0.4cm}

\begin{thebibliography}{}
\end{thebibliography}

\vspace{-0.4cm}

\begin{thebibliography}{}

\vspace{-0.6cm}

\bibitem[Doran(2003)]{Clifford}C. Doran and A. Lasenby, {\it Geometric Algebra for Physicists} (Cambridge University Press, Cambridge, 2003).

\bibitem[Dorst(2007)]{Dorst}L. Dorst, D. Fontijne, and S. Mann, {\it Geometric Algebra for Computer Science} (Elsevier, Amsterdam, 2007).

\bibitem[Christian(2017)]{RSOS}J. Christian, {\it Quantum correlations are weaved by the spinors of the Euclidean primitives}, R. Soc. Open Sci., {\bf 5}, 180526 (2018); \href{https://doi.org/10.1098/rsos.180526}{https://doi.org/10.1098/rsos.180526}; See also \href{https://arxiv.org/abs/1806.02392}{https://arxiv.org/abs/1806.02392} (2018).

\bibitem[Kenwright(2012)]{Kenwright}B. Kenwright, {\it A beginners guide to dual-quaternions: what they are, how they work, and how to use them for 3D character hierarchies}, in Proceedings of the 20th International Conferences on Computer Graphics, Visualization and Computer Vision, 1--10 (2012).

\bibitem[Hurwitz(1898)]{Hurwitz-1898}A. Hurwitz, {\it \"Uber die Composition der quadratischen Formen von beliebig vielen Variabeln}, Nachr. Ges. Wiss. G\"ottingen, {\bf 1898}, 309--316 (1898).


\bibitem[Hurwitz(1923)]{Hurwitz-1923}A. Hurwitz, {\it \"Uber die Komposition der quadratischen Formen}, Math. Ann., {\bf 88} (1–2), 1--25 (1923).

\bibitem[Baez(2002)]{Baez}J. C. Baez, {\it The octonions}, Bull. Am. Math. Soc., {\bf 39}, 145--205 (2002).

\bibitem[Milnor(1956)]{Milnor}J. W. Milnor, {\it Topology from the Differentiable Viewpoint} (Princeton University Press, Princeton, New Jersey, 1997).

\bibitem[Lasenby(2004)]{Lasenby-1} A. Lasenby, {\it Recent applications of conformal geometric algebra}, in Computer Algebra and Geometric Algebra with Applications, 298--328 (Springer, New York, 2004).

\bibitem[Lasenby(2011)]{Lasenby-2}A. Lasenby, {\it Rigid body dynamics in a constant curvature space and the ‘1D up’ approach to conformal geometric algebra}, in Guide to Geometric Algebra in Practice, 371--389 (Springer, New York, 2011).

\bibitem[Jordan(1933)]{Jordan}P. Jordan,
{\it \"Uber die Multiplikation quantenmechanischer Gr\"o\ss en}, Zeitschrift f\"ur Physik, {\bf 80}, 285--291 (1933).

\bibitem[Dirac(1939)]{Dirac}P. A. M. Dirac, {\it The relation between mathematics and physics}, Proceedings of the Royal Society (Edinburgh), {\bf 59}, Part II, 122--129 (1939).

\bibitem[Lounesto(2001)]{Lounesto}P. Lounesto, {\it Octonions and triality}, Advances in Applied Clifford Algebras, {\bf 11}, 191 (2001).

\bibitem[Forbenius(1878)]{Forbenius}F. G. Frobenius, {\it \"Uber lineare Substitutionen und bilineare Formen}, Journal für die reine und angewandte Mathematik, {\bf 84}, 1--63 (1878).

\bibitem[Dray(1999)]{Dray}T. Dray and C. A. Manogue, {\it The Geometry of the Octonions} (World Scientific, Singapore, 2015), Chapter 5.

\bibitem{RSOS-2}J. Christian, {\it Response to ‘Comment on “Quantum correlations are weaved by the spinors of the Euclidean primitives”’}, R. Soc. Open Sci., {\bf 9}, 220147 (2022); \href{https://doi.org/10.1098/rsos.220147}{https://doi.org/10.1098/rsos.220147}. See also \href{https://arxiv.org/abs/2211.09867}{https://arxiv.org/abs/2211.09867}.

\bibitem{Local}J. Christian, {\it Local origins of quantum correlations rooted in geometric algebra}, https://doi.org/10.48550/arXiv.2205.11372.

\end{thebibliography}
\end{document}